\documentclass[reqno, 12pt]{amsart}

\usepackage{enumerate}
\usepackage{amssymb}
\usepackage[nocompress]{cite}

\usepackage{xcolor}
\usepackage{hyperref}
\hypersetup{
    colorlinks,
    linkcolor={red!50!black},
    citecolor={blue!50!black},
    urlcolor={blue!80!black},
    linktoc=page
}

\usepackage{tikz-cd}
\usepackage[margin=3.14cm]{geometry}
\usepackage{mathtools}
\usepackage{bm, bbm}
\usepackage{textcomp}
\usepackage{stmaryrd}
\usepackage{tcolorbox}
\usepackage[normalem]{ulem}

\usepackage{aliascnt}
\usepackage[capitalize]{cleveref}

\newtheorem{lemma}{Lemma}[section]

\newaliascnt{theorem}{lemma}
\newtheorem{theorem}[theorem]{Theorem}
\aliascntresetthe{theorem}

\newtheorem{maintheorem}{Theorem}

\newtheorem*{theorem*}{Theorem}

\newcounter{thmAast}
\newtheorem{theoremAast}[thmAast]{Theorem}
\crefname{thmAast}{Theorem}{Theorems}

\newaliascnt{claim}{lemma}

\aliascntresetthe{claim}

\newaliascnt{corollary}{lemma}
\newtheorem{corollary}[corollary]{Corollary}
\aliascntresetthe{corollary}

\newtheorem*{corollary*}{Corollary}
\newtheorem*{restatedbeautiful}
  {Theorem~\ref{thm:beautiful_main} (restated)}

\newaliascnt{proposition}{lemma}
\newtheorem{proposition}[proposition]{Proposition}
\aliascntresetthe{proposition}

\newtheorem*{proposition*}{Proposition}
\newtheorem*{JIC}{Joint Intersectivity Conjecture}

\newtheorem{problem}{Problem}
\newtheorem*{problem*}{Problem}
\newtheorem*{question*}{Question}
\theoremstyle{definition}

\newtheorem*{claim*}{Claim}

\newaliascnt{definition}{lemma}
\newtheorem{definition}[definition]{Definition}
\aliascntresetthe{definition}

\newaliascnt{example}{lemma}
\newtheorem{example}[example]{Example}
\aliascntresetthe{example}

\newtheorem*{remark}{Remark}

\newcommand		{\p}[1]	{\left(#1\right)}

\newcommand		{\abs}[1]{\left|#1\right|}

\newcommand 	{\norm}[1]{\left\lVert#1\right\rVert}
\newcommand 	{\nnorm}[1]{\lvert\!|\!| #1|\!|\!\rvert}
\newcommand 	{\short}{\hspace{-0.1em}}
\newcommand		{\snorm}[1]{\left|\short\left|\short\left|#1\right|\short\right|\short\right|}
\newcommand    	{\R} {\mathbb R}
\newcommand    	{\C} {\mathbb C}
\newcommand    	{\N} {\mathbb N}
\DeclarePairedDelimiterX{\inp}[2]{\langle}{\rangle}{#1, #2}
\DeclareMathOperator{\E}{\mathbb E}

\newcommand    	{\Q} {\mathbb Q}
\newcommand    	{\Z} {\mathbb Z}

\newcommand    	{\ve} {\varepsilon}

\renewcommand    	{\d} {\, \mathrm d}
\newcommand    	{\1} {\mathbbm 1}
\newcommand    	{\T} {\mathbb T}

\renewcommand 	{\emph}[1]{\textit{#1}}

\newcommand 	{\con}[1]{\overline{#1}}

\renewcommand   {\mod} {\mathrm{mod}\ }

\newcommand    	{\m} {\mathrm m}

\newcommand 	{\core} {\mathrm{core}}

\newcommand 	{\acts} {\curvearrowright}

\newcommand 	{\nilbor} {\text{Nil-Bohr$_0$}}

\newcommand 	{\gen}[1]{\langle#1\rangle}

\newcommand{\CA}{{\mathcal A}}
\newcommand{\CB}{{\mathcal B}}
\newcommand{\CC}{{\mathcal C}}
\newcommand{\CD}{{\mathcal D}}

\newcommand{\CF}{{\mathcal F}}
\newcommand{\CI}{{\mathcal I}}

\newcommand{\CM}{{\mathcal M}}
\newcommand{\CN}{{\mathcal N}}

\newcommand{\CX}{{\mathcal X}}
\newcommand{\CY}{{\mathcal Y}}
\newcommand{\CW}{{\mathcal W}}
\newcommand{\CZ}{{\mathcal Z}}
\newcommand{\CG}{{\mathcal G}}

\newcommand{\CU}{{\mathcal U}}

\begin{document}
\title[Pairwise Independent Intersective Polynomials]{Multidimensional Szemer\'edi
Theorem for Pairwise Independent Jointly Intersective Polynomials}
\author{Nikos Frantzikinakis}
\author{Davide Sclosa}
\address[Nikos Frantzikinakis]{University of Crete, Department of mathematics and applied mathematics, Voutes University Campus, Heraklion 70013, Greece} \email{frantzikinakis@gmail.com}
\address[Davide Sclosa]{King's College London, Department of mathematics,
Strand Campus, London, WC2R 2LS, United Kingdom
} \email{davide.sclosa@gmail.com}

\begin{abstract}
We develop a new approach to the Joint Intersectivity Conjecture of Bergelson, Leibman, and Lesigne for commuting transformations. The conjecture asserts that the necessary divisibility condition is the only obstruction to polynomial multiple recurrence and, consequently,  to the multidimensional polynomial Szemer\'edi theorem. We establish the conjecture
for pairwise independent polynomials, obtaining the first general result
for sequences of arbitrary length beyond the classical one-transformation setting. Our argument is inspired by the Furstenberg-Katznelson approach and its polynomial enhancement due to Bergelson-Leibman. Its main novelty, however, lies in the development of new techniques that allow us to bypass the use of the polynomial van der Waerden theorem, on which these methods traditionally rely, and for which no suitable analogue is currently available in our setting.
In particular, we prove a Nil-Bohr returns theorem for compact extensions in systems with nilpotent structure, as well as a nil-plus-null decomposition for polynomial multicorrelation sequences that preserves joint intersectivity and positivity.
\end{abstract}

\thanks{
The authors were supported by the Research Grant ELIDEK HFRI-NextGenerationEU-15689. The second author was also partially supported by the Research Grant EPSRC UKRI1021.
}




	
	\maketitle
	
	\setcounter{tocdepth}{1}
	\tableofcontents
	
    
	\section{Introduction and main results}

    Furstenberg's ergodic proof of Szemer\'edi's theorem~\cite{furstenberg1977ergodic}, together with the multidimensional extension of Furstenberg and Katznelson~\cite{furstenberg1978ergodic}, revealed a deep and fruitful connection between additive combinatorics and ergodic theory. A striking further development was obtained by Bergelson and Leibman~\cite{bergelson1996polynomial}, who proved a polynomial multidimensional extension: if $\Lambda \subseteq \mathbb{Z}^d$ has positive upper density and
$p_1,\ldots,p_\ell$ are \emph{integral polynomials}
(i.e.~rational polynomials taking integer values on the integers)
with zero constant term, then for every $v_1,\ldots,v_\ell \in \mathbb{Z}^d$ there exist $m\in \mathbb{Z}^d$ and~$n\in \mathbb{N}$ such that
$$
m,\ m+p_1(n)v_1,\ \ldots,\ m+p_\ell(n)v_\ell \in \Lambda.
$$
Their argument proceeds through an ergodic multiple recurrence theorem, and despite substantial subsequent efforts, no proof of the general result that avoids ergodic theory is currently known.

A natural next question is whether one can obtain the full polynomial Szemerédi phenomenon without assuming that the polynomials vanish at zero. At this point, however, genuine arithmetic obstructions appear. For instance, if $\Lambda=r\mathbb{Z}$, then any recurrence statement of this type forces the values $p_1(n),\ldots,p_\ell(n)$ to be simultaneously divisible by $r$. This leads to the notion introduced in~\cite{bergelson2008intersective}: a family
of integral polynomials $p_1,\ldots,p_\ell$ is called \emph{jointly intersective} if for every $r \in \mathbb{N}$ there exists $n \in \mathbb{N}$ such that $p_1(n),\ldots,p_\ell(n)$ are all divisible by $r$.
Bergelson, Leibman, and Lesigne conjectured in~\cite{bergelson2008intersective} that this
obvious obstruction is the only one.

\begin{JIC} \label{Con:BLL}
Let $p_1,\ldots,p_\ell$ be jointly intersective and $T_1,\ldots,T_\ell$ be commuting measure-preserving transformations on $(X,\mathcal{X},\mu)$. Then for every $A \in \mathcal{X}$ with $\mu(A)>0$, one has
$$
\liminf_{N\to\infty}\frac{1}{N}\sum_{n=1}^N
\mu\bigl(A \cap T_1^{-p_1(n)}A \cap \cdots \cap T_\ell^{-p_\ell(n)}A\bigr)>0.
$$
\end{JIC}

In other words, this conjecture seeks to completely characterize the polynomial families that are good for multiple recurrence.
By Walsh's norm convergence theorem~\cite{walsh2012norm},
the limit inferior above is in fact a limit.

We remark that the conjecture can be stated, and all our subsequent results can be proved, with the uniform averages
$\liminf_{N-M\to\infty}\frac{1}{N-M}\sum_{n=M+1}^N$
in place of Ces\`aro averages, without any essential changes to the proofs.

The conjecture remains widely open. For general $\ell$, a positive answer is known when $T_1,\ldots,T_\ell$ are powers of a single transformation (see \cite[Theorem~F]{frantzikinakis2008multiple} for $\ell=2,3$ and \cite{bergelson2008intersective} for general $\ell$) and the proof relies on the Host-Kra theory of characteristic factors~\cite{host2005nonconventional}, which is not available in the setting of several commuting transformations. Moreover, the Bergelson-Leibman argument for the polynomial Szemerédi theorem with zero constant term~\cite{bergelson1996polynomial} does not seem adaptable to this setting, since it depends on a coloristic input that is still unknown for intersective polynomials. Very recently, the conjecture was established in the case $\ell=2$ in~\cite{frantzikinakis2025ergodic}, by combining a broad range of ideas developed over the past two decades in ergodic theory.
For $\ell \geq 3$, however, allowing nontrivial linear relations among pairwise independent polynomials introduces substantial new difficulties, and this general case remained beyond the reach of existing methods.

The main goal of this article is to introduce a new approach to the Joint Intersectivity Conjecture and to prove the first general result in this direction for polynomial families of arbitrary length outside the classical one-transformation setting.

We say that the polynomials~$p_1,\ldots,p_\ell$ are
\emph{pairwise independent} if~$1,p_i,p_j$ are linearly
independent whenever~$i\neq j$.
When~$\ell=1$, we require only that $p_1$ be non-constant. We emphasize that this condition is much weaker than linear independence (under which the result is already known; see the remark following~\cite[Corollary~2.11]{frantzikinakis2025joint}), as the polynomials we consider may satisfy many linear relations.

Our main result is the following.

\begin{maintheorem}\label{thm:main_}
The Joint Intersectivity Conjecture has a positive answer
for jointly intersective families
of pairwise independent polynomials.
\end{maintheorem}

The core new ingredient in the proof of Theorem A is a verification of the Joint Intersectivity Conjecture for a class of structured systems called $\CZ_s$-systems (defined in \cref{sec:prel_HK}). To prove this structured recurrence theorem, we develop new tools that make it possible to implement the Furstenberg–Katznelson multiple recurrence strategy without relying on the coloring results that underlie the classical polynomial Szemer\'edi theorem. We expect this mechanism to have further applications to multiple recurrence for commuting transformations. See \cref{SS:Beyond} for more details. 

A family of polynomials is jointly intersective if and only if all of its members are divisible by a common intersective polynomial~\cite[Proposition~6.1]{bergelson2008intersective}. While the most familiar examples of intersective polynomials $p$ with $p(0)\neq 0$ are products of linear factors, such as $p(n)=n^2-1$,
there also exist intersective polynomials with no linear factor at all,
such as~$p(n)=(n^2-2)(n^2-17)(n^2-34)$.
See \cite{bergelson2007van} for further examples.

\cref{thm:main_}
yields many genuinely new instances that were inaccessible to previous methods.
A representative new case is the following.
Let $\ell=3$, let $p$ be any non-constant intersective polynomial,
let $q,r \in \mathbb{Z}[t]$
be such that $q/r$ is non-constant,
and consider the polynomials
$$
p_1=pq,\quad p_2=pr,\quad p_3=p(q+r).
$$
Then $p_1,p_2,p_3$ are jointly intersective
and pairwise independent,
and hence fall within the scope
of our main result.
The conjecture remains open even in the degenerate case
of three equal intersective polynomials.
On the other hand, it now seems likely that a solution
of the equal-polynomial case for an arbitrary number of commuting transformations
could be combined with the techniques developed here
to resolve the full Joint Intersectivity Conjecture.

Combining \cref{thm:main_} with Furstenberg's correspondence principle~\cite{furstenberg1978ergodic}
in its multidimensional form,
we obtain the following combinatorial result.

\begin{corollary*}
If $\Lambda \subseteq \mathbb{Z}^d$ has positive upper density and $p_1,\ldots,p_\ell$ are jointly intersective pairwise independent polynomials,
then for every $v_1,\ldots,v_\ell \in \mathbb{Z}^d$ there exist
an integer~$n\in \N$ and a vector~$m\in \Z^d$
such that~$m,\ m+p_1(n)v_1,\ \ldots,\ m+p_\ell(n)v_\ell \in \Lambda.$
\end{corollary*}

\subsection*{Notation}
We collect here the notation used most often in the paper;
most of the following objects are defined in the preliminaries, \cref{sec:prel}.
\begin{itemize}
\item 
$\N=\{1,2,3,\ldots\}$.
Integers~$\Z$: $m,n$, sometimes~$\alpha,\beta$. Rationals~$\Q$.
For an integer~$N\geq 1$, $[N]:=\{1,\ldots,N\}$ and~$\E_{n\in [N]} := \frac{1}{N}\sum_{n=1}^N$.

\item 
Functions: $f,g,h$. Indicator functions:~$\1_A,\1_B$.

\item  We write $[\![t]\!]:=\lfloor t+1/2\rfloor$ for the closest integer to $t$.

\item 
Measures: $\mu,\nu,\rho$. Disintegrations: $(\mu_y)_{y\in Y}$.
Probability spaces: $X,Y,Z,\Omega,\Lambda$.
Expectation and conditional expectation: $\E$.

\item 
Measure-preserving systems: $(X,\CX,\mu,T)$, $(X,T)$, or simply~$X$.
Factors: $Y$, $Z$, $W$, with sub-$\sigma$-algebras~$\CY,\CZ,\CW$ respectively;
factor maps: $\pi,\kappa$.
Invariant factor: $I(X)$, with sub-$\sigma$-algebra~$\CI(X)$.
Order-$s$ Host--Kra factor: $Z_s(X)$, with sub-$\sigma$-algebra~$\CZ_s(X)$
and seminorm~$\snorm{\cdot}$.

\item 
Polynomials: $p,q,r,p_i,q_i$.
All polynomials in this paper are integral:
they have rational coefficients and take integer values on the integers.

\item 
Nilmanifolds: $G/\Gamma$.
Subgroups of~$G$: $H,L,M$, or~$\Gamma,\Lambda$ if discrete; elements: $g,h,\ell,a,b$.
Quotient, index and generated subgroup: $G/H$, $[G\colon H]$ and~$\gen{G,H}$.
Points of~$G/\Gamma$: $g\cdot 1_{G/\Gamma}$ or~$g \Gamma$.
Nilsequences: $\CN,\CM$. Dual sequences: $\CD$.
\end{itemize}

\subsection*{Acknowledgement}
The first author thanks Redmond McNamara for helpful discussions during the early stages of this project in 2024.

\smallskip
\noindent\textbf{Use of AI.}
During the final stages of preparing this article, the authors used
AI-assisted tools for proofreading and language editing.
The mathematical content of this manuscript is entirely the work of
the authors, who take full responsibility for its correctness.

    
\section{Overview of the proof}
The goal of this section is to provide a roadmap for the proof of our main result. For technical reasons, we need a weighted strengthening of \cref{thm:main_}
that fits the inductive structure of our argument, with weights~$\CN$
drawn from the class of measurable varying positive nilsequences,
whose definition is given in \cref{sec:prel_nil_definitions}.

\begin{theoremAast} \label{thm:main}
Let~$(X,\mu,T_1,\ldots,T_\ell)$ be a commuting system and~$A$ a set of positive measure.
Let~$q_1,\ldots,q_m,p_1,\ldots,p_\ell$ be jointly intersective polynomials
with~$p_1,\ldots,p_\ell$ pairwise independent.
Let~$\CN\colon X\times \Z^m \to [0,1]$ be a measurable varying positive nilsequence.
Then
\[
\liminf_{N\to\infty} \frac{1}{N}\sum_{n=1}^N
\int_X
\CN(x,q_1(n),\ldots,q_m(n)) \cdot \1_A(x) \cdot \1_A( T_1^{p_1(n)} x) \cdots \1_A(T_\ell^{p_\ell(n)} x)
\d\mu(x) > 0.
\]
\end{theoremAast}

Clearly \cref{thm:main} implies \cref{thm:main_}.
We now present the strategy to prove \cref{thm:main} and isolate the two novel ingredients,
which we believe may find further applications in multiple recurrence.
We also briefly discuss the main obstacles that had to be overcome
and several approaches that did not work.

\subsection{A first reduction} As noted earlier, the original method of Bergelson and Leibman used to prove the polynomial Szemer\'edi theorem cannot be applied in our setting, since it relies on a coloristic result that is not available here. Our starting point is therefore the following weighted version
of a result by the first author and Kuca~\cite[Theorem~2.8]{frantzikinakis2025joint}, which allows us to make an important reduction that introduces additional structure into the problem.

	\begin{proposition} \label{cor:partial_factor_use}
	Let~$(X,\mu,T_1,\ldots,T_\ell)$, $\CN$, $q_1,\ldots,q_m$, $p_1,\ldots,p_\ell$ be as in \cref{thm:main}.
    Let~$d:=\max_{1\leq i\leq \ell} \deg p_i$ and let~$r$
    denote
    the step of~$\CN(x,q_1(n),\ldots,q_m(n))$.
	Then there is~$s$, dependent only on~$d, \ell, r$, such that
	\[
	\lim_{N\to\infty} \frac{1}{N}\sum_{n=1}^N
	\int \CN(x,q_1(n),\ldots,q_m(n))
		\cdot f_0(x) \cdot f_1(T_1^{p_1(n)} x) \cdots f_\ell(T_\ell^{p_\ell(n)} x)	\d\mu(x)
		= 0
	\]
	whenever~$f_0,f_1,\ldots,f_\ell \in L^\infty(\mu)$
    satisfy~$f_i\perp \CZ_s(T_i)$ for some~$i=1,\ldots,\ell$. 
\end{proposition}

It follows that in order to prove \cref{thm:main}, we may replace $\1_A$ by $\E(\1_A \mid \CZ_s(T_i))$ in the corresponding places. Hence, \cref{thm:main} is reduced to the following statement.
 
\begin{proposition}\label{prop:mainZ}
	\cref{thm:main} holds if we 
	replace  $\1_A$ with 
	$\E(\1_A|\CZ_s(T_i))$ in the respective places. Moreover, under this hypothesis, the conclusion holds for all families of jointly intersective polynomials (not necessarily pairwise independent). 
\end{proposition}
Replacing the indicator functions by their conditional expectations is a crucial simplification, since it effectively reduces the proof of \cref{thm:main} to the case in which
$
\mathcal{Z}_s(T_1)=\cdots=\mathcal{Z}_s(T_\ell)=\mathcal{X}.
$
In this setting, one can hope to exploit the Host-Kra
structure theorem \cite[Chapter 16, Theorem 1]{host2018nilpotent}, which asserts that the ergodic
components of $\CZ_s$-systems are inverse limits of nilsystems.
Unfortunately, the resulting multiple recurrence problem remains formidable.
One reason is that we cannot  reduce
matters to the case where all systems are ergodic,
while a direct attack on non-ergodic $\CZ_s$-systems
appears to be out of reach.
In fact, without making further use of the commutativity
of the transformations, we do not know how to prove
\cref{prop:mainZ}, or  whether it is true, even in
some of its simplest special cases, such as the following one.
\begin{problem}
	Let $T,S$ be measure-preserving transformations on a probability space $(X,\CX,\mu)$ and suppose that for some $s\in \N$ we have $\CZ_s(T)=\CZ_s(S)=\CX$. Is it then true that for every $A\in \CX$ with $\mu(A)>0$ there exists $n\in \N$ such that 
	$$
	\mu(A\cap T^{-n}A\cap S^{-n}A)>0 \, ?
	$$ 
	\end{problem}

The next natural question is how one can once again exploit
commutativity.
One approach, used in \cite{frantzikinakis2025ergodic}, was to prove a limit formula
reducing the recurrence problem to one involving polynomials
with zero constant term. More precisely, in \cite{frantzikinakis2025ergodic},
averages along $p(k!n+n_k)$ are compared, for suitable $n_k\in\N$,
with averages along $k!n$, as $k\to\infty$.
However, we were not able to make this approach work
in the present setting.
Furthermore, when all the transformations are ergodic, commutativity together with the results of \cite{frantzikinakis2025ergodic} implies that they are simultaneously isomorphic to ergodic pro-nilsystems. This is a major simplification, and in this case the main result can be proved using the approach of \cite{bergelson2008intersective}. Unfortunately, no useful variant of the \cite{frantzikinakis2025ergodic} result seems to be available in the non-ergodic setting.

In the next subsections we explain an approach that does
work.

\subsection{The main ingredients needed to prove \cref{prop:mainZ}} 
The initial idea is to return to the original method of Furstenberg-Katznelson, later used by Bergelson-Leibman, and to establish multiple recurrence results via the Furstenberg-Katznelson structural theorem, which decomposes the underlying commuting action into compact and weak-mixing extensions, or combinations of the two; see \cite[Theorem 2.4]{furstenberg1978ergodic}.

Using this result, it turns out that in the setting
of \cref{prop:mainZ}, the essential difficulty is to show that
the multiple recurrence property passes through compact extensions. More precisely, assuming that \cref{prop:mainZ} holds for
a common factor $\CY$ of the system, and that the
extension $\CX\to\CY$ is compact for each of
$T_1,\ldots,T_\ell$, one must prove that \cref{prop:mainZ} also
holds for $\CX$.
This turns out to be a pretty difficult  task, and it occupies most of
the present article.
Our main technical innovations are designed to address precisely this problem. They include two new results, which we believe are of independent interest and which we now describe.

\subsubsection{$\nilbor$-returns property for compact extensions of $\CZ_s$-systems}
The next result is  proved in Sections \ref{sec:iso_extension_ergodic} and \ref{sec:iso_Zs}, and is what allows us, in the setting of $\CZ_s$-systems, to upgrade the return-times structure of compact extensions to a nilstructure. This stronger information is what ultimately enables us to avoid the use of coloristic arguments, which play a crucial role in the proof of the polynomial Szemer\'edi theorem of Bergelson-Leibman.

\begin{theorem} [$\nilbor$ Returns] \label{thm:beautiful_main}
Let~$X$ be a $\CZ_s$-system and~$X\to Y$ a compact extension.
Then there is a countable subset~$\CF\subseteq  L^\infty(X)$, dense in~$L^2(X)$,
with the following property:
for every~$f\in \CF$ and~$\ve>0$
there is a measurable varying $\nilbor$ sequence~$\CM\colon Y\times \Z\to \{0,1\}$ satisfying
\[
	\CM(y,n)=1 \qquad\implies\qquad \E(\abs{T^n f - f}^2\mid Y)(y) \leq\ve
\]
for almost every $y\in Y$ and every $n\in \N$.
\end{theorem}

In other words, the return-times of~$f$ contain
 $\nilbor$ sets.

\begin{example}
	To get a better feel for the result, consider the following explicit  example. Let $T(x,y,z):=(x,y+h_1(x),z+k(x)y+h_2(x))$, where $k\colon \T\to\Z$ and $h_1,h_2\colon \T\to \T$ are measurable. Then $$ T^n(x,y,z)=(x,y+nh_1(x),z+nk(x)y+\binom{n}{2}k(x)h_1(x)+nh_2(x)), $$
	and $X:=(\T^3, m_{\T^3}, T)$ is a $\CZ_2$-system.  We also let $Y$  be the factor defined by the functions that depend only on $x,y$, and as $\CF$ the linear combinations of functions of the form $e(kx+ly+mz)$, $k,l,m\in \Z$ (or the set of continuous functions $C(\T^3)$). The extension $X\to Y$ is compact. 
	For fixed $x,y\in \T$, the set of $\ve$-return times contains the intersection of the sets $\{n\in \N\colon \norm{nk(x)y+\binom{n}{2}k(x)h_1(x)+nh_2(x)}_\T \leq \delta\}$ and $\{n\in \N\colon \norm{nh_1(x)}_\T \leq \delta \}$, for $\delta=\delta(\varepsilon)$ sufficiently small, which is a Nil-Bohr$_0$ set. If instead we took as $Y$ the factor defined by the functions that depend only on $x$, then the conclusion would no longer hold: assuming the extension compact is necessary.
\end{example}
Despite the simplicity of the previous example, the proof of \cref{thm:beautiful_main} in its full generality is surprisingly delicate and is divided into three steps, each with its own challenges.
\smallskip

$(i)$ In the ergodic nilsystem case, the task is to turn the abstract notion of a compact extension into a concrete algebraic description. The key finding is the following proposition:
in reduced form, the extension $G/\Gamma\to G/L$ is compact
if and only if the identity component of $L$ is contained
in the center of $G$ (see \cref{thm:iso_central}). Once this is established,
the $\nilbor$ control follows from the behavior of nilcharacters
along central fibers (see \cref{prop:key_lem_nilsystem}).  

\medskip 

$(ii)$ In the ergodic pro-nilsystem case, the main obstruction is that
a compact extension $X \to Y$ does not automatically induce
compact extensions on the approximating nilsystem factors.
To overcome this, we must make suitable choices of the approximants
$X_n$ that satisfy the ``good intersection property'' with respect
to $Y$ (see \cref{sec:good}). This is achieved in \cref{prop:pronil_rectified},
after establishing the crucial fact that, by slightly enlarging the nilsystem factors $X_n$,
one can ensure that the good intersection property holds.
This crucial fact is proven in~\cref{prop:Znil_Ypro}:
via a ``killing the denominators'' construction
(which ultimately reduces to a claim about rationally equivalent lattices),
we are able to guarantee the algebraic condition for good intersection
established in~\cref{lem:LM_expected}.
The idea of the proof is presented in the paragraph
right before~\cref{prop:Znil_Ypro}.

\medskip 

$(iii)$ In the general case of $\CZ_s$-systems, we plan to use part $(ii)$ to extract useful information from the ergodic components, and the main difficulty is to find a way to exploit this information at the level of the non-ergodic system. A first obstacle is measurability: each ergodic component is isomorphic
to a pronilsystem, but it is not immediately clear whether these isomorphisms can be chosen measurably across components.
We rely on recent work of Jamneshan and Machado~\cite{jamneshan2026non}, who show that the nilsystems approximating the individual ergodic components can be glued together in measurable bundles to simultaneously approximate the whole non-ergodic system.
A more substantial obstacle is that the $\nilbor$ property does not pass to arbitrary mixtures of  compact extensions of ergodic $\CZ_s$-systems, because of the potential uncountability of several choices involved. To overcome this difficulty, we establish a crucial ``countable control'' result for compact extensions of nilsystem bundles,
see \cref{ex:countable_control_I}, \cref{ex:countable_control_II}
and \cref{countable_control_simple_nilbundle}.
Countable control,
together with a non-ergodic generalization of the killing-the-denominators
construction, \cref{killing_nonergodic},
gives the full \cref{thm:beautiful_main}.

\subsubsection{Nil-plus-null decomposition that respects intersectivity}
Nil-plus-null decompositions for multiple correlation sequences
were first established by Bergelson, Host, and Kra
\cite{BHK-nilsequences}. 
Leibman in~\cite{leibman2015nilsequences}  proved that an integral combination
of nilsequences can be decomposed into the sum of a nilsequence and a nullsequence.
We will need a strengthening of this decomposition result that preserves
joint intersectivity and a certain positivity property.

	\begin{theorem} \label{thm:nil+null}
	Let~$G/\Gamma$ be a nilmanifold,
	$f \colon G/\Gamma \to [0,1]$ be a continuous function with~$f(1_{G/\Gamma})>0$, 
	$\rho$ be a Borel probability measure on~$G^\ell$ and  $p_1,\ldots,p_\ell$
	be jointly intersective polynomials.
	Then there are integers~$\alpha,m\in \N$ and~$\beta\geq 0$
	such that 
	\begin{equation} \label{eq:nilnul_goal}
		\int_{G^\ell} f(a_1^{p_1(\alpha n+\beta)}\cdots a_\ell^{p_\ell(\alpha n+\beta)} \cdot 1_{G/\Gamma}) 
		\d\rho(a_1,\ldots,a_\ell)
		\geq g({b_1}^{q_{1}(n)} \cdots {b_m}^{q_m(n)} \cdot 1_{H/\Lambda})
		+ e(n)
	\end{equation}
	for every $n\in \N$, 	where $H$ is a rational subgroup of $G$ and $H/\Lambda$ is a nilmanifold, $b_1,\ldots,b_m\in H$, 
	$g \colon H/\Lambda\to [0,1]$ is a continuous function
	with~$g(1_{H/\Lambda})>0$,
	$q_1,\ldots,q_m$ are polynomials
	such that the family $p_1(\alpha n+\beta),\ldots,p_\ell(\alpha n+\beta),q_1(n),\ldots,q_m(n)$ is jointly intersective,
	and~$e(n)$ a nullsequence with values in $[-1,1]$.
\end{theorem}

What we actually use is a generalization of \cref{thm:nil+null},
namely \cref{P:nil+nul=useful},
in which the function~$f$ is not fixed but varies measurably,
and the measure itself is also not fixed but
dependent on a parameter.


Leibman proved in~\cite{leibman2015nilsequences} a similar result with equality, provided that one allows the right-hand side to be an arbitrary nilsequence. His argument turns out to be insufficient for our purposes, as it does not control
for the additional properties we require for the right hand side.
In contrast to~\cite{leibman2015nilsequences},
our argument proceeds directly at the polynomial level, rather than via Furstenberg's lifting trick. This makes our setting considerably more delicate, since at a crucial point (see part 7 in the proof of \cref{thm:nil+null}) one must explicitly identify the generic orbit closures of certain parametric polynomial sequences, which at first appears to be a daunting task. In the setting of \cite{leibman2015nilsequences}, this difficulty does not arise, because the relevant sequences there are linear, and the corresponding
orbit closure can be described directly.
The decisive technical maneuver that greatly simplifies the problem is to characterize the generic orbit closure indirectly, yet with enough precision for our purposes, through two simple structural properties that can be verified by hand: it contains the diagonal in $X \times X$, and it is symmetric with respect to the diagonal. A Goursat-type lemma (see \cref{prop:fiber_product}) then shows that any connected subnilmanifold of $X \times X$ with these two properties must be a fiber product over a quotient, and this is exactly the structural information we need (parts 5 and 6 in the proof of \cref{thm:nil+null}).
A second novel aspect of our approach is the
borrowing of basic ideas from real algebraic geometry
(\cref{sec:algebraic_geometry}).

\subsection{Proof-sketch of the Main Theorem given the main ingredients}
We seek to model the Furstenberg-Katznelson approach to multiple
recurrence for commuting transformations, which lifts recurrence
through extensions supplied by a structure theorem and passes
to inverse limits; see
\cite[Theorem~2.4 and Lemma~3.2]{furstenberg1978ergodic}.
Passing to inverse limits is more delicate in our setting,
because the recurrence property needed for compact lifting
involves measurably varying nilsequence weights.
The abstract lemma of Bergelson, Leibman, and Lesigne,
\cite[Lemma~2.1]{bergelson2008intersective}, allows unweighted
recurrence to pass from an increasing sequence of factors
to its inverse limit, but does not directly yield the
weighted statement required here.

The model case in \cref{P:MotivationTower} illustrates how to
avoid this passage to a limit. When $X=\CZ_s(T_i)$ for every
$i\in[\ell]$, the system can be obtained from the point system
through a finite tower of extensions, each compact for all
the transformations. Repeated application of \cref{C:compactlift}
then yields the required weighted recurrence.
The main difficulty lies in establishing this lifting property
for compact extensions, and \cref{sec:proof_lifting} is devoted
to its proof.

The next example provides a very simple model for testing the compact extension step. It illustrates, in a particularly transparent setting, how the main ingredients are combined to carry out this crucial part of the argument.

\begin{example}\label{L:baby2}
Let $A\subseteq \T^2$ with $m_{\T^2}(A)>0$, $\gamma\in \R$, 
and let $\phi\colon \T\to [0,1]$ be continuous, with $\phi(0)>0$. 
Let also $\rho\colon \T^2\to\T$ be measurable. Our goal is to show that (all limits below can easily be shown to exist)
\begin{equation}\label{E:phi-positive}
\lim_{N\to\infty} \frac{1}{N}\sum_{n=1}^N  \int \phi(n\rho(y,z)) \cdot {\1}_A(y,z)\cdot {\1}_A(y+n\gamma, z+ny+n^2\gamma)\, \d y \d z>0
\end{equation}
assuming that a similar property holds for Kronecker systems.
The map~$(y,z,n)\to \phi(n\rho(y,z))$ is a positive measurable varying nilsequence (the nilmanifold is $\T$).
	
As a first step, we make crucial use of a very simple case of the decomposition result in \cref{thm:nil+null}, which yields the following:
let~$A':=\{y\in \T\colon \int {\1}_A(y,z) \d z>0\}$,
there exists $a>0$ and $B\subseteq  A'$
of measure at least $m_\T(A')/2$, such that for~$y\in B$ the set
\begin{equation}\label{E:Bax}
B_{a,y}:=\Big\{n\in \N\colon \int \phi(n\rho(y,z))\cdot {\1}_{A}(y,z) \d z\geq a \Big\}
\end{equation}
contains a Bohr$_0$-set modulo a set of zero density.
Up to shrinking~$a$, we may suppose~$a\leq m_\T(A')/4$:
the inequality in~\eqref{E:Bax} remains even more true.

We denote~$\norm{f}^2_{L^2(m_y)}:=\int \abs{f(y,z)}^2 \d z$.
Next, we use a simple special case of \cref{thm:beautiful_main}: there exists a $1$-bounded function $g\in C(\T^2)$ such that 	
\begin{equation}\label{E:global}
	\norm{{\1}_A-g}_{L^2(m_{\T^2})}\leq a^2/64,
\end{equation}
and $\delta=\delta(a,g)>0$ such that for every $y\in \T$ we have
\begin{equation}\label{E:approxg}
\norm{g(y+n\gamma, z+ny+n^2\gamma)-g(y,z)}_{L^2(m_y)}\leq a/32, \quad \text{ if } \norm{(n\gamma,ny+n^2\gamma)}_{\T^2}\leq \delta.  
\end{equation}
It is crucial here that this approximation is valid for a $\nilbor$ set of $n\in \N$; this would not be the case if we were working with a general compact extension.

We deduce from \eqref{E:global} that there exists a set $G_a$ of $y\in \T$ such that
$$
m_\T(G_a)\geq 1-a
$$
and
\begin{equation}\label{E:local}
	\norm{{\1}_A-g}_{L^2(m_y)}\leq a/32, \quad y\in G_a.
\end{equation}
Up to replacing~$G_a$ by~$G_a\cap B$,
which satisfies~$m_\T(G_a\cap B)\geq m_{\T}(A')/4>0$ by our choice of~$a$,
we may suppose that for $y\in G_a$ both \eqref{E:Bax} and \eqref{E:local} hold.
From \eqref{E:local} we get that
\begin{equation}\label{E:localn}
\norm{{\1}_A(y+n\gamma,z+ny+n^2\gamma)-g(y+n\gamma,z+ny+n^2\gamma)}_{L^2(m_y)}\leq a/32, \quad y+n\gamma\in G_a.
\end{equation}

In order to show that \eqref{E:phi-positive} holds, it suffices to show that, for $(y,z)$ in a set of positive measure, we have
\begin{multline*}
\lim_{N\to\infty} \frac{1}{N}\sum_{n=1}^N\, {\1}_{B_{a,y}}(n)\cdot {\1}_{[-\delta,\delta]^2}(n\gamma,ny+n^2\gamma)\cdot {\1}_{G_a}(y)\cdot {\1}_{G_a}(y+n\gamma)\cdot \\ \int \phi(n\rho(y,z)) \cdot {\1}_A(y,z)\cdot {\1}_A(y+n\gamma,z+ny+n^2\gamma)\, dz>0.
\end{multline*}
Using \eqref{E:global}-\eqref{E:local}, we get that for those $n\in \N$ that are relevant for the average, we may lower-bound the integral by
$\int \phi(n\rho(y,z)) \cdot g^2(y,z) \d z - a/4$. Using \eqref{E:local} once again, it therefore suffices to show that
\begin{multline*}
\lim_{N\to\infty} \frac{1}{N}\sum_{n=1}^N\, {\1}_{B_{a,y}}(n)\cdot {\1}_{[-\delta,\delta]^2}(n\gamma,ny+n^2\gamma)\cdot  {\1}_{G_a}(y)\cdot  {\1}_{G_a}(y+n\gamma)\cdot  \\ \Big(\int \phi(n\rho(y,z)) \cdot {\1}_A(y,z) \d z -a/2\Big)>0.
\end{multline*}	
Using \eqref{E:Bax}, we get that for $n\in B_{a,y}$ the integral is greater than $a$, and hence the last expression is greater than
$$
\lim_{N\to\infty} \frac{1}{N}\sum_{n=1}^N\, {\1}_{B_{a,y}}(n)\cdot {\1}_{[-\delta,\delta]^2}(n\gamma,ny+n^2\gamma)\cdot  {\1}_{G_a}(y)\cdot {\1}_{G_a}(y+n\gamma)\cdot \frac{a}{2}.
$$	
So it remains to show that the previous limit is positive for a set of $y\in \T$ of positive measure.
Integrating over $y\in \T$ and using Fatou's lemma, it is easy to see that it suffices to prove
$$
\lim_{N\to\infty} \frac{1}{N}\sum_{n=1}^N\int \CN(y,n) \cdot 
  {\1}_{G_a}(y)\cdot {\1}_{G_a}(y+n\gamma)\, dy>0
$$
for some positive measurable varying nilsequence $\CN(y,n)$ (the underlying nilmanifold is the direct product of $\T$ and the Heisenberg nilmanifold) satisfying
$$
\CN(y,n)\leq \frac{1}{2}\cdot  {\1}_{B_{a,y}}(n)\cdot {\1}_{[-\delta,\delta]^2}(n\gamma,ny+n^2\gamma)
$$ 
for every $y\in B$ and all $n\in\N$
outside a set of density zero depending on $y$.
These exceptional sets do not affect the limiting averages.
The resulting positivity property now concerns Kronecker systems, and thus it holds by our standing assumption. 

Note that, although this argument has reduced the complexity of the system under consideration, the measurable varying nilsequence has become more complicated. If we were to continue with this strategy, we would
eventually reduce matters to the point system, where the positivity
property takes the form
$$
\lim_{N\to\infty} \frac{1}{N}\sum_{n=1}^N \tilde{\CN}(n) >0,
$$
for some positive nilsequence~$\tilde \CN$. This positivity property holds,
in the setting of intersective polynomials, by \cite[Proposition 4.3]{bergelson2008intersective}.
\end{example}

\subsection{Beyond polynomial iterates}\label{SS:Beyond}
The tools developed in this article provide a flexible framework for establishing multiple recurrence results 
along more general sequences $a_1,\ldots, a_\ell\colon \N\to \Z$. Such extensions require suitable variants of the weighted characteristic-factor result
in \cref{cor:partial_factor_use} and the nil-plus-null decomposition in \cref{thm:nil+null},
with the iterates $a_1,\ldots, a_\ell$ in place of $p_1,\ldots,p_\ell$. These variants must accommodate the measurable weights and auxiliary sequences
arising in the proof, while preserving positivity and the required multiple recurrence
on nilmanifolds. The Nil-Bohr$_0$ returns theorem, \cref{thm:beautiful_main}, requires no change.

Closest-integer iterates of Hardy field sequences provide a natural setting in which this approach seems promising, leading in particular to the
following problem.

  \begin{problem}\label{Pro:2}
Let $a_1,\ldots,a_\ell$ be non-zero fractional polynomials,
meaning finite real linear combinations of $n^\alpha$ with $\alpha> 0$,
and suppose that $a_i/a_j$ is non-constant whenever $i\neq j$.
Then for every commuting system 
$(X,\mathcal X,\mu, T_1,\ldots, T_\ell)$ and every $A\in\mathcal X$ with $\mu(A)>0$, we have 
$$
\liminf_{N\to\infty} \frac{1}{N}\sum_{n=1}^N \mu\big(
A\cap T_1^{-[\![a_1(n)]\!]}A
\cap\cdots\cap T_\ell^{-[\![a_\ell(n)]\!]}A
\big)>0.
$$
\end{problem}
The assumption of pairwise independence is imposed only to make our
methods applicable; we expect the recurrence conclusion to hold
without it.
When $T_1=\cdots=T_\ell$, the problem follows from
\cite[Theorem~B(iii)]{BMR2024}, even without the assumption
of pairwise independence.
For general commuting transformations, the problem remains open
already for the family $n^a, n^b,n^a+n^b$,
with distinct positive non-rational exponents $a,b$. 
The corresponding statement with integer parts in place of closest
integers fails in the generality of \cref{Pro:2}, even for a single transformation;
see~\cite[Example~1.14]{BMR2024}.

\section{Preliminaries} \label{sec:prel}

This section collects the conventions and structural facts used
throughout the paper.

\subsection{Measure preserving systems} \label{sec:prel_systems}
We specify our measure-theoretic conventions for systems, factors,
and conditional expectations.

\subsubsection{Systems}
A \emph{standard Borel space} is a measurable space~$(X,\CX)$ isomorphic
to~$(Z,\CB(Z))$ for some Polish space~$Z$ and
Borel sigma-algebra~$\CB(Z)$, and a
\emph{standard probability space} is a standard Borel space~$(X,\CX)$
together with a probability measure~$\mu$ on~$\CX$.
Throughout this paper,
a \emph{measure preserving system},
or simply \emph{system}, is a tuple~$(X,\CX,\mu,T)$,
where~$(X,\CX,\mu)$ is a standard probability space and~$T\colon X\to X$
is a measurable bijection with measurable inverse
satisfying~$\mu(T^{-1}A)=\mu(A)$ for all~$A\in\CX$.


\subsubsection{Factors}
Given sub-$\sigma$-algebras~$\CA,\CB\subseteq \CX$, we
write~$\CA\subseteq _\mu\CB$ if for every~$A\in\CA$ there
is~$B\in\CB$ with~$\mu(A\triangle B)=0$, and~$\CA=_\mu\CB$
if~$\CA\subseteq _\mu\CB$ and~$\CB\subseteq _\mu\CA$.
As we always work modulo nullsets, when only one measure is involved
we simply write~$\CA=\CB$.

A \emph{factor} of a system~$(X,\CX,\mu,T)$ is a
system~$(Y,\CY,\nu,S)$ together with a measurable
map~$\pi\colon X\to Y$ (defined on a $\mu$-conull set) such
that~$\pi_*\mu=\nu$ and~$\pi\circ T = S\circ \pi$ $\mu$-almost everywhere.
We often leave implicit the additional data and simply say that~$Y$
is a factor of~$X$, or that~$\pi\colon X\to Y$
is an \emph{extension}.

Every factor map~$\pi$ gives the
sub-$\sigma$-algebra~$\pi^{-1}\CY=\{\pi^{-1}B : B\in\CY\}\subseteq \CX$,
which is \emph{$T$-invariant}, meaning~$T^{-1}\pi^{-1}\CY=\pi^{-1}\CY$.
Conversely, for every sub-$\sigma$-algebra~$\CA\subseteq \CX$
with~$T^{-1}\CA=\CA$ there is a factor~$\pi\colon X\to Y$
with~$\pi^{-1}\CY=\CA$, and it is unique up to isomorphism of
factors~\cite{einsiedler2011ergodic}.
This gives a bijection between factors up to isomorphism and
$T$-invariant sub-$\sigma$-algebras up to nullsets.
We use both descriptions throughout,
often identifying~$\CY$ with~$\pi^{-1}\CY$.

As we work modulo nullsets, a factor determines its sub-$\sigma$-algebra only
up to $\mu$-nullsets, and for most purposes the choice of a representative is
irrelevant. It becomes relevant when several measures are involved at once, as in \cref{sec:iso_Zs}.
In such situations we represent factors by countably generated
sub-$\sigma$-algebras (for instance~$\pi^{-1}\CY$ for a factor map~$\pi$).
The details of this technicality are collected in \cref{app:measure}.


\subsubsection{Disintegration and conditional expectation}
Let~$\pi\colon (X,\CX,\mu,T)\to (Y,\CY,\nu,S)$ be an extension.
Then there exists a family~$(\mu_y)_{y\in Y}$ of Borel probability
measures on~$(X,\CX)$, called the \emph{disintegration} of~$\mu$
over~$Y$ and written~$\mu=\int_Y \mu_y \d\nu(y)$, such that
\begin{enumerate}[(i)]
    \item for every~$A\in\CX$ the map~$y\mapsto\mu_y(A)$ is measurable
        and~$\int_Y \mu_y(A)\d\nu(y)=\mu(A)$;
    \item $\mu_y\bigl(\pi^{-1}(y)\bigr)=1$ for $\nu$-almost every~$y\in Y$.
\end{enumerate}
The family is unique up to $\nu$-nullsets: any other such
family~$(\mu'_y)$ satisfies~$\mu'_y=\mu_y$ for $\nu$-almost every~$y$.
Consequently, for every~$f\in L^1(\mu)$ (with any Borel
representative) we have~$f\in L^1(\mu_y)$ for $\nu$-almost every~$y$,
the map~$y\mapsto\int f\d\mu_y$ is measurable,
and~$\int f\d\mu=\int_Y\int f\d\mu_y\d\nu(y)$.

In terms of the invariant sub-$\sigma$-algebra~$\pi^{-1}\CY$,
the disintegration computes conditional expectations:
for every~$f\in L^1(\mu)$,
\[
    \E(f\mid \pi^{-1}\CY)(x) = \int f\d\mu_{\pi(x)}
\]
for $\mu$-almost every~$x\in X$.
Conversely, this identity characterizes the disintegration
up to $\nu$-nullsets.
We will use the notation~$\norm{f}_{L^2(\mu_y)}^2=\int \abs{f}^2 \d\mu_y$.



\subsection{Compact extensions} \label{sec:compact_extensions}
We discuss three characterizations of compact extensions; many more can be
found in \cite[Chapter~9]{glasner2003ergodic}, \cite[Section~6.3]{furstenberg2014recurrence}
and \cite[Corollary~3.7]{haase2025precompactness}.
Let~$(X,\CX,\mu,T)$ be an invertible measure-preserving system, $Y$~a factor,
and~$(\mu_y)_{y\in Y}$ a corresponding disintegration.
Following Furstenberg
\cite[Section~6.3, property~$C_3$]{furstenberg2014recurrence},
we say that the extension $\pi\colon X\to Y$ is \emph{compact} if
\begin{enumerate}[(i)]
\item for every $f\in L^2(X)$ and $\ve,\delta>0$ there are a set $B\subseteq  Y$ with
	$\nu(B)>1-\ve$ and finitely many $g_1,\ldots,g_r\in L^2(X)$ such that,
	writing $f_B:=f\cdot\1_{\pi^{-1}(B)}$, for every $n\in\N$
	\[
		\min_{j=1,\ldots,r}\norm{T^nf_B-g_j}_{L^2(\mu_y)}<\delta
		\qquad\text{for almost every } y\in Y.
	\]
\end{enumerate}
Compactness is also equivalent to the following property:
\begin{enumerate}[(i)]
\setcounter{enumi}{1}
\item for every $f\in L^\infty(X)$, we have $f=0$ if and only if
\begin{equation}\label{eq:cwm}
	\limsup_{N\to\infty}\frac1N\sum_{n=1}^N
	\E\p{\abs{\E(T^nf\cdot\con f\mid Y)}^2}=0.
\end{equation}
\end{enumerate}
This equivalence is essentially contained in \cite{tao2009poincare}, but
we give the details for completeness.
Let $L^2(X\mid Y)$ denote the space of $f\in L^2(X)$ with
$\E(\abs f^2\mid Y)\in L^\infty(\mu)$, so that
$L^\infty(X)\subseteq  L^2(X\mid Y)\subseteq  L^2(X)$.
Let $X\to Z\to Y$ be the maximal intermediate factor such that $Z\to Y$ is compact.
By \cite[Proposition~2.14.10]{tao2009poincare}, a function $f\in L^2(X\mid Y)$
satisfies~\eqref{eq:cwm} if and only if $\E(f\con g\mid Y)=0$ for all $Z$-measurable
$g\in L^2(X\mid Y)$.
If $X\to Y$ is compact, then~$Z=X$; hence if $f\in L^\infty(X)$ satisfies~\eqref{eq:cwm},
we may take~$g=f$ and obtain~$\E(\abs{f^2}\mid Y)=0$,
thus~$f=0$.
Conversely, if $X\to Y$ is not compact, then $Z\neq X$, so there is a set $E$
with $\1_E\neq\E(\1_E\mid Z)$ \cite[Exercise~2.13.6]{tao2009poincare}.
The function $f:=\1_E-\E(\1_E\mid Z)$ is bounded and non-zero, and for every
$Z$-measurable $g\in L^2(X\mid Y)$ we have
$\E(f\con g\mid Y)=\E(\E(f\mid Z)\con g\mid Y)=0$;
hence $f$ satisfies~\eqref{eq:cwm}.

A third characterization is via invariant modules.
Given $f_1,\ldots,f_r\in L^2(X)$,
we call $M:=\sum_{i=1}^r L^\infty(Y) f_i$ a
\emph{finitely generated $L^\infty(Y)$-module}.
If moreover $TM\subseteq  M$, we say that $M$ is \emph{$T$-invariant}.
Then $X\to Y$ is compact if and only if
\begin{enumerate}[(i)]
\setcounter{enumi}{2}
\item the union of all finitely generated $T$-invariant $L^\infty(Y)$-modules
	has dense linear span in $L^2(X)$.
\end{enumerate}
A proof with~$Y$ ergodic is obtained
in~\cite[Theorem~22]{zorin2019compact},
and in full generality
in \cite[Corollary~3.7]{haase2025precompactness} and \cite{jamneshan2023uncountable}.
The complication for non-ergodic $Y$
is that $\dim_y(M)$ may vary with $y\in Y$, so that orthonormal bases have to be
replaced by suborthonormal ones.
If $Y$ is ergodic, then~$\dim_y(M)$ is constant
and~$M$ has a 
conditional orthonormal basis~\cite[Section~9.2]{glasner2003ergodic}:
there are $h_1,\ldots,h_r\in M$ generating $M$ such that
\[
	\E(h_i\,\con{h_j}\mid Y)=\delta_{ij}.
\]
Moreover, writing $h=(h_1,\ldots,h_r)$ and $\Lambda_{ij}:=\E(Th_i\,\con{h_j}\mid Y)$, we have
\[
	h(Tx)=\Lambda(\pi(x))\,h(x),
	\qquad\text{i.e.}\qquad
	h_i(Tx)=\sum_{j=1}^r\Lambda_{ij}(\pi(x))\,h_j(x)\quad(i=1,\ldots,r),
\]
for $\mu$-a.e.\ $x\in X$, where $\Lambda=(\Lambda_{ij})\colon Y\to\C^{r\times r}$ is
measurable and $\Lambda(y)$ is unitary for almost every $y\in Y$
(see~\cite[Definition~9.10]{glasner2003ergodic}).

A useful general property of compact extensions,
that will be used throughout, is the following
\cite[Lemma 9.12]{glasner2003ergodic}: if
\[
    X\to Y \to Z
\]
is a chain of extensions with~$X\to Z$ compact,
then the intermediate extensions~$X\to Y$ and~$Y\to Z$
are also compact.



\subsection{Nilmanifolds and nilsystems} \label{sec:prel_nilsystem}

In this section we collect known results on nilmanifolds and nilsystems,
mostly following \cite{host2018nilpotent}.

\subsubsection{Nilmanifolds}
Let~$G$ be a nilpotent Lie group and~$\Gamma$ a discrete cocompact subgroup.
The quotient~$X:=G/\Gamma$ is called a \emph{nilmanifold}.
The set~$X$ is endowed with the quotient topology and with the natural
action of~$G$ on left cosets.
If~$G$ is $s$-step, we say that~$X$ is $s$-step.
We denote by $1_G$ and~$1_X:=1_G \Gamma$
the identity elements in~$G$ and~$X$ respectively.

\subsubsection{Haar measure}
The \emph{Haar measure} on~$X$, denoted~$\m_X$,
is the unique probability measure that is invariant
under the natural action of~$G$ on~$X$.

\subsubsection{Connectedness}

Let~$G^0$ denote the connected component of~$1_G$ in~$G$.
The subgroup~$G^0\Gamma$ has finite index in~$G$.
The nilmanifold~$X$ is connected if and only if~$G=G^0\Gamma$.

\subsubsection{Simple connectedness}
We say that~$G$ is simply connected if~$G^0$ is simply connected.
Any $s$-step nilmanifold~$X=G/\Gamma$ can be represented as~$\tilde G/\tilde \Gamma$,
where~$\tilde G$ is $s$-step and simply connected.

\subsubsection{Subnilmanifolds}
Let~$H\subseteq  G$ be a closed subgroup. We say that~$H$ is \emph{$\Gamma$-rational}
(or simply \emph{rational}) if~$H\Gamma$ is closed in~$G$,
or equivalently if~$(H\Gamma)/\Gamma$ is compact.
A \emph{subnilmanifold} is a subset~$Y\subseteq  X$
of the form~$Y=gH1_X$ for some~$g\in G$
and~$H$ a rational subgroup. Equivalently, it is a subset
of the form~$Y=H x$ for some~$x\in X$ and some rational subgroup~$H\subseteq  G$
(if~$x=g1_X$, then~$Hg1_X = g\cdot H' 1_X$ with~$H':=g^{-1}Hg$).
We will often switch between these two representations of a subnilmanifold.

\subsubsection{Quotients}
Let~$L\subseteq  G$ be a normal rational subgroup. Then~$L\backslash X$
can be canonically identified with~$G/(L\Gamma)$,
and with~$(G/L)/(L\Gamma/L)$. The latter formulation shows that~$L\backslash X$
has a natural nilmanifold structure.

\subsubsection{Malcev basis}
Suppose that~$G$ is connected and simply connected.
Then there are elements~$e_1,\ldots,e_d$
such that every~$g\in G$ can be represented, in a unique way,
as~$g=\prod_{i=1}^d e_i^{s_i}$, with~$s_i\in \R$.
We call~$e_1,\ldots,e_d$ a \emph{Malcev basis}.

One can always choose~$e_1,\ldots,e_d\in \Gamma$,
and in such a way that the elements of~$\Gamma$ are precisely those with integer coordinates.
We call such a basis \emph{$\Gamma$-rational} (or simply rational).
In other words, we have~$G=\prod_{i=1}^d e_i^\R$ and~$\Gamma=\prod_{i=1}^d e_i^\Z$.
Group multiplication in a rational basis is given by rational polynomials
in the coordinates.

If~$H\subseteq  G$ is a subgroup, we say that the basis \emph{passes through~$H$}
if, up to reordering, the subgroup~$H$ is spanned by~$e_1,\ldots,e_{\dim H}$.
If~$H$ is a connected rational subgroup, one can always choose a rational
basis passing through~$H$~\cite[Theorem 5.1.6]{corwin1990representations}.
We remark that if~$H$ is not normal, such basis may lack the ``nested property''
typical of a basis passing through the lower central series,
compare~\cite[Proposition 1.2.7]{corwin1990representations}
			with~\cite[Proposition 1.2.8]{corwin1990representations} for details.

Suppose now that~$G$ is simply connected and not connected. Let~$d^0:=\dim G^0$.
Then there are elements~$e_1,\ldots,e_d\in G$ such that every~$g\in G$
can be represented, in a unique way,
as~$g=\prod_{1\leq i\leq d^0} e_i^{s_i} \cdot \prod_{d^0<i\leq d} e_i^{k_i}$, with~$s_i\in \R$
and~$k_i\in \Z$.
In general, we cannot guarantee~$e_i\in \Gamma$
(in fact, in general the best one can do is~$e_i^\delta \in \Gamma$
where~$\delta$ is the index of~$G^0\Gamma$ in~$G$,
take for example~$G=\frac{1}{2}\Z$ and~$\Gamma=\Z$ and~$\delta=2$).
Suppose, however, that~$X=G/\Gamma$ is connected, that is~$G=G^0\Gamma$.
In such situation we can choose~$e_1,\ldots,e_d\in \Gamma$
and guarantee~$\Gamma=\prod_{1\leq i\leq d} e_i^\Z$.
See~\cite[Section 1.5]{leibman2006rational} for details.

\subsubsection{Sections} \label{prel:sections}
Let~$H \subseteq  G$ be any closed subgroup.
By the Borel cross-section theorem~\cite[Theorem 5.4.3]{srivastava1998course},
the canonical projection~$G\to G/H$ admits a Borel
section~$G/H\to G$.
If~$H$ is a rational subgroup, one can obtain such a section, in an alternative way,
by using a rational Malcev basis passing through~$H$.
The advantage of the second approach is that the resulting section is not
only Borel, but Riemann (continuous in the interior of a fundamental domain).
However, we will not need this extra regularity in this paper.

\subsubsection{Rationality}
Fix a $\Gamma$-rational basis.
An element of~$G$ is called \emph{rational} if it has rational coordinates
in the basis. An element of~$X=G/\Gamma$ is called \emph{rational} if
it has the form~$g\Gamma$ for some rational element~$g\in G$.
For example, the elements~$1_G$ and~$1_X$ are rational.
Recall that a subnilmanifold~$Y\subseteq  X$ is (the image in~$X$)
of a set of the form~$gH1_X$, for some~$g\in G$ and some rational subgroup~$H\subseteq  G$.
The subnilmanifold is called \emph{rational} if~$g$ is rational.
Alternatively, a rational subnilmanifold has the form~$H x$ with~$x\in X$ a rational
point and~$H\subseteq  G$ a rational subgroup.
Rational points, subgroups, and subnilmanifolds, are studied
in detail in~\cite{leibman2006rational}.

\subsubsection{Polynomial sequences}
Let~$a_1,\ldots,a_\ell \in G$. Let~$p_1,\ldots,p_\ell$ be integral polynomials,
i.e. rational polynomials taking integer values on the integers.
Let~$n$ vary in~$\Z$.
We call~$a_1^{p_1(n)}\cdots a_\ell^{p_\ell(n)}$ a \emph{polynomial sequence} in~$G$
and~$a_1^{p_1(n)}\cdots a_\ell^{p_\ell(n)} 1_{X}$ a \emph{polynomial sequence} in~$X=G/\Gamma$. The following result, due to Leibman,
can be found in \cite[Chapter~14, Theorem~18]{host2018nilpotent}.

\begin{theorem}
\label{thm:poly_orbit_closure}
Let~$g(n):=a_1^{p_1(n)}\cdots a_\ell^{p_\ell(n)} 1_{X}$ be a polynomial sequence
in~$X=G/\Gamma$. Then there exist an integer~$k\geq 1$, a rational subgroup~$H\subseteq  G$,
and points~$x_0,\ldots,x_{k-1} \in X$, such that the sets~$H x_j$, for~$j=0,\ldots,k-1$,
are pairwise disjoint or equal, such that
\[
	\overline{\{g(n)\colon n\in \Z\}}
	= \bigcup_{j=0}^{k-1} {H x_j},
\]
and such that for~$j\in \{0,\ldots,k-1\}$ the sequence~$(g(n))_{n\in k\Z+j}$
is equidistributed in~$H x_j$.
\end{theorem}

\subsubsection{Nilsystems}

A \emph{nilsystem} is a measure preserving system~$(X,\mu,T)$ 
that is measurably isomorphic to a system of the form~$(G/\Gamma,\m_{G/\Gamma},\tau)$,
where~$G/\Gamma$ is a nilmanifold and~$\tau \in G$ acts by left multiplication.
We will often omit~$\m_{G/\Gamma}$ and write~$(G/\Gamma,\tau)$.
We think of~$(X,\mu,T)$ as an ``abstract'' nilsystem (an equivalence class up to measurable
isomorphism), and of~$(G/\Gamma,\m_{G/\Gamma},\tau)$ as a ``concrete'' topological model.
Many abstract measure theoretical properties of nilsystems are translated into
concrete algebraic properties of~$(G/\Gamma,\m_{G/\Gamma},\tau)$, and proved in this context.
Clearly, the representation~$(G/\Gamma,\m_{G/\Gamma},\tau)$ is not unique,
and for a specific task
some representations may be more effective than others.

\subsubsection{Reduced representations} \label{sec:preliminaries_reduced}
As in~\cite[Chapter 13]{host2018nilpotent},
we will sometimes work with \emph{reduced representations}.
We say that~$(G/\Gamma,\tau)$ is \emph{reduced} if the following holds:
\begin{enumerate}
\item $\Gamma$ contains no non-trivial normal subgroups of~$G$;
\item $G$ is spanned by the identity component~$G^0$ and the element~$\tau\in G$.
\end{enumerate}
Every ergodic nilsystem has a reduced
representation~\cite[Chapter 11, Section 1.2]{host2018nilpotent}.

\subsubsection{Factors of nilsystems}\label{SS:nilfactors}
Factors of an $s$-step ergodic nilsystem are
$s$-step ergodic nilsystems \cite[Chapter 13, Theorem 11]{host2018nilpotent}.
Consider a nilsystem~$(G/\Gamma,\tau)$ and let~$\alpha \in G$.
Then the map~$\pi\colon (G/\Gamma, \tau)\to (G/\Gamma,\alpha \tau \alpha^{-1})$,
given by~$\pi(g\Gamma)=\alpha g \Gamma$, is an isomorphism.
In particular~$(G/\Gamma,\tau)$ and~$(G/\Gamma,\alpha \tau \alpha^{-1})$
are two representations of the same ``abstract'' nilsystem.
Note that if~$(G/\Gamma,\tau)$ is reduced, then~$(G/\Gamma,\alpha \tau \alpha^{-1})$
is also reduced.

Let~$(G/\Gamma,\tau)$ and~$(G'/\Gamma',\tau')$ be nilsystems.
Let~$\alpha \in G'$.
Let~$\varphi\colon G\to G'$ be a surjective continuous group homomorphism
satisfying~$\varphi(\Gamma)\subseteq  \Gamma'$ and~$\varphi(\tau)=\alpha^{-1} \tau' \alpha$.
Then~$\pi\colon (G/\Gamma,\tau) \to (G'/\Gamma',\tau')$, given
by~$\pi(g 1_{G/\Gamma})=\alpha \varphi(g) 1_{G'/\Gamma'}$,
is a factor map.

Now let~$\pi \colon X\to Y$
be any factor map between ergodic nilsystems.
Fix reduced representations~$(G/\Gamma,\tau)$ and~$(G'/\Gamma',\tau')$ for~$X$ and~$Y$.
Then~$\pi$ takes the form~$\pi(g 1_{G/\Gamma})=\alpha \varphi(g) 1_{G'/\Gamma'}$
for some~$\alpha$
and~$\varphi$ as above~\cite[Chapter 13, Theorem 12]{host2018nilpotent}.

Note that~$\pi$ may not preserve the origin,
since~$\pi(1_{G/\Gamma}) = \alpha 1_{G'/\Gamma'}$.
However, up to replacing~$\pi$ by~$R_{\alpha^{-1}} \circ \pi$
and~$(G'/\Gamma',\tau')$ by~$(G'/\Gamma',\alpha^{-1} \tau' \alpha)$,
which is another representation of~$Y$, we can force the origin to be preserved.
Alternatively, we can replace~$\pi$ with~$\pi \circ R_\beta$,
where~$\beta\in G$ is any element
such that~$\varphi(\beta) = \alpha^{-1}$ (recall that~$\varphi$ is surjective).
We refer to such procedures as ``change of base point''.

Given an origin-preserving factor map~$\pi(g 1_{G/\Gamma})=\varphi(g) 1_{G'/\Gamma'}$,
let~$L:=\varphi^{-1}(\Gamma')$. Then~$L$ is a closed subgroup of~$G$ containing~$\Gamma$.
We can represent~$\pi\colon X\to Y$ as~$\pi\colon G/\Gamma \to G/L$, given
by~$\pi(g\Gamma) = g L$,
see~\cite[Chapter 13, Theorem 13]{host2018nilpotent} for details.
Now let~$L^0$ be the identity component of~$L$.
Since~$L$ contains~$\Gamma$, the identity component~$L^0$
is normal in~$G$ by~\cite[Chapter 13, Lemma 2]{host2018nilpotent}.
The factor map~$G/\Gamma \to G/L$
factorizes as
\[
	G/\Gamma \to G/(L^0 \Gamma) \to G/L.
\]
Note that~$L^0$ is a normal connected rational subgroup and~$L/(L^0\Gamma)$ is finite.
In particular, the left map is the quotient by a normal connected rational
subgroup, and the right map is finite-to-one, see Examples~$1$ and~$2$
in \cite[Chapter 13, Section 1.3]{host2018nilpotent}.


\subsection{Host--Kra seminorms, dual functions, and the structure theorem}
\label{sec:prel_HK}

We collect the facts on Host--Kra seminorms, characteristic factors,
dual functions and dual sequences, and the Host--Kra structure theorem that
are used throughout the paper. Our main reference is
\cite[Chapters~8, 9, 12 and~16]{host2018nilpotent}.

\subsubsection{Host--Kra seminorms and characteristic factors}

Let $(X,\mu,T)$ be a system and $f\in L^\infty(\mu)$. Let $\CI$ denote the
$\sigma$-algebra of $T$-invariant sets, and define
$\snorm{f}_1:=\norm{\E(f\mid \CI)}_{L^2(\mu)}$; if the system is ergodic,
then $\snorm{f}_1=\abs{\int f\d\mu}$. For $k\geq 1$ define inductively
\[
	\snorm{f}_{k+1}^{2^{k+1}} := \lim_{H\to \infty}
		\frac{1}{H} \sum_{h=0}^{H-1} \snorm{f\cdot T^h \con f}_k^{2^k}.
\]
The limit exists, and for every $k\geq 1$ the map $\snorm{\cdot}_k$ is a
seminorm on $L^\infty(\mu)$, called the \emph{Host--Kra seminorm of order~$k$},
and $\snorm{f}_k\leq \snorm{f}_{k+1}$
for every $k\geq1$ \cite[Chapter~8]{host2018nilpotent}.

For every $k\geq 0$ there exists a unique $T$-invariant sub-$\sigma$-algebra
$\CZ_k$ of $X$ such that, for every $f\in L^\infty(\mu)$,
\[
	\E(f\mid \CZ_k)= 0 \iff \snorm{f}_{k+1}=0,
\]
see \cite[Chapter~9]{host2018nilpotent}. The corresponding factor of $X$ is
denoted by $Z_k(X)$ and is called the \emph{Host--Kra characteristic factor}
(or \emph{structure factor}) of order~$k$. When no confusion may arise, we use
the same symbol for a factor and the corresponding invariant
sub-$\sigma$-algebra. Note the shift of index: $\CZ_k$ is characterized by the
seminorm of order $k+1$. For instance, $\CZ_0=\CI$, and if $X$ is ergodic,
then $Z_1(X)$ is the Kronecker factor of $X$. By the monotonicity of the
seminorms, $\CZ_0\subseteq  \CZ_1\subseteq  \CZ_2\subseteq \cdots$.

A system $(X,\mu,T)$ is a \emph{system of order~$k$}, or a
\emph{$\CZ_k$-system} for short, if $Z_k(X)=X$; equivalently,
$\snorm{\cdot}_{k+1}$ is a norm on $L^\infty(\mu)$.
For every system $X$, the factor $Z_k(X)$ is a system of order~$k$, that is,
$Z_k(Z_k(X))=Z_k(X)$ \cite[Chapter 9, Corollary 9]{host2018nilpotent}.
Throughout the paper we repeatedly use the following properties:
\begin{enumerate}[(i)]
\item\label{it:Zk-increasing}
	every $\CZ_k$-system is a $\CZ_\ell$-system for every $\ell\geq k$
    \cite[Chapter 9, Corollary 17]{host2018nilpotent};
\item\label{it:Zk-factor}
	every factor of a $\CZ_k$-system is a $\CZ_k$-system
    \cite[Chapter 9, Corollary 17]{host2018nilpotent};
\item\label{it:Zk-invlimit}
	an inverse limit of $\CZ_k$-systems is a $\CZ_k$-system
    \cite[Chapter 9, Theorem 20]{host2018nilpotent}. 
\end{enumerate}

\subsubsection{Dual functions} \label{subsubsec:dual-functions}

We denote elements of $\{0,1\}^k$ by binary strings
$\epsilon=\epsilon_1\cdots\epsilon_k$, and write
$\abs{\epsilon}:=\sum_{i=1}^k \epsilon_i$ and
$\m\cdot \epsilon := \sum_{i=1}^k m_i \epsilon_i$ for
$\m=(m_1,\ldots,m_k)\in\N^k$. We write~$\CC$ for the complex conjugation
operator, $\CC f:=\con f$, so that $\CC^{j}f$ equals $f$ if $j$ is even and
$\con f$ if $j$ is odd.
Let $(X,\mu,T)$ be a system, $k\geq 1$, and $f\in L^\infty(\mu)$. Define
\[
	D_k(f) := \lim_{M\to\infty}
		\frac{1}{M^k}
		\sum_{\m\in \{1,\ldots,M\}^k}\,
		\prod_{\substack{\epsilon \in \{0,1\}^k\\ \epsilon \neq 00\cdots 0}}
		\CC^{\abs{\epsilon}}\, T^{\m\cdot \epsilon}f,
\]
where the limit exists in $L^2(\mu)$; for the proof of convergence and
equivalent definitions see \cite[Chapter~8]{host2018nilpotent}.
We call any function of the form $D_k(f)$ a \emph{dual function of level~$k$}.
For any~$f\in L^\infty(\mu)$ we have
\[
	\int f \cdot D_k(f) \d\mu = \snorm{f}_k^{2^k},
\]
and the following facts are equivalent: $\E(f\mid \CZ_{k})\neq 0$, $\snorm{f}_{k+1}\neq 0$,
$f$ correlates with some dual function of level~$k+1$.
In fact, the linear span of $\{D_{k+1}(g) : g\in L^\infty(\mu)\}$ is dense in
$L^1(X,\CZ_k,\mu)$ \cite[Chapter~9, Theorem~12]{host2018nilpotent}.

Now let $(X,\mu,T)$ be an ergodic nilsystem.
If $f\in C(X)$, then the averages defining $D_k(f)$ converge for every
$x\in X$, and the limit is a continuous function on~$X$
\cite[Chapter~12]{host2018nilpotent}.
Moreover, if $X$ is a $(k-1)$-step nilmanifold, then the linear
span of $\{D_k(f) : f\in C(X)\}$ is dense in $C(X)$ with respect to the
uniform norm \cite[Chapter~12, Theorem~25]{host2018nilpotent}.

\subsubsection{Dual sequences}

We also use the companion notion of dual sequence,
as defined in~\cite{frantzikinakis2015} by the first author.
A \emph{dual sequence of level~$k$} is a sequence $(\CD(n))_{n\in \N}$ of the form
\[
	\CD(n)=\lim_{M\to\infty} \frac{1}{M^k}
	\sum_{\m\in \{1,\ldots,M\}^k}\,
	\prod_{\substack{\epsilon \in \{0,1\}^k\\ \epsilon \neq 00\cdots 0}}
	\CC^{\abs{\epsilon}}\, d(n+\m\cdot \epsilon),
\]
where $(d(n))_{n\in \N}$ is a bounded sequence of complex numbers such that the
above limit exists for every $n\in \N$.

Dual functions on nilsystems give rise to dual sequences: let
$(X=G/\Gamma,\mu,T)$ with $Tx=bx$ be an ergodic nilsystem, let $f\in C(X)$,
and let $F:=D_k(f)$ be the continuous dual function of level $k$ from
Section~\ref{subsubsec:dual-functions}. Then for every $x\in X$ (in particular
for $x=1_X$, the image of the identity of $G$), the sequence
$\CD(n):=F(b^n x)$ is a dual sequence of level~$k$. Indeed, since the averages
defining $D_k(f)$ converge at every point, we can take $d(n):=f(b^n x)$.

\subsubsection{The Host--Kra structure theorem}

We use in an essential way the following
result~\cite[Chapter 16, Theorem 1]{host2018nilpotent}.

\begin{theorem}[Ergodic structure theorem]
\label{thm:HK-structure}
Let $k\in \N$. An ergodic system is a system of order~$k$ if and only if it is
an inverse limit of a sequence of ergodic $k$-step nilsystems.
\end{theorem}


\subsection{Nilsequences, $\nilbor$ sequences, and nullsequences} \label{sec:prel_nil_definitions}
In this section we define nilsequences, $\nilbor$ sequences,
and variants used in the paper.

\subsubsection{Nilsequences and their step}
Let $X=G/\Gamma$ be an $s$-step nilmanifold, let $a\in G$, and let~$F\in C(X)$.
The sequence $\CN(n):=F(a^n 1_X)$, $n\in\mathbb{Z}$, is called a
\emph{basic $s$-step nilsequence}, and $s$ is its \emph{step}; see
\cite[Section~2]{leibman2015nilsequences}. More generally, for
$a_1,\ldots,a_\ell\in G$,
\[
\CN(n_1,\ldots,n_\ell):=F\p{a_1^{n_1}\cdots a_\ell^{n_\ell}1_X},\qquad
(n_1,\ldots,n_\ell)\in\mathbb{Z}^\ell,
\]
is called an ($\ell$-dimensional) $s$-step nilsequence. If, in addition,
$p_1,\ldots,p_\ell\in\mathbb{Z}[n]$ have degree at most $d$, then
$g(n)=a_1^{p_1(n)}\cdots a_\ell^{p_\ell(n)}$ is a \emph{polynomial sequence}
in $G$.
The resulting sequence $\CN(p_1(n),\ldots,p_\ell(n))=F(g(n)1_X)$
is itself a basic nilsequence, of step at most $2sd$
\cite[Section~2]{leibman2015nilsequences}.
Leibman \cite{leibman2015nilsequences} defines a general nilsequence as a
uniform limit of basic nilsequences. Since most nilsequences in this paper are
basic, we usually omit the word ``basic''.

\subsubsection{Nullsequences}
A bounded sequence~$e\colon\mathbb{N}\to\mathbb{C}$
is a \emph{nullsequence} if it converges to~$0$ in density, that is
\[
\lim_{N\to\infty}\frac{1}{N}\sum_{n=1}^{N}\abs{e(n)}=0,
\]
or, equivalently, if the same holds with $\abs{e(n)}^2$ in place of
$\abs{e(n)}$.

\subsubsection{Positive nilsequences and $\nilbor$ sequences}
Nil-Bohr$_0$ sets were first considered by Host and
Kra~\cite{host2011nilbohr}.
We will consider nilsequences satisfying the following positivity property.

\begin{definition} \label{def:pos_nilsequence}
Let~$G/\Gamma$ be a nilmanifold and~$a_1,\ldots,a_m\in G$.
Let~$f \colon G/\Gamma\to [0,1]$ be a continuous function such that~$f(1_{G/\Gamma})>0$.
We call \emph{positive nilsequence} a sequence~$\CN \colon \Z^m\to [0,1]$ of the form
\[
    \CN(n_1,\ldots,n_m)=f(a_1^{n_1}\cdots a_m^{n_m} 1_{G/\Gamma}).
\]
Let~$U$ be a non-empty open neighborhood of~$1_{G/\Gamma}$.
We call \emph{Nil-Bohr$_0$ sequence} a sequence~$\CM \colon \Z^m\to \{0,1\}$ of the form
\[
    \CM(n_1,\ldots,n_m)=\1_U(a_1^{n_1}\cdots a_m^{n_m} 1_{G/\Gamma}),
\]
and we call \emph{Nil-Bohr$_0$ set} a set of the
form~$\{n\in \Z^m\colon \CM(n)=1\}$.
\end{definition}

In other words, a $\nilbor$ sequence is the indicator function of a $\nilbor$ set.

A key observation is that the product of finitely many $\nilbor$ sequences
is again a $\nilbor$ sequence. Equivalently, the intersection of finitely
many $\nilbor$ sets is a $\nilbor$ set.
If $\nilbor$ sequences~$(\CM_i(n_1,\ldots,n_m))_{i=1,\ldots,k}$ are
defined on the nilmanifolds~$(G_i/\Gamma_i)_{i=1,\ldots,k}$ respectively,
then the pointwise product~$\prod_{i=1}^k \CM_i(n_1,\ldots,n_m)$ can be defined on the
cartesian product~$\prod_{i=1}^k G_i/\Gamma_i$.

The following observation implies that, for the purposes of our paper,
positive nilsequences and $\nilbor$ sequences are essentially interchangeable.

\begin{lemma} \label{lem:nilsequence=nilbohr}
The following facts are true:
\begin{enumerate} [(i)]
    \item If~$\CM$ is a $\nilbor$ set, then we have
    	\[
			\CM(n_1,\ldots,n_m) \geq \CN(n_1,\ldots,n_m)
		\]
        for some positive nilsequence~$\CN$.
    \item If~$\CN$ is a positive nilsequence, then for every~$0< \ve < \CN(0,\ldots,0)$
    	we have
    	\[
			\CN(n_1,\ldots,n_m) \geq \ve \cdot \CM(n_1,\ldots,n_m)
		\]
		for some $\nilbor$ set~$\CM$.
\end{enumerate}
\end{lemma}
\begin{proof}
$(i)$~If~$U$ is a non-empty open neighborhood of~$1_{G/\Gamma}$,
there is a continuous function~$f \colon G/\Gamma \to [0,1]$ with~$f(1_{G/\Gamma})>0$
such that~$\1_{U}(x) \geq f(x)$ for all~$x\in G/\Gamma$.
$(ii)$~If~$f \colon G/\Gamma \to [0,1]$ is a continuous function
with~$f(1_{G/\Gamma})>0$, then for every~$0<\ve< f(1_{G/\Gamma})$
the set~$U:=f^{-1}(\ve,1]$ is open and contains~$1_{G/\Gamma}$.
Let~$\CM:=\1_{(\ve,1]}\circ \CN$.
\end{proof}

\subsubsection{Measurable varying nilsequences}
We will need a notion of measurable varying nilsequence, but not the most general one. Requiring, as we do in the next definition,
that everything lives on a single nilmanifold suffices for most
of our paper, and simplifies notation.

\begin{definition} \label{def:mv_nilsequence}
Let~$X$ be a standard measure space.
Let~$G/\Gamma$ be a nilmanifold and~$a_1,\ldots,a_m \colon X\to G$ measurable maps.
Let~$f \colon X\times G/\Gamma\to [0,1]$ be a measurable function
such that~$f(x,\cdot)$ is continuous and~$f(x,1_{G/\Gamma})>0$ for almost every~$x\in X$.
We call \emph{measurable varying positive nilsequence}
a sequence~$\CN \colon X\times \Z^m\to [0,1]$ of the form
\[
    \CN(x,n_1,\ldots,n_m)=f(x,a_1(x)^{n_1}\cdots a_m(x)^{n_m}1_{G/\Gamma}).
\]
Let~$(U(x))_{x\in X}$ be a measurable varying
non-empty open neighborhood of~$1_{G/\Gamma}$.
We call \emph{measurable varying Nil-Bohr$_0$ sequence}
a sequence~$\CM \colon X\times \Z^m\to \{0,1\}$ of the form
\[
    \CM(x,n_1,\ldots,n_m)=\1_{U(x)}(a_1(x)^{n_1}\cdots a_m(x)^{n_m}1_{G/\Gamma}).
\]
\end{definition}

\cref{lem:nilsequence=nilbohr} has the following measurable varying version.

\begin{lemma} \label{lem:nilsequence=nilbohr_measurable}
The following facts are true.
\begin{enumerate} [(i)]
    \item If~$\CM$ is a measurable varying $\nilbor$ set, then we have
    	\[
			\CM(x,n_1,\ldots,n_m) \geq \CN(x,n_1,\ldots,n_m)
		\]
		almost everywhere,
        for some measurable varying positive nilsequence~$\CN$.
    \item If~$\CN$ is a measurable varying positive nilsequence
    	and~$\CN(x,0,\ldots,0)\geq c$ almost everywhere for some~$c>0$,
    	then for every~$0< \ve < c$ we have
    	\[
			\CN(x,n_1,\ldots,n_m) \geq \ve \cdot \CM(x,n_1,\ldots,n_m)
		\]
		almost everywhere, for some measurable varying $\nilbor$ set~$\CM$.
\end{enumerate}
\end{lemma}
\begin{proof}
It follows from the proof of \cref{lem:nilsequence=nilbohr}.
\end{proof}

Measurable varying $\nilbor$ sequences occur naturally as follows.

\begin{lemma} \label{lem:F_nilbohr}
Let~$X$ be a standard measure space.
Consider a nilmanifold~$G/\Gamma$,
a measurable varying point~$\sigma \colon X\to G/\Gamma$,
a measurable varying element~$\tau \colon X\to G$,
and a continuous function~$F \colon G/\Gamma \to \C$. Then for every~$\ve>0$ the set
\begin{equation}\label{eq:F_nilbohr}
	\CM(x,n):= \1_{[0,\ve)}(\abs{F(\tau(x)^n \sigma(x))-F(\sigma(x))})
\end{equation}
is a measurable varying $\nilbor$ set~$\CM \colon X\times\Z\to \{0,1\}$.
\end{lemma}
\begin{proof}
Choose a measurable lift $s\colon X\to G$ such that  $s(x)\Gamma=\sigma(x)$.
Define~$f \colon X\times G/\Gamma\to \C$ as~$f(x,z):=\abs{F(s(x) z)-F(\sigma(x))}$.
Then for every~$x\in X$ the function~$f(x,\cdot)$ is continuous
and~$f(x,1_{G/\Gamma})=0$.
Thus~$U(x):=\{z\in G/\Gamma\colon f(x,z)\leq \ve\}$ is a measurable varying open neighborhood
of~$1_{G/\Gamma}$.
Moreover, define~$a(x):=s(x)^{-1}\tau(x)s(x)$.
Then
\[
    \abs{F(\tau(x)^n \sigma(x))-F(\sigma(x))}
    = f(x,a(x)^n 1_{G/\Gamma}).
\]
In particular~$\CM(x,n)=\1_{[0,\ve)}(f(x,a(x)^n 1_{G/\Gamma}))=\1_{U(x)}(a(x)^n 1_{G/\Gamma})$
is a measurable varying $\nilbor$ set.
\end{proof}

Each previous definition involves a single nilmanifold. In rare instances,
we will deal at once with a finite or countable collection of nilmanifolds,
varying measurably. We give one final definition to include this case.

\begin{definition} \label{def:composite}
Let~$X$ be a standard measure space.
We call~$\CN\colon X\times \Z^m \to [0,1]$
a \emph{composite measurable varying positive nilsequence}
if there is
a measurable partition~$X=\bigsqcup_{i \in \N} X_i$ such that
the restriction of~$\CN$ to~$X_i\times \Z^m$
is a measurable varying positive nilsequence
(defined on a nilmanifold that may depend on~$i$).
\emph{Composite measurable varying $\nilbor$ sequences} and \emph{sets}
are defined in a similar way.
\end{definition}

Given a composite measurable varying nilsequence~$\CN \colon X\times \Z^m \to [0,1]$,
i.e. involving up to countably many nilmanifolds,
and a subset~$A\subseteq  X$ of positive measure,
by the pigeonhole principle
one can always find a subset~$A'\subseteq  A$
of positive measure such that the restriction of~$\CN$ to~$A'$ involves a single nilmanifold.
This shows that the veracity of~\cref{thm:main} remains unaffected 
if~$\CN$ is allowed to be composite.

\begin{definition} [Step] \label{def:nil_meas_step}
If~$\CN\colon X\times \Z^m\to [0,1]$ is a measurable varying nilsequence
(or measurable varying $\nilbor$ sequence), we say that~$\CN$ is $s$-step
if~$(\CN(x,n))_{n\in \N}$ is $s$-step for almost every~$x\in X$.
\end{definition}

\begin{remark} \label{rem:along_poly}
Throughout the paper,
we will often evaluate measurable varying
nilsequences~$\CN \colon X\times \Z^m\to [0,1]$
along polynomials, as
\[
	\CN'(x,n) := \CN(x,p_1(n),\ldots,p_m(n)).
\]
The resulting sequence~$X\times \Z \to [0,1]$
is also a measurable varying nilsequence.
If~$\CN(x,n)$ is $s$-step and~$\deg(p_i)\leq d$ for all~$i$,
then~$\CN'(x,n)$
is $(2sd)$-step \cite[Section~2]{leibman2015nilsequences}.
\end{remark}




	\section{Nil-plus-null decomposition}
	
	In~\cite{leibman2015nilsequences}, Leibman proved that an integral combination
	of nilsequences is equal to a nilsequence plus a nullsequence.
The purpose of this section is to prove the variant of this result given in \cref{thm:nil+null} and its corollaries in \cref{SS:cor} that are better suited for our purposes.

%
	
	\subsection{Algebraic sets} \label{sec:algebraic_geometry}

    We will borrow some basic ideas from real algebraic geometry.
	All algebraic sets in this paper will be subsets of~$\R^d$.
	
	\begin{definition}
		An \emph{algebraic set} in~$\R^d$
		is the common zero set of finitely many polynomials with real coefficients.
		An algebraic set is \emph{irreducible} if it is
		not empty  and not the union of two proper algebraic subsets.
		An algebraic set is \emph{rational} if it can be defined by polynomials
		with rational coefficients.\footnote{
			In algebraic geometry,
			the term ``rational'' usually refers to a different property,
			while the property considered here is referred to as ``defined over~$\Q$''.
		}
	\end{definition}

    If~$V\subseteq  \R^d$ is a rational algebraic subset,
    we call~$W\subseteq V$
    a \emph{rational algebraic subset} if~$W=V\cap R$
    for some rational algebraic subset~$R\subseteq  \R^d$.
    Note that there are only countably many.
    
    \begin{example}
    The algebraic set~$\{(x,y)\in \R^2 \colon xy=1\}$ is irreducible,
    despite being disconnected in the Euclidean topology.
    The algebraic set~$\{x\in \R \colon x^2=0\}$ is irreducible,
    despite the polynomial~$x^2$ being not.
    The algebraic set~$\{x\in \R\colon x^2=2\}$ is rational and not irreducible,
    the algebraic set~$\{x\in \R\colon x=\sqrt{2}\}$ is irreducible and not rational,
    and the algebraic set~$\{x\in \R\colon x^3=2\}$ is rational and irreducible.
    \end{example}
    
	The ring of polynomials with real coefficients has the Noether property, that is,
	any ideal is finitely generated.
	Many elementary facts in algebraic geometry follow from the Noether property,
	such as the following:
	\begin{itemize}
	\item Any arbitrary intersection of algebraic sets is actually a finite intersection,
	and thus an algebraic set;
	\item Any algebraic set is a union of finitely many irreducible
	algebraic sets.
	\end{itemize}
	We will use the previous facts without further mention.
	The following lemma is standard.
	
	\begin{lemma} \label{lem:var_sections}
		Let~$V\subseteq  \R^d$ be an algebraic set and~$W\subsetneqq V\times V$ a proper algebraic subset.
		Then~$V':=\{v\in V\colon \{v\}\times V \subseteq  W\}$
		is a proper algebraic subset of~$V$ and
		for all~$v\in V\setminus V'$ the intersection~$(\{v\}\times V)\cap W$
		is a proper algebraic subset of~$\{v\}\times V$.
	\end{lemma}
	\begin{proof}
		The only non-trivial claim is that
		the set~$V'$ is algebraic.
		For each~$u\in V$ the set~$V_u:=\{v\in V\colon (v,u)\in W\}$ is algebraic,
		thus~$V'=\bigcap_{u\in V} V_u$ is algebraic.
	\end{proof}
	
	\begin{corollary} \label{cor:product_irreducible}
	If~$V$ is irreducible, then~$V\times V$ is irreducible.
	\end{corollary}
	\begin{proof}
	If~$V\times V$ is  a union of two proper algebraic subsets~$W_i$, $i=1,2$,
	then~$V$ is a union of two proper algebraic
	subsets~$V_i:=\{v\in V\colon \{v\}\times V \subseteq  W_i\}$, $i=1,2$.
	\end{proof}
	
	To get a taste of how we plan to use algebraic geometry,
	we anticipate here the following corollary.
	
	\begin{corollary} \label{cor:irreducible_fubini}
	Let~$\rho$ be a finite Borel measure supported on an algebraic set~$V$
	and such that~$\rho(W)=0$ for every proper algebraic subset~$W\subsetneqq V$.
	Let~$\tilde V\subseteq  V\times V$ be the union of proper rational
	algebraic subsets of~$V\times V$. Then~$(\rho\times\rho)(\tilde V)=0$.
	\end{corollary}
	\begin{proof}
		By countable additivity,
		it suffices to prove the claim under the assumption that~$\tilde V$ is a single
		proper algebraic subset of~$V\times V$. 
		By \cref{lem:var_sections},
		the set~$V':=\{v\in V\colon\{v\}\times V \subseteq  \tilde V\}$
		is a proper  
        algebraic subset of~$V$. Then~$\rho(V')=0$,
		and for every~$v\in V\setminus V'$
		the section~$(\{v\}\times V)\cap \tilde V=:\{v\}\times \tilde V_v$
		satisfies~$\rho(\tilde V_v)=0$.
		By the Fubini Theorem, we conclude that
		\[
			(\rho\times\rho)(\tilde V)
			= \int_{V} \rho(\tilde V_v) \d \rho(v)
			= \int_{V\setminus V'} \rho(\tilde V_v) \d \rho(v) = 0. \qedhere
		\]
	\end{proof}
	
	We recall the notion of \emph{algebraic dimension}:
	the dimension of an irreducible algebraic set~$V$ is the largest integer~$k$
	such that there exists a
	sequence~$V_0\subsetneqq V_1 \subsetneqq\cdots \subsetneqq V_{k}=V$
	where each~$V_i$ is a nonempty irreducible algebraic set, and the dimension of
	a general algebraic set is the largest among the dimensions of the irreducible subsets.
	In particular, if~$V$ is an irreducible algebraic set
	and~$W\subsetneqq V$ is a proper algebraic subset, then~$\dim(W) < \dim(V)$.
	
	The following lemma is the only non-elementary result we will need
	from real algebraic geometry. It ultimately relies on the fact that the real
	field has characteristic~$0$, and that the Euclidean topology is complete.
	
	\begin{lemma} [Baire-Type Lemma] \label{lem:var_baire}
		An irreducible algebraic set in~$\R^d$ is not contained in a countable union of proper
		algebraic subsets.
	\end{lemma}
	\begin{proof}
		Let~$V\subseteq  \R^d$ be an irreducible algebraic set
		and~$W\subseteq  V$ a countable union of proper algebraic subsets of~$V$.
		The natural idea is to use the Baire Category Theorem to prove that
		not only~$V\setminus W$ is non-empty, but dense in~$V$. However, in general
		such a density claim is only true if we restrict ourselves to a neighborhood of
		a nonsingular point.\footnote{
			Definition and basic properties of
			nonsingular points are given in~\cite[Section 3.3]{bochnak2013real}.
			Existence of nonsingular points relies on the fact that the real
			field is separable~\cite[Proposition 3.3.2]{bochnak2013real}.
		}
		Let~$k$ be the algebraic dimension of~$V$.
		Let~$v\in V$ be a nonsingular point.
		Then there is a non-empty open neighborhood~$\CU$ of~$v$ in~$\R^d$
		in the Euclidean topology such that~$V\cap \CU$ is a smooth submanifold of~$\R^d$
		of dimension~$k$.
		We will prove that
		the set~$(V\setminus W)\cap \CU$ is dense in~$V\cap \CU$
		in the Euclidean topology.
		By the Baire Category Theorem, we can assume that~$W$
		consists of a single irreducible algebraic set.
		Since~$W\subsetneqq V$ and $W$ is irreducible, 
		the algebraic dimension of~$W$ is less than~$k$.
		By~\cite[Theorem 2.3.6, Corollary 2.8.9]{bochnak2013real}, the set~$W$ does not contain
		non-empty open $k$-dimensional Euclidean balls. 
		In particular~$W\cap \CU$ has empty interior in~$V\cap \CU\cong \R^k$,
		and thus $(V\setminus W)\cap \CU$ is dense in $V\cap \CU$.
		Again, we remark that intersecting with~$\CU$ is crucial.\footnote{
			Let~$V$ be an irreducible algebraic set and~$W$ a proper algebraic
			subset. Despite~$\dim(W)<\dim(V)$,
			it is false, in general, that~$V\setminus W$ is dense in~$V$ in the Euclidean
			topology. An example
			is the \emph{stick} in \emph{Whitney's umbrella}~\cite[Example 3.1.2]{bochnak2013real}.
		}
	\end{proof}
	
	
	
	\subsection{Algebraic sets in nilsystems}
	
	Let~$G/\Gamma$ be a nilmanifold with~$G$ connected and simply connected.
	After choosing a $\Gamma$-rational Malcev basis (i.e. a basis with generators in~$\Gamma$),
	we can identify~$\R^d\equiv G$, and talk about \emph{algebraic sets}
	and \emph{rational algebraic sets} in~$G$.
	Since two systems of $\Gamma$-rational coordinates differ by a polynomial diffeomorphism
	with rational coefficients, these notions do not depend on the choice
	of coordinates.
	Rational points and rational subgroups, which are particular cases
	of rational algebraic sets in~$G$, are extensively
	studied in~\cite{leibman2006rational}.
	We remark that if~$H$ is a rational subgroup, then the projection of~$H$
	in~$G/\Gamma$ is closed (it is a subnilmanifold),
	but the projection of a general rational algebraic set
	need not  be closed.
	
	Let~$V\subseteq  \R^d$ be an algebraic set. We call a map~$a \colon V\to G$
	\emph{polynomial} if it is polynomial for some $\Gamma$-rational Malcev basis
	of~$G$ (and thus for any), and \emph{rational} if such polynomials can be taken rational.
	As an example, if~$e_1,\ldots,e_d \in G$ is a $\Gamma$-rational Malcev basis,
	the map~$a(v)=e_1^{v_1}\cdots e_d^{v_d}$ is a rational polynomial (since it is the identity
	in this basis),
	and so is the map~$a(v)=(e_1^{v_1}\cdots e_d^{v_d})^n$ for any~$n\in \Z$
	(by the Baker–Campbell–Hausdorff formula).
	
	\begin{lemma} \label{lem:polynomial_belonging}
		Let~$V$ be an
		algebraic set and~$a_1,\ldots,a_\ell\colon V\to G$ rational polynomial maps, 
		$p_1,\ldots,p_\ell$ be integral polynomials, and
		define $a(v,n)=a_1(v)^{p_1(n)} \cdots a_\ell(v)^{p_\ell(n)}$.
		Let~$Z \subseteq  X$ be a finite union of rational subnilmanifolds.
		Then, for any fixed~$n\in \Z$,
		the set~$V(Z,n):=\{v\in V\colon a(v,n)\Gamma\in Z\}$ is a countable union of rational
		algebraic subsets of~$V$.
	\end{lemma}
	\begin{proof}
		The set~$Z$ is a finite union of connected rational subnilmanifolds.
		The preimage of a connected rational subnilmanifold of~$X=G/\Gamma$ in~$G$
		has the form~$gH\Gamma$
		for some connected rational subgroup~$H$ and rational element~$g\in G$.
		The subgroup~$H$ is a rational algebraic set, and multiplication is given
		by rational polynomials, thus for every~$\gamma\in \Gamma$
		the set~$gH\gamma$ is a rational algebraic subset.
		Therefore, the preimage of~$Z$ in~$G$ is countable union of rational
		algebraic sets.  
		The maps~$v\mapsto a_i(v)$
		and~$(a_1,\ldots,a_\ell)\mapsto a_1^{p_1(n)}\cdots a_\ell^{p_\ell(n)}$
		are rational and polynomial. Since composition of rational polynomials
		is a rational polynomial, the claim follows.
	\end{proof}
	\begin{lemma} \label{lem:intersective_rational}
		Let~$G/\Gamma$ be a nilmanifold.
		Let~$a(n)=a_1^{p_1(n)}\cdots a_\ell^{p_\ell(n)}\in G$ be a polynomial sequence
		with~$p_1,\ldots,p_\ell$ jointly intersective.
		Then the orbit closure of~$a(n)\Gamma$ is a finite union of connected
		rational subnilmanifolds, visited uniformly periodically,
		and one of which contains the origin.
	\end{lemma}
	\begin{proof}
		The orbit closure of a polynomial sequence in a nilmanifold
		has finitely many connected components, each equal to a subnilmanifold,
		and these components
		are visited uniformly periodically~\cite[Chapter 14, Theorem 18]{host2018nilpotent}.
		By~\cite[Theorem 3.2]{leibman2006rational}, all components are
		rational if one component is rational.
		By joint intersectivity~\cite[Proposition 2.4]{bergelson2008intersective},
		one component must contain the origin~$1_{G/\Gamma}$
		(which is a rational point), and thus it must be rational
		by~\cite[Proposition 2.10]{leibman2006rational}. We conclude that all components are rational.
	\end{proof}
	
	The introduction of algebraic sets in our context is motivated by the following lemma,
	which is a variation of~\cite[Theorem 2.1]{leibman2007},
	and characterizes the generic orbit of a polynomial sequence parametrized
	by an algebraic set.
	
	\begin{lemma} [Generic Orbits] \label{lem:generic}
		Assume the hypotheses of \cref{lem:polynomial_belonging},
		and furthermore assume that~$p_1,\ldots,p_\ell$ are jointly intersective.
		Let~$Y$ be the smallest finite union of rational subnilmanifolds of~$X$
		containing~$a(v,n)\Gamma$ for all~$v\in V$ and~$n\in \Z$,
		and let~$V_0$ be the countable union of proper rational algebraic subsets of~$V$.
		Then for every~$v\in V\setminus V_0$
		we have~$\overline{\{a(v,n)\Gamma\colon n\in \Z\}}=Y$.
	\end{lemma}
	
	Note that if $V$ is irreducible, then $V\setminus V_0$ is non-empty by \cref{lem:var_baire}. 
	
	\begin{proof}
		Fix~$v\in V$ and consider the orbit closure~$Z(v)=\overline{\{a(v,n)\Gamma\colon n\in \Z\}}$.
		By \cref{lem:intersective_rational}, the set~$Z(v)$
		is  a finite union of rational subnilmanifolds of~$X$.
		In particular~$Z(v)\subseteq  Y$.
		Let~$v\in V$ be such that~$Z(v) \subsetneqq Y$.
		Then, by definition of~$Y$,
		there is~$n\in \Z$ such that~$a(w,n)\Gamma\notin Z(v)$
		for some~$w\in V$, and thus by \cref{lem:polynomial_belonging}
		the set~$V(Z(v),n):=\{w\in V\colon a(w,n)\Gamma\in Z(v)\}$
		is a countable union of \emph{proper} rational algebraic subsets
		of~$V$. Since~$v\in V(Z(v),n)$, in particular~$v\in V_0$.
		This proves that~$V_0$ has the desired property.
	\end{proof}

	
	\subsection{Factorization of polynomial sequences}
	
	Suppose that a polynomial sequence $a(n)=a_1^{p_1(n)}\cdots a_\ell^{p_\ell(n)}$
	has values in a subnilmanifold, say~$a(n)\Gamma \in H\Gamma$ for all~$n\in \Z$,
	for some rational subgroup~$H$. In general, it is not true that~$a_i\in H$,
	nor that~$a_i\Gamma \in H\Gamma$: cancellations can occur.
	However, as used before in~\cite[Proposition 4.3]{bergelson2008intersective},
	one can ``factorize'' such a polynomial sequence inside~$H$, in a sense that will
	be made precise below.
	First, we give an example.
	
	\begin{example}
		Let~$G/\Gamma$ be the Heisenberg nilmanifold. Consider the polynomial sequence
		\begin{equation} \label{eq:ex_factorization}
			a_1^{n} a_2^{n^2} a_3^{-n^2-n} \Gamma
			:=
			\begin{pmatrix} 1 & v_{1} & v_{3} \\ 0 & 1 & v_{2} \\ 0 & 0 & 1
				\end{pmatrix}^{\hspace{-.3em} n}
			\begin{pmatrix} 1 & v_{4} & v_{6} \\ 0 & 1 & v_{5} \\ 0 & 0 & 1
				\end{pmatrix}^{\hspace{-.3em} n^2}
			\begin{pmatrix} 1 & v_{7} & v_{9} \\ 0 & 1 & v_{8} \\ 0 & 0 & 1
				\end{pmatrix}^{\hspace{-.3em} -n^2-n}
			\Gamma
		\end{equation}
		and suppose that~$(a_1,a_2,a_3)$ belongs to the rational algebraic set defined by
		\[
		v_1=v_4=v_7, \quad v_3=v_6=v_8=v_9=0, \quad v_1v_2 = 1.
		\]
		Then, an explicit computation reveals that~\eqref{eq:ex_factorization} is equal
		(for every~$n\in \Z$) to
		\[
		\begin{pmatrix} 1 & 0     &  0    \\ 0 & 1 & v_{2} \\ 0 & 0 & 1
			\end{pmatrix}^{\hspace{-.3em} n}
		\begin{pmatrix} 1 & 0     &  0    \\ 0 & 1 & v_{5} \\ 0 & 0 & 1
			\end{pmatrix}^{\hspace{-.3em} n^2}
		\begin{pmatrix} 1 & 0     & v_1 v_5\\ 0 & 1 & 0    \\ 0 & 0 & 1
			\end{pmatrix}^{\hspace{-.3em} n^3}
		\begin{pmatrix} 1 & 0     & \frac{1}{2} v_1 v_5\\ 0 & 1 & 0 
			\\ 0 & 0 & 1
			\end{pmatrix}^{\hspace{-.3em} n^4-n^2} \Gamma
		=: b_1^n b_2^{n^2} b_3^{n^3} b_4^{n^4-n^2} \Gamma.
		\]
		Let~$H\subseteq  G$ be the subgroup of matrices that are zero in the middle
		coordinate of the first row.
		While~\eqref{eq:ex_factorization} is contained in~$H\Gamma$,
		the elements~$a_1,a_2,a_3$ do not belong to~$H$
		(since~$v_1v_2=1$, we cannot have~$v_1=0$). However, we have rewritten the sequence
		using elements~$b_1,b_2,b_3,b_4\in H$. Note that the coordinates~$b_i$
		are rational polynomials in the coordinates of~$a_1,a_2,a_3$.
		Moreover, note that the union of the new family of polynomials~$(n,n^2,n^3,n^4-n^2)$
		with the old family of polynomials~$(n,n^2,-n^2-n)$ forms a jointly intersective family.
	\end{example}

	\begin{proposition} [Factorization of Polynomial Sequences] \label{prop:factorization}
		Let~$X=G/\Gamma$ be a connected nilmanifold with~$G$ simply connected.
		Let~$H$ be a connected $\Gamma$-rational subgroup of~$G$.
		Let~$a_1^{p_1(n)}\cdots a_\ell^{p_\ell(n)} \in G$ be a polynomial sequence
		such that~$a_1^{p_1(n)}\cdots a_\ell^{p_\ell(n)} \Gamma \in H\Gamma$
		for all~$n\in \Z$. Then there are elements~$b_1,\ldots,b_m\in H$
		and integral polynomials~$q_1(n),\ldots,q_m(n)$ such that:
		\begin{enumerate} [(i)]
			\item $a_1^{p_1(n)}\cdots a_\ell^{p_\ell(n)} \Gamma
			= b_1^{q_1(n)} \cdots b_m^{q_m(n)} \Gamma$ for all~$n\in \Z$;
			\item in any $\Gamma$-rational Malcev basis, the coordinates of each~$b_i$
			are rational polynomials in the coordinates of~$a_1,\ldots,a_\ell$;
			\item if~$p_1,\ldots,p_\ell$ are jointly intersective,
			then~$p_1,\ldots,p_\ell,q_1,\ldots,q_m$ are jointly intersective.
		\end{enumerate}
	\end{proposition}
	\begin{proof}
		We adapt the proof of~\cite[Proposition 4.1]{bergelson2008intersective},
		which essentially gives~$(i)$, and show that it also gives~$(ii)$ and~$(iii)$.
		Choose a $\Gamma$-rational Malcev basis~$e_1,\ldots,e_{d^0} \in \Gamma$ for~$G^0$
		such that~$e_1,\ldots,e_{\dim H}$ span~$H$
		(the existence of such a basis is proven
		in~\cite[Theorem 5.1.6]{corwin1990representations}).
		If~$G$ is not connected, extend the basis from~$G^0$ to~$G$ by adding
		discrete generators~$e_{d^0+1},\ldots,e_{d} \in \Gamma$
		(the assumption that~$X$ is
		connected is crucial here\footnote{
		If~$G$ is connected, the basis~$e_i$ is always chosen so that~$G=\prod_i e_i^\R$
		and~$\Gamma=\prod_i e_i^\Z$.
		The second identity cannot be guaranteed if~$G$ is not connected.
		In fact,
		the best one can do is~$G=\prod_{i\leq d^0} e_i^\R \prod_{i>d^0} e_i^\Z$
		with~$e_i^\delta\in \Gamma$, where~$\delta$ is the index of~$G^0\Gamma$ in~$G$.
		For example, let~$G=\frac{1}{2}\Z$,
		$\Gamma=\Z$, and~$e=\frac{1}{2}$, then~$\Gamma=e^{2\Z}$.
		Assuming~$X=G/\Gamma$ connected, that is~$G= G^0\Gamma$,
		allows one to take~$e_i\in \Gamma$.
		Then~$\prod_{i\leq d^0} e_i^\Z \prod_{i>d^0} e_i^\Z \subseteq  \Gamma$,
		and conversely~$\prod_{1\leq i\leq d^0} e_i^{s_i} \prod_{i>d^0} e_i^{k_i} \in \Gamma$
		implies~$\prod_{1\leq i\leq d^0} e_i^{s_i} \in G^0\cap \Gamma$,
		and thus~$s_i\in \Z$, since~$G^0$ is connected.
		For more details, see~\cite[Section 1.5]{leibman2006rational}.
		}).
		We now write the elements~$a_i$ in this basis,
		say~$a_i=\prod_{j=1}^{d} e_j^{x_{i,j}}$ for~$i=1,\ldots,\ell$,
		and let~$x=(x_{i,j})_{1\leq i\leq \ell,1\leq j\leq d}$.
		Then, since multiplication in a Malcev basis is given by rational polynomials\footnote{
			We remark that
			the subgroup~$H$ might be not normal.
			The multiplication is given by a rational polynomial,
			but it lacks the nested property typical, for example, of a Malcev basis
			adapted to the lower central series. For details,
			compare~\cite[Proposition 1.2.7]{corwin1990representations}
			with~\cite[Proposition 1.2.8]{corwin1990representations}.
		}, we have
\begin{equation}\label{eq:prop:factorization}
		\prod_{i=1}^\ell a_i^{p_i(n)}
		= \prod_{j=1}^d e_j^{P_j(p_1(n),\ldots,p_\ell(n),x)}
\end{equation}
		for some rational polynomials~$P_j$ such that~$P_j(0,\ldots,0,x)=0$.
		By hypothesis, for any~$n$ the group
		element~\eqref{eq:prop:factorization} belongs to~$H\Gamma$.
		Since~$H\Gamma = \prod_{i\leq \dim H} e_i^\R
			\cdot \prod_{i>\dim H} e_i^\Z$\footnote{
		Since~$H=\prod_{i\leq \dim H} e_i^\R$,
		$\Gamma = \prod_{i} e_i^\Z$ (see previous footnote), and~$H$ is a subgroup,
		we obtain~$H\Gamma
			\subseteq  \prod_{i\leq \dim H} e_i^\R \cdot \prod_{i>\dim H} e_i^\Z$.
		The converse holds since~$e_i\in \Gamma$ for all~$i$ and~$e_i\in H$ for~$i\leq \dim H$.
		}, since representation in the basis is unique by simple connectedness,
		and since~$e_i\in \Gamma$, it follows that the image of~\eqref{eq:prop:factorization}
		modulo $\Gamma$ is equal to
		$
		\prod_{j=1}^{\dim H} e_j^{P_j(p_1(n),\ldots,p_\ell(n),x)} \Gamma.
		$
		By interpreting each expression $P_j(n_1,\ldots,n_\ell,x)$ as a polynomial in $\Q[x][n_1,\ldots,n_\ell]$,
		we have
		\[
		\prod_{j=1}^{\dim H} e_j^{P_j(p_1(n),\ldots,p_\ell(n),x)}
		=
		\prod_{j=1}^{\dim H}
		\prod_{k\in \N}
		\prod_{h_1,\ldots,h_k\in [\ell]}
		(e_j^{Q_{j,h_1,\ldots,h_k}(x)})^{p_{h_1}(n)\cdots p_{h_k}(n)},
		\]
		with~$Q_{j,h_1,\ldots,h_k}\in \Q[x]$, of which
		all but finitely many are zero.
		Let~$b_1,\ldots,b_m$ and $q_1,\ldots,q_m$
		be the elements~$e_j^{Q_{j,h_1,\ldots,h_k}(x)}$ and
		the corresponding polynomials $p_{h_1},\ldots, p_{h_k}$
		appearing in the last expression.
		By construction~$b_i\in H$, and it is clear that~$(i)$ and~$(iii)$ hold.
		We also have~$(ii)$ for a specific
		basis, and thus for all bases.
	\end{proof}

	
	\subsection{Fiber products}
	Let us recall from~\cite{leibman2015nilsequences}
	the definition and basic properties of fiber products of nilmanifolds.
	Let~$X=G/\Gamma$ be a nilmanifold with~$G$ connected.
	Let~$L$ be a connected normal rational subgroup of~$G$. Then
	the quotient nilmanifold~$L\backslash X \equiv G/(L\Gamma)$ is well-defined.
	We define the fiber product as
	\[
	X\times_{L\backslash X} X:= \{(x,y)\in X\times X\colon \ Lx=Ly\}.
	\]
	This is a rational subnilmanifold of~$X\times X$,
	corresponding to the connected rational
	subgroup~$G\times_{G/L} G := \{(g,h)\in G\times G\colon gL = hL\}$.
	
	For a continuous function~$f$ on~$X$, let~$\E(f\mid L\backslash X)$
	be the conditional expectation with respect to the projection~$X\to L\backslash X$.
	Explicitly, if~$\m_{Z}$ denotes the normalized Haar measure on a subnilmanifold~$Z$,
	we have~$\E(f\mid L\backslash X)(x)=\int_{Lx} f \d\m_{Lx}$ for~$x\in X$.
	By~\cite[Chapter 10, Proposition 16]{host2018nilpotent},
	the map~$\E(f\mid L\backslash X)$
	is continuous on~$X$ and descends to a continuous function on the nilmanifold~$L\backslash X$.
	Moreover, for any two continuous functions~$f_1,f_2$ on~$X$ we have
	\begin{equation} \label{eq:integral_fiber_product}
	\int_{X\times_{L\backslash X} X} f_1(x) f_2(y) \d \m_{X\times_{L\backslash X} X}(x,y)
	= \int_X f_1(x) \cdot \E(f_2\mid L\backslash X)(x) \d \m_X(x).
	\end{equation}

	
	The following result
	is reminiscent of the Goursat Lemma for groups~\cite[Proposition 1.6]{serre2016finite}.
	
	\begin{lemma} [Goursat-Type Lemma] \label{prop:fiber_product}
		Let~$X=G/\Gamma$ be a nilmanifold with~$G$ connected.
		Let~$Y\subseteq  X\times X$ be a connected subnilmanifold
		containing the diagonal of~$X\times X$.
		Then~$Y=X\times_{L\backslash X} X$ for some connected normal
		rational subgroup~$L$ of~$G$.
		\end{lemma}
	\begin{proof}
		Let~$K\subseteq  G\times G$ be the unique connected rational subgroup
		such that~$(K \cdot (\Gamma\times \Gamma))/ (\Gamma\times \Gamma) = Y$.
		We claim that~$K$ contains the diagonal of~$G\times G$.
		Since the subnilmanifold~$Y$ contains the diagonal of~$X\times X$,
		the set~$K \cdot (\Gamma\times \Gamma)$ contains the diagonal
		of~$G\times G$.
		In particular, the identity component of~$K \cdot (\Gamma\times \Gamma)$,
		which is equal to~$K$, must contain the identity component of the diagonal,
		which is the diagonal itself since~$G$ is connected.\footnote{
			Connectedness of~$G$ is essential here.
			Take for example~$G=\Gamma=\Z/2\Z$
			and~$Y=X\times X$ (which consists of a single point).
			Clearly~$K=\{(0,0)\}$
			satisfies~$K \cdot \Gamma\times\Gamma/\Gamma\times \Gamma = Y$,
			but it does not contain the diagonal of~$G\times G$.
		}
		
		Let~$L\subseteq  G$ be defined by~$L\times \{1\} := K\cap (G\times \{1\})$,
		that is~$L:=\{g\in G \colon (g,1)\in K\}$.
		It is clear that~$L$ is a rational subgroup.
		Furthermore, the subgroup~$L$ is normal. To see this,
		let~$g\in G$ and~$\ell\in L$, then~$(g,g)\in K$ and~$(\ell,1)\in K$, thus
		$
		(g^{-1} \ell g, 1) = (g,g)^{-1} (\ell,1) (g,g) \in K,
		$
		from which~$g^{-1} \ell g\in L$.
		We claim that~$K=G\times_{G/L} G$.
		If~$(g,h)\in G\times_{G/L} G$, then~$(g,h)=(h,h)(\ell,1)$ for some~$\ell\in L$,
		and thus~$(g,h)\in K$ since~$(h,h)$ and~$(\ell,1)$ are elements of~$K$.
		Conversely, if~$(g,h)\in K$,
		then~$(gh^{-1},1) = (g,h)(h,h)^{-1} \in L\times \{1\}$,
		showing that~$(g,h)\in G\times_{G/L} G$.
		This proves~$K=G\times_{G/L} G$.
		To see that~$L$ is connected,
		note that the map~$K\to L$ given by $(g,h)\mapsto gh^{-1}$ is continuous
		and surjective.

        Recall that~$K$ is the subgroup corresponding to~$Y$.
        We showed that~$K=G\times_{G/L} G$. It remains to show that~$Y=X\times_{L\backslash X} X$.
        Consider the map
		\[
			G\times_{G/L} G \to (G/\Gamma) \times_{L\backslash (G/\Gamma)} (G/\Gamma),
			\quad (g,h)\mapsto (g\Gamma,h\Gamma).
		\]
		The map is well defined: if $g,h\in G$ satisfy~$Lg=Lh$,
		clearly~$Lg\Gamma = Lh\Gamma$. The map is surjective:
		let~$g\Gamma, h\Gamma \in G/\Gamma$
		be such that~$Lg\Gamma= Lh\Gamma$. Then~$g=\ell h \gamma$ for
		some~$\ell\in L$ and~$\gamma\in \Gamma$, and
		thus~$(g\Gamma,h\Gamma)=(\ell h, h) \cdot(\Gamma\times \Gamma)$.
		Clearly, the kernel is~$
        (G\times_{G/L} G)\cap (\Gamma\times \Gamma)$. We conclude that
        \[
            (G\times_{G/L} G \cdot \Gamma\times\Gamma)/(\Gamma\times \Gamma)
            \cong (G/\Gamma) \times_{L\backslash (G/\Gamma)} (G/\Gamma).
        \]
        The left hand
        side is~$Y=K\cdot (\Gamma\times\Gamma)/(\Gamma\times \Gamma)$,
        the right hand side is~$X\times_{L\backslash X} X$.
		\end{proof}

	
	\subsection{Proof of \cref{thm:nil+null}}
	
	Parts 1 to 3 in the following proof are standard reductions (similar to those
	in the proof of~\cite[Proposition 4.3]{bergelson2008intersective}),
	but are here carried out
	for a collection of polynomial sequences at once, which requires some extra care.
	Part~4 is similar in spirit to~\cite[Lemma 4.5]{leibman2015nilsequences},
	but with ``subnilmanifolds'' replaced here by general algebraic subsets.
	In contrast to the argument by Leibman, at this step we completely separate
	the parametric space, from which the polynomial sequences are sampled,
	from the nilmanifold, to which the polynomial sequences are mapped.
	The parametric space is an algebraic set, and it is here
	where the inputs from real algebraic geometry play a crucial
	role. Parts~5 and 6 contain the essence of our argument,
	based on the Baire-Type \cref{lem:var_baire}
	and the Goursat-Type \cref{prop:fiber_product}.
    
	\begin{proof} [Proof of \cref{thm:nil+null}]
	We break down the proof into several parts.
    
    	\textit{1. Reduction to~$X$ connected, $G$ simply connected.}
		We use the notations~$a=(a_1,\ldots,a_\ell) \in G^\ell$,
		and~$a(n)=a_1^{p_1(n)}\cdots a_\ell^{p_\ell(n)}$.
		Up to moving to a universal cover~\cite[Chapter~10, Section~1.5]{host2018nilpotent},
		we can assume~$G^0$, and thus~$G$,
		simply connected. We claim that we can suppose~$X=G/\Gamma$ connected, which
		is the same as~$G=G^0\Gamma$.
		By discreteness and compactness, the subgroup~$G^0\Gamma$ has finite index in~$G$,
		say index~$k$. By~\cite[Proposition 3.2]{bergelson2008intersective},
		there are integers~$\alpha \in \N$ and~$\beta\geq 0$ such that the polynomials
		\begin{equation} \label{eq:modified_intersective}
			p_1(\alpha n+\beta),\ldots,p_\ell(\alpha n+\beta)
		\end{equation}
		are jointly intersective, and such that~$k\mid p_i(\alpha n+\beta)$ for all~$n\in \Z$
		and~$i=1,\ldots,\ell$. Up to replacing~$n$ by~$\alpha n+\beta$ and~$p_i(n)$ by~$p_i(\alpha n+\beta)$,
		we can suppose~$k\mid p_i(n)$ for all~$n\in \Z$ and~$i=1,\ldots,\ell$.
		Up to replacing each~$a_i$ by~$a_i^k$,
		each~$p_i(n)$ by~$\frac{1}{k} p_i(n)$
        (note that replacing~$p_i(n)$ by~$\frac{1}{k} p_i(n)$ preserves
        joint intersectivity, and
        recall that we are considering polynomials with integer values, not
        necessarily integer coefficients),
		and equipping~$(G^0\Gamma)^\ell$ with the pushforward of~$\rho$
		under~$(g_1,\ldots,g_\ell)\to (g_1^k,\ldots,g_\ell^k)$,
		we can assume~$G=G^0\Gamma$.
		
		\textit{2. Reduction to a connected generic orbit closure.}
		We first present the reduction for a fixed~$a$, then let~$a$ vary. Fix~$a$.
		By \cref{lem:intersective_rational}, the
		orbit closure of~$a(n)\Gamma$ is a finite union of connected components,
		each a rational subnilmanifold and one of them containing the origin,
		visited uniformly periodically. Let~$\alpha \in \N$ denote the period.
		By~\cite[Proposition 3.1]{bergelson2008intersective}, there is an integer~$\beta \geq 0$
		such that the polynomials~\eqref{eq:modified_intersective}
		are jointly intersective. Necessarily,
		by joint intersectivity~\cite[Proposition 4.3]{bergelson2008intersective},
		the sequence~$a(\alpha n+\beta)\Gamma$ 
		is contained in the component containing the origin.
		Therefore, up to replacing~$n$ by~$\alpha n+\beta$, we can
		suppose~$\overline{\{a(n)\Gamma\colon n\in \Z\}}=H\Gamma/\Gamma$
		for some connected rational subgroup~$H$.
		This concludes the reduction for a fixed polynomial sequence~$a(n)$.
		Now let~$a$ vary. For all integers~$\alpha \in \N$ and~$\beta \geq 0$, and
		for every connected rational subgroup~$H$, the
		set~$\{a\in G^\ell\colon \forall n\in \Z,\, a(\alpha n+\beta)\Gamma \in H\Gamma\}$
		is  Borel. Since there are only countably many triples~$(\alpha,\beta,H)$,
		and since we only need to prove the inequality~\eqref{eq:nilnul_goal},
		up to restricting~$\rho$ to a set of positive measure, we can suppose a single
		triple works $\rho$-almost everywhere. Therefore,
		we can assume~$\overline{\{a(n)\Gamma\colon n\in \Z\}}=H\Gamma/\Gamma$
		for $\rho$-almost
		every~$a\in G^\ell$.
		
		\textit{3. Reduction to~$X$ as orbit closure.}
		By \cref{prop:factorization}, once we fix a $\Gamma$-rational
		Malcev basis for~$G$, we obtain a measurable (in fact, polynomial)
		map~$G^\ell\to H^m$ (defined $\rho$-almost everywhere),
		mapping~$(a_1,\ldots,a_\ell)$ to~$(b_1,\ldots,b_m)$,
		such that~$a(n)\Gamma = b(n)\Gamma$, and where~$b(n)=b_1^{q_1(n)}\cdots b_m^{q_m(n)}$
		are as given by the proposition.
		Since~$b(n)\in H$, the isomorphism $H\Gamma/\Gamma\to H/(H\cap \Gamma)$ of $H$-spaces
		preserves the polynomial sequence. Therefore, up to replacing~$G,\Gamma,a(n)$
		by~$H,H\cap \Gamma,b(n)$
		and taking the pushforward measure
		we can suppose the following: the group~$G$ is connected,
		and for $\rho$-almost every~$a\in G^\ell$ the sequence~$a(n)\Gamma$
		is uniformly distributed in~$G/\Gamma$.
		
		\textit{4. Reduction to an algebraic set of positive measure.}
		Fix a $\Gamma$-rational Malcev basis for~$G$. Then~$G$ is identified with~$\R^d$
		for some integer~$d\geq 0$.
		We write~$a(v)$ for the element of~$G^\ell$ with coordinate~$v\in \R^{d\ell}$.
		From now on, we think of~$\rho$ as a measure on~$\R^{d\ell}$,
		rather than~$G^\ell$. By the following iterative process,
		we will reduce to the case in which~$\rho$ has support on an irreducible algebraic set, and vanishes
		on any proper algebraic subset (equivalently, on any proper irreducible algebraic subset).
		By definition, the set~$\R^{d\ell}$ is an irreducible algebraic set.
		If~$\rho$ vanishes on any proper algebraic subset~$V\subseteq  \R^{d\ell}$, we stop.
		Otherwise, we can choose
		an irreducible algebraic set~$V$ with $\rho(V)>0$ and~$\dim(V)<\dim(\R^{d\ell})$.
		If~$\rho$ vanishes on any proper algebraic subset of~$V$, we stop. Otherwise, we repeat
		the process. The process terminates in at most~$d\ell$ steps.
		We have thus reduced to an integral of the form
		$
		\int_V f(a_1(v)^{p_1(n)}\cdots a_\ell(v)^{p_\ell(n)} \Gamma) \d \rho(v),
		$
		with~$v\mapsto a_i(v)$ rational polynomial maps,
		and  such that
		for $\rho$-almost every~$v$ the orbit closure of~$a(v,n)\Gamma$ is~$G/\Gamma$, where  $a(v,n):=a_1(v)^{p_1(n)}\cdots a_\ell(v)^{p_\ell(n)}$.
		
		\textit{5. Reduction to a connected generic orbit closure for the joint sequence.}
		Consider the joint sequences~$(a(v,n)\Gamma,a(w,n)\Gamma)$ for~$v,w\in V$.
		Let~$\tilde V$ be the countable union of proper rational algebraic subsets
		of~$V\times V$.
		By \cref{cor:product_irreducible}, the algebraic set~$V\times V$
		is irreducible. By \cref{cor:irreducible_fubini},
		we have~$(\rho\times\rho)(\tilde V)=0$.
		By \cref{lem:generic}, applied to the joint sequences,
		there is~$Y\subseteq  X\times X$, a finite union of rational subnilmanifolds,
		such that~$(a(v,n)\Gamma,a(w,n)\Gamma)$ 
        is dense
		on~$Y$ for every~$(v,w)\in (V\times V)\setminus \tilde V$. 

		Next, we show that we can assume~$Y$ connected and containing the origin.
        Fix an arbitrary~$(\bar v,\bar w)\in (V\times V)\setminus \tilde V$.
        Then we can label the connected components of~$Y$ as~$Y^1,\ldots,Y^\alpha$ (possibly with repetition),
        in such a
        way that for each~$i\in \{1,\ldots,\alpha\}$ the
        subsequence~$(a(\bar v,\alpha n+i)\Gamma,a(\bar w,\alpha n+i)\Gamma)_{n\in \N}$
        has values in~$Y^i$ and is equidistributed
        on~$Y^i$, see~\cite[Chapter 14, Theorem 18]{host2018nilpotent}
        or preliminaries.
        
        By~\cite[Proposition 3.2]{bergelson2008intersective}, there is~$\beta$
        such that joint intersectivity is preserved along~$\alpha n+\beta$,
        and by~\cite[Proposition 2.4]{bergelson2008intersective}
        the component~$Y^\alpha$ must contain the origin.
        For every~$n\in \N$, define
        \[
            V_{n,i} := \{(v,w)\in V\times V \colon (a(v,\alpha n+i)\Gamma,a(w,\alpha n+i)\Gamma) \in Y^i\}.
        \]
        By \cref{lem:polynomial_belonging}, the set~$V_{n,i}$ is a countable
        union of algebraic sets.
        For every~$n\in \N$ we have~$V\times V = \bigcup_{i=1}^\alpha V_{n,i}$.
        Since~$V\times V$ is irreducible,
        by \cref{lem:var_baire}
        we must have~$V\times V=V_{n,\sigma(n)}$
        for some~$\sigma(n)\in \{1,\ldots,\alpha\}$.
        Since~$(a(\bar v,\alpha n+\beta)\Gamma,a(\bar w,\alpha n+\beta)\Gamma) \in V_{\alpha n+\beta,\beta}$,
        we must have~$\sigma(\alpha n+\beta)=\beta$ for all~$n$. That is,
        we have~$(a(v,\alpha n+\beta)\Gamma,a(w,\alpha n+\beta)\Gamma) \in Y^\alpha$ for all~$(v,w) \in V\times V$ for all~$n$.
        Moreover, for~$(v,w)\in V\times V\setminus \tilde V$
        the sequence~$(a(v,\alpha n+\beta)\Gamma,a(w,\alpha n+\beta)\Gamma)_{n\in \N}$ must be equidistributed
        in~$Y^\alpha$ by~\cref{lem:generic} and definition of~$\tilde V$.
        Therefore, up to replacing~$n$ by~$\alpha n+\beta$
        and~$Y$ by~$Y^\alpha$,
        we can suppose the following: the subnilmanifold~$Y$ is connected,
		it contains the origin, it contains~$(a(v,n)\Gamma,a(w,n)\Gamma)$
		for every~$v,w\in V$ and every~$n$,
		and it coincides with the orbit closure of~$(a(v,n)\Gamma,a(w,n)\Gamma)$
		for all~$(v,w)\in V\times V\setminus \tilde V$.
		
		\textit{6. Goursat-Type argument.}
		We claim that the subnilmanifold~$Y\subseteq  X\times X$ contains the
		diagonal of~$X\times X$.
		By construction, the set~$Y$ contains~$(a(v,n)\Gamma,a(w,n)\Gamma)$
		for every~$v,w\in V$ and every~$n$, and it is the generic orbit closure.
		Fix an arbitrary~$x\in X$. Choose~$v\in V$ such that the sequence~$a(v,n)\Gamma$
		is dense in~$X$ (by construction, the set of such~$v\in V$
		has full $\rho$-measure, in particular it is non-empty).
		Choose~$n$ such that~$a(v,n)\Gamma$ is sufficiently close to~$x$.
		Then~$(a(v,n)\Gamma,a(v,n)\Gamma)\in Y$ is sufficiently  close to~$(x,x)$.
		Since~$Y$ is closed, it follows that~$Y$ contains the diagonal.		
		By \cref{prop:fiber_product}, there is a connected normal rational
		subgroup~$L\subseteq  G$ such that~$Y=X\times_{L\backslash X} X$.
		
		\textit{7. Conclusion of the proof.}		
		We claim that
		\begin{equation} \label{eq:nil+null_identity}
		\int_V f(a(v,n)\Gamma) \d \rho(v)
		= \rho(V)\cdot \E(f\mid L\backslash X)(a(w,n)\Gamma) + e(n)
		\end{equation}
		for any choice of~$w\in V$ and for some nullsequence~$e(n)$.
		First, note that for all~$v,w\in V$ we have~$(a(v,n)\Gamma,a(w,n)\Gamma)\in Y$,
		and thus~$L a(v,n) \Gamma = L a(w,n)\Gamma$,
		and thus the value~$\E(f\mid L\backslash X)(a(w,n)\Gamma)$ does not depend on~$w\in V$.
		From now on, we consider~$w\in V$ fixed.
		It is enough to prove that 
		\[
		\lim_{N\to\infty}
		\frac{1}{N}\sum_{n=1}^N
		\abs{
			\int_V (f-\E(f\mid L\backslash X))(a(v,n)\Gamma) \d \rho(v)
		}^2 = 0.
		\]
		Up to replacing~$f$ by~$f-\E(f\mid L\backslash X)$,
		we can suppose~$\E(f\mid L\backslash X)=0$ and prove
		\[
		\lim_{N\to\infty}
		\frac{1}{N}\sum_{n=1}^N
		\int_{V\times V} (f\otimes \con{f}) (a(v,n)\Gamma,a(w,n)\Gamma) \d(\rho\times\rho)(v,w)
		= 0.
		\]
		Since~$Y$ is the $(\rho\times\rho)$-generic orbit closure
		and~$Y=X\times_{L\backslash X} X$, the previous limit is equal to
		\begin{align*}
			\int_{X\times_{L\backslash X} X} f\otimes \con{f} \d\m_{X\times_{L\backslash X} X}
			= \int_X \E(f\mid L\backslash X) \cdot \con{f} \d \m_{X} = 0,
		\end{align*}
		where in the first equality we used the identity~\eqref{eq:integral_fiber_product}.
		
		It remains to show that the sequence~$\E(f\mid L\backslash X)(a(w,n)\Gamma)$
		has the desired properties.
		Since~$f$ is continuous, non-negative, and positive at the origin,
		 its average along the $L$-orbit of the origin is also positive,
		that is~$\E(f\mid L\backslash X)(1\Gamma)>0$.
        The quotient space~$Y:=L\backslash X=L\backslash (G/\Gamma)$
        is naturally identified with the nilmanifold~$(G/L)/(\Gamma L/L)$
        \cite[Chapter~10, Section~3.3]{host2018nilpotent}.
        Let~$H:=G/L$ and~$\Lambda:= (\Gamma L)/L$.
		We identify~$\E(f\mid L\backslash X)$ with the corresponding function~$g$
		on the nilmanifold~$Y=H/\Lambda$.
		Moreover, let~$b_i \in H$ be the projection of~$a_i(w)\in G$
		modulo~$L$ (recall that~$w$ is fixed). Since~$L$ is a normal subgroup,
		this projection is a group homomorphism. In
		particular~$a_1(w)^{p_1(n)}\cdots a_\ell(w)^{p_\ell(n)}$
		is mapped to the polynomial sequence~$b_1^{p_1(n)} \cdots b_\ell^{p_\ell(n)}$.
	\end{proof}
	
	
	\subsection{Corollaries and generalizations}\label{SS:cor}    
	We now state and prove various corollaries of \cref{thm:nil+null} and of its proof.
	In \cref{thm:nil+null} the function~$f$ is fixed.
	\cref{cor:nil+nul+f} below
	shows that we can allow the function~$f$ to vary measurably.
	In the notation of \cref{sec:prel_nil_definitions},
	and up to passing to a subsequence~$\alpha n+\beta$,
	\cref{cor:nil+nul+f} gives a nil+null decomposition of the form
	\[
		\int_X \CN(x,p_1(n),\dots,p_\ell(n)) \d \mu(x) \geq \CN'(q_1(n),\ldots,q_m(n)) + null
	\]
	for any measurable varying positive nilsequence~$\CN$, for some
	positive nilsequence~$\CN'$.
	By \cref{lem:nilsequence=nilbohr_measurable}, the same holds
	with ``positive nilsequence'' replaced by ``$\nilbor$ sequence''.

	\begin{corollary} \label{cor:nil+nul+f}
	Let~$G/\Gamma$ be a nilmanifold and~$p_1,\ldots,p_\ell$ be jointly intersective polynomials.
	Let~$(X,\mu)$ be a standard probability space. Let~$f\colon X\times G/\Gamma \to [0,1]$
	be a Borel function such that~$f(x,\cdot)$ is continuous on~$G/\Gamma$
	and~$f(x,1_{G/\Gamma})>0$ for $\mu$-almost every~$x$.
	Let~$a_1,\ldots,a_\ell \colon X\to G$ be Borel measurable maps.
	Then the conclusion of \cref{thm:nil+null} holds for
	\[
	\int_{X}
	f(x,a_1(x)^{p_1(n)}\cdots a_\ell(x)^{p_\ell(n)}\cdot 1_{G/\Gamma})
	\d\mu(x).
	\]
\end{corollary}
    
	\begin{proof}
		Fix a metric on~$G/\Gamma$ inducing the topology.
		For every~$k\in \N$, let~$f_k$ be a non-negative continuous function on~$G/\Gamma$
		supported on an $(1/k)$-neighborhood
		of~$1_{G/\Gamma}$ and with maximum value~$1/k$ obtained at~$1_{G/\Gamma}$. Define
		\[
			X_k := \{x \colon f(x,g)\geq f_k(g), \text{ for all } g\in G/\Gamma\}.
		\]
		Since it suffices to check the previous inequality on a countable dense set of points $g\in G/\Gamma$, the sets~$X_k$ are Borel and satisfy~$X=\bigcup_{k=1}^\infty X_k$ up to sets of $\mu$-measure zero.
		So there exists $k\in \N$ such that~$\mu(X_k)>0$.
        Since we are only interested in a lower bound, we may restrict the integral
to some $X_k$ with positive measure. After doing so and normalizing the restricted measure, which only introduces an additional positive constant that can be absorbed into $g$, we may suppose that $f(x,\cdot)\geq f_k(\cdot)$
		everywhere. Again, since we are only interested in a lower bound,
		we can replace~$f$ by~$f_k$, that is, assume~$f$ independent of~$x\in X$.
		\cref{thm:nil+null} now applies with~$\rho$ equal to the pushforward of~$\mu$
		under the map~$(a_1,\ldots,a_\ell) \colon X\to G^\ell$.
	\end{proof}

\subsection{Extension of the nil+null result to varying measures}
We will actually need a variant of \cref{cor:nil+nul+f}, that allows for variable measures $\mu$; the exact statement is given next and is proved after we prove a suitable Borel selection theorem that is adapted to our setting.
The variable measures arise from a disintegration~$\mu=\int \mu_y \d\nu(y)$
with respect to a factor, defined in \cref{sec:prel_systems}.

\begin{proposition} \label{P:nil+nul=useful}
	Let~$G/\Gamma$ be a nilmanifold, 
	$(X,\mu)$ be a standard probability space and let $f\colon X\times G/\Gamma \to [0,1]$
	be a Borel function such that $f(x,\cdot)$ is continuous on $G/\Gamma$
	and $f(x,1_{G/\Gamma})>0$ for $\mu$-almost every~$x\in X$. Let also $m\in \N$, $q_1,\ldots, q_m$ be a family of jointly intersective polynomials,  and 
	let $a_1,\ldots,a_m \colon X\to G$ be Borel  maps.
	Furthermore, let  $\mu=\int \mu_y \d\nu(y)$ be the disintegration of $\mu$ with respect to a factor~$(Y,\nu)$.
	
	Then there are disjoint Borel sets $Y_0,Y_1,\ldots$,  where some of the $Y_i$ may be empty,  such that $Y=\bigcup_{i\geq 0} Y_i$, $\nu(Y_0)=0$,
	with the following property:  For every $i\in \N$ such that  $Y_i\neq \emptyset$, 
	there exist   integers $\alpha_i\geq 1$, $k_i \geq 1$, $\beta_i\geq 0$,  a nilmanifold~$H_i/\Lambda_i$, polynomials $r_{i,1},\ldots,r_{i,k_i}$,
	Borel functions~$g_i \colon Y_i\times H_i/\Lambda_i \to [0,1]$, where for every $y\in Y_i$ the function $g_i(y,\cdot)$  is continuous and
    satisfies~$g_i(y,1_{H_i/\Lambda_i})>0$,
	Borel maps~$b_{i,1},\ldots,b_{i,k_i}\colon Y_i\to H_i$,
	and a Borel varying nullsequence~$e_i\colon Y_i\times \N\to [-1,0]$,
	such that
	\begin{multline}\label{E:muy}
		\int_X f(x,a_1(x)^{q_1(\alpha_in+\beta_i)}\cdots a_m(x)^{q_m(\alpha_in+\beta_i)} \cdot 1_{G/\Gamma})
		\d \mu_y(x)
		\geq\\
		g_i(y,b_{i,1}(y)^{r_{i,1}(n)}\cdots b_{i,k_i}(y)^{r_{i,k_i}(n)} \cdot1_{H_i/\Lambda_i})
		+ e_i(y,n)
	\end{multline}
	for  every~$y\in Y_i$ and $n\in \N$.  Moreover, for every $i\in \N$  the family
	$$
	q_1(\alpha_in+\beta_i),\ldots,q_m(\alpha_in+\beta_i), r_{i,1}(n),\ldots,r_{i,k_i}(n)
	$$
	is jointly intersective.
\end{proposition}
\begin{remark}
	The result also applies to bounded Borel measures upon normalizing them and  introducing a constant on the right hand side. 

    Moreover, if $p_1,\ldots, p_\ell$ is any family of integer polynomials such that the family $p_1,\ldots, p_\ell, q_1,\ldots, q_m$ is jointly intersective, then we can 
add to our conclusion that the extended family 
$$
	p_1(\alpha_in+\beta_i),\ldots,p_\ell(\alpha_in+\beta_i),q_1(\alpha_in+\beta_i),\ldots,q_m(\alpha_in+\beta_i), r_{i,1}(n),\ldots,r_{i,k_i}(n)
	$$
	is jointly intersective. To see this,  apply the result with $q_1,\ldots, q_m, p_1,\ldots, p_\ell$ in place of  $q_1,\ldots, q_m$, and simply insert the additional sequences $a_{m+1}(x)=\cdots=a_{m+\ell}(x)=1_G$, raised respectively to the powers $p_1(\alpha_in+\beta_i),\ldots, p_\ell(\alpha_in+\beta_i)$, into the left hand side of \eqref{E:muy}. Since these additional terms are trivial, the expression on the left-hand side remains unchanged.
\end{remark}
\subsubsection{An abstract Borel selection theorem} We start with a Borel selection theorem that we state and prove in a somewhat abstract setting.  If $(Y,\nu)$ is a standard Borel probability space, we say that  a set  $Z$ is {\em $\nu$-null} if it is contained in a Borel set $Z'\subseteq Y$ with $\nu(Z')=0$. We will use the following classical Borel selection theorem (see for instance \cite[Theorem~18.1]{kechris2012classical}).

\begin{theorem}[Jankov-von Neumann]\label{T:SelectionJV}
Let $Y,\Theta$ be standard Borel spaces, let  $R\subseteq  Y\times \Theta$ be Borel, and let 
$$
E:=\{y\in Y\colon \text{ there exists } \theta\in \Theta \text{  such that } (y,\theta)\in R\}.
$$ 
Let $\nu$ be a Borel probability measure on $Y$. Then  there exist a Borel set $E'\subseteq  E$ such that $E\setminus E'$ is $\nu$-null and  a Borel map $\phi\colon E'\to \Theta$ such that 
$$
(y,\phi(y))\in R \quad \text{ for every } y\in E'. 
$$
\end{theorem}

We remark that the set~$E$ is analytic, but not necessarily Borel.


\begin{lemma}\label{L:SelectionAbstract}
	Let  $(Y,\nu)$ be a standard Borel probability space.
	For  $j\in \N$ let  $\Theta_j$
	 be a standard Borel space, and  let 
	 $$
	 \Phi_j\colon Y\times \N\to [0,1], \qquad \Psi_j\colon\Theta_j\times \N \to [0,1]
	 $$ be Borel maps.    Suppose that there is a  Borel set $Y'$ such that $Y\setminus Y'$ is $\nu$-null and such that for every $y\in Y'$ there exist $j=j_y\in \N$, $\theta=\theta_y\in \Theta_j$, and a nullsequence $e(y,\cdot)\colon \N\to [-1,1]$, such that 
	$$
	 \Phi_j(y,n)\geq \Psi_j(\theta,n)+e(y,n), \qquad \text{for every } n\in \N. 
	$$
	 Then there  exist disjoint Borel sets $Y_0,Y_1,\ldots$, some of which may be empty, such that 
	 $$
	 Y=\bigcup_{i\geq 0} Y_i, \qquad  \nu(Y_0)=0, 
	 $$
  and for $i\in \N$ there exist  Borel maps $\theta_i\colon Y_i\to \Theta_{i}$,  and Borel varying  nullsequences $e_i\colon Y_i\times \N\to [-1,0]$, such that 
	 \begin{equation}\label{E:Fj}
\Phi_i(y,n)\geq \Psi_i(\theta_i(y),n)+e_i(y,n), \qquad \text{for every } n\in \N \text{ and }  y\in Y_i.
	 \end{equation}
\end{lemma}
\begin{remark}
A useful observation is the following. Fix $i\in\N$, $y\in Y_i$,
and $\theta_i(y)\in\Theta_i$. Then there exists a nullsequence
$e:\N\to[-1,0]$ such that
$$
\Phi_i(y,n)\geq \Psi_i(\theta_i(y),n)+e(n)
\qquad\text{for every }n\in\N, 
$$
if and only if
$$
\lim_{N\to\infty}
\E_{n\in[N]}
\big(\Psi_i(\theta_i(y),n)-\Phi_i(y,n)\big)_+=0,
$$
where~$(\cdot)_+$ denotes the positive part.
\end{remark}
\begin{proof}
	Let  $j\in \N$. For $N\in \N$, $y\in Y$,  $\theta \in \Theta_j$, let 
	$$
	A_{j,N}(y,\theta):=\E_{n\in [N]} \big(\Psi_j(\theta,n) - \Phi_j(y,n)\big)_+.
	$$  
	Let also 
	$$
	R_j:=\big\{(y,\theta)\in Y\times \Theta_j\colon 	\lim_{N\to \infty} A_{j,N}(y,\theta)=0
	\big\}
	$$
	and 
	$$
	E_j:=\big\{y\in Y\colon 	\text{ there exists } \theta\in \Theta_j \text{ such that } (y,\theta)\in R_j
	\big\}.
	$$
	Clearly the set $R_j$ is Borel (the set $E_j$ is not necessarily Borel, but this will not bother us). 
	Our hypotheses imply that for every $y\in Y'$, there exist $j=j_y$ and $\theta=\theta_y\in \Theta_j$ such that
$(\Psi_j(\theta,n)-\Phi_j(y,n))_+\leq |e(y,n)|$
	for every $n\in\N$ for some nullsequence $e$, hence 	
	 $(y,\theta)\in R_j$ by the previous remark.  Hence, 
	$$
Y'\subseteq  \bigcup_{j\geq 1} E_j.
	$$
	
	We apply the selection result of \cref{T:SelectionJV} to the Borel sets $R_j$. We get that there exist Borel sets $E'_j\subseteq E_j$ with $ E_j\setminus E'_j$  a $\nu$-nullset (not necessarily Borel),  and Borel maps
	$\phi_j\colon E'_j\to \Theta_j$ such that $(y,\phi_j(y))\in R_j$ for every $y\in E'_j$.  If  $E''_j:=E'_j\setminus \bigcup_{h<j}E'_h$, the sets $E''_j$ are disjoint Borel sets (some of them perhaps empty) and have the same union as the sets $E_j$ up to a $\nu$-null  set. We set
	$$
	Y_i:=E''_i \text{ for }\  i\in \N \quad \text{ and } \quad   Y_0:=Y\setminus \bigcup_{i\geq 1} Y_i,
	$$
	and 
	$$
	\theta_i :=\phi_i|_{Y_i}, \quad i\in \N.
	$$
	Then $
	Y=\bigcup_{i\geq 0} Y_i$, the sets $Y_0,Y_1,\ldots$ are disjoint, and since   $Y_0\subseteq  (Y\setminus Y')\cup \bigcup_{j\geq 1}(E_j\setminus E_j')$, we have $\nu(Y_0)=0$. Moreover, for $y\in Y_i$ we let 
	$$
	e_i(y,n):=-\big(\Psi_i(\theta_i(y),n)-\Phi_i(y,n)\big)_+.
	$$
	Then $y\mapsto e_i(y,n)$ is Borel in $y$ for every fixed $n$, takes values in $[-1,0]$, and  it is a nullsequence for every $y\in Y_i$. 
It follows  that \eqref{E:Fj} holds for all $i\in \N$. 
\end{proof}

\subsubsection{Proof of \cref{P:nil+nul=useful}}
We prove \cref{P:nil+nul=useful} by first using \cref{cor:nil+nul+f}  for  fixed values of $y$ and then using  \cref{L:SelectionAbstract} to make a Borel selection of the output data.
	
	\smallskip 
	
\textit{1. Organizing the output of \cref{cor:nil+nul+f}.} We first note the following consequence of   \cref{cor:nil+nul+f}.
Let us fix the nilmanifold  $G/\Gamma$, the standard probability space $X$, and the intersective family of polynomials $q_1,\ldots,q_m$, but vary the measure~$\mu$ on~$X$.
The output data we obtain from \cref{cor:nil+nul+f}
(defined in \cref{thm:nil+null}),
as we vary~$\mu$, can be organized into a ``discrete'' part and
a ``continuous'' part as follows.
The ``discrete'' part is a countable collection 
\begin{equation}\label{E:CountableData}
(\alpha_j,\beta_j, k_j, H_j/\Lambda_j, r_{j,1},\ldots, r_{j,k_j})_{j\in \N},
\end{equation}
where $\alpha_j,m_j,k_j\geq 1$, $\beta_j\geq 0$, and   $H_j/\Lambda_j$ is a nilmanifold chosen from a countable collection
depending only on~$(G,\Gamma)$, and the family 
\begin{equation}\label{E:ji}
q_1(\alpha_jn+\beta_j),\ldots, q_m(\alpha_jn+\beta_j),
r_{j,1}(n),\ldots, r_{j,k_j}(n)
\end{equation}
is jointly intersective.
The countability comes from the fact that the rational subgroups
and  nilmanifold models (containing the origin) that occur in the proof form countable
families, and from the countability of the set of integer-valued polynomials.
The ``continuous'' part is the remaining output data, namely  the continuous function $g\in C(H_j/\Lambda_j)$ with values in $[0,1]$ and positive at $1_{H_j/\Lambda_j}$, and the elements $b_1,\ldots, b_{k_j}\in H_j$, which can vary in an uncountable set.

Our goal is to make a Borel selection of these data using \cref{L:SelectionAbstract}. As
 a preparation, we 
organize the output data as follows: For each $j\in \N$, let 
$$
C_j^+:=\{g\in C(H_j/\Lambda_j,[0,1])\colon g(1_{H_j/\Lambda_j})>0\}, 
$$
and we equip this set with the $\norm{\cdot}_\infty$-norm. We let 
$$
\Theta_j:=C_j^+\times H_j^{k_j},
$$
which is a standard Borel space with the product topology. 
An element in $\Theta_j$ has the form 
\begin{equation}\label{E:ContinuousData}
\theta=(g,b_1,\ldots, b_{k_j}), \qquad g\in C^+_j, \quad b_1,\ldots, b_{k_j}\in H_j, 
\end{equation}
and for any such $\theta$ and $n\in \N$ we let 
\begin{equation}\label{E:Psii}
\Psi_j(\theta,n):=g(b_1^{r_{j,1}(n)}\cdots b_{k_j}^{r_{j,k_j} (n)} 1_{H_j/\Lambda_j} ).
\end{equation}
Next, for $y\in Y$ and $j,n\in \N$, we define
\begin{equation}\label{E:Phii}
\Phi_j(y,n):=\int f(x,a_1(x)^{q_1(\alpha_jn+\beta_j)}\cdots a_m(x)^{q_m(\alpha_jn+\beta_j)} 1_{G/\Gamma})\d\mu_y(x). 
\end{equation}
Since the integrand is a Borel function of $x$ and the map $y\mapsto \mu_y$ is Borel, the function~$y\mapsto \Phi_j(y,n)$ is Borel. Thus, $\Phi_j\colon Y\times \N\to [0,1]$ is Borel. 
Finally, we  let 
$$
X_0:=\{x\in X\colon f(x,1_{G/\Gamma})=0\}, \quad Y':=\{y\in Y\colon \mu_y(X_0)=0\}.
$$
The set $X_0$ is Borel, and since the map $y\mapsto\mu_y(X_0)$ is Borel we get that  $Y'$ is Borel.  By our hypothesis we have $0=\mu(X_0)=\int \mu_y(X_0)\, d\nu(y)$, hence $\mu_y(X_0)=0$ for $\nu$-almost every $y\in Y$, which implies that   $\nu(Y\setminus Y')=0$. 
	
\smallskip 

\textit{2.  Applying  \cref{cor:nil+nul+f}.} For $y\in Y'$ we  have $f(x,1_{G/\Gamma})>0$  for $\mu_y$-almost every~$x\in X$, hence  \cref{cor:nil+nul+f} applies for  the nilmanifold $G/\Gamma$, the  intersective family of polynomials $q_1,\ldots, q_m$, the probability measure $\mu_y$ in place of $\mu$, the Borel function $f$ and the Borel maps $a_1,\ldots, a_m$. The output consists of the ``discrete'' data that appears in the countable list \eqref{E:CountableData}, say at some index $i=i_y\in \N$, and  the ``continuous'' data that consists of a point $\theta=\theta_y \in \Theta_i$ as in \eqref{E:ContinuousData}, 
for which the following estimate holds 
$$
\Phi_i(y,n)\geq \Psi_i(\theta,n)+e(y,n), \qquad \text{for every } n\in \N,
$$ 
where $\Phi_i$ is as in \eqref{E:Phii},  $\Psi_i$ is as in \eqref{E:Psii}, and   $e(y,n)$ is a nullsequence with values in~$[-1,1]$.
	
\smallskip

\textit{3.  Borel selection via \cref{L:SelectionAbstract}.} The hypotheses of \cref{L:SelectionAbstract} are satisfied for 
the standard Borel probability space $(Y,\nu)$, the standard Borel spaces $(\Theta_j)_{j\in\N}$, the Borel maps $(\Phi_j,\Psi_j)_{j\in \N}$, and the previously defined set $Y'$. We deduce that there 
exist disjoint Borel sets $Y_0,Y_1,\ldots$, some of which may be empty, such that 
$$
Y=\bigcup_{i\geq 0} Y_i, \qquad  \nu(Y_0)=0, 
$$
and for $i\in \N$ there exist  Borel maps $\theta_i\colon Y_i\to \Theta_{i}$,  and Borel varying  nullsequences $e_i\colon Y_i\times \N\to [-1,0]$, such that 
\[
	\Phi_i(y,n)\geq \Psi_i(\theta_i(y),n)+e_i(y,n), \qquad \text{for every } n\in \N \text{ and }  y\in Y_i.
\]
By \eqref{E:Phii} the left hand side of this estimate coincides with the one in \eqref{E:muy}.
It remains to unpack the definitions in order to realize that 
$\Psi_i(\theta_i(y),n)$ too has the asserted form. Note that for $i\in \N$ the Borel map $\theta_i$ has the form  
$$
\theta_i(y)=(g_i(y,\cdot), b_{i,1}(y),\ldots, b_{i,k_i}(y))\in C^+_i\times H^{k_i}_i, \qquad y\in Y_i, 
$$
hence $g_i\colon Y_i\times H_i/\Lambda_i\to [0,1]$ is Borel, for every $y\in Y_i$ the function $g_i(y,\cdot)$ is continuous and $g_i(y,1_{H_i/\Lambda_i})>0$ for all $y\in Y_i$, the maps $b_{i,1},\ldots, b_{i,k_i}\colon Y_i\to H_i$ are Borel, and 
$$ 
\Psi_i(\theta_i(y),n)=g_i(y,b_{i,1}(y)^{r_{i,1}(n)}\cdots b_{i,k_i}(y)^{r_{i,k_i}(n)}1_{H_i/\Lambda_i}), \qquad y\in Y_i, 
$$
which coincides with the expression that appears on the right hand side of 
\eqref{E:muy}. Lastly, the joint intersectivity assertion of the polynomials follows from \eqref{E:ji}. 



\section{Compact extensions in ergodic $\CZ_s$-systems}
\label{sec:iso_extension_ergodic}

The goal of this section is to prove \cref{thm:beautiful_main}
under the extra assumption of ergodicity.
We begin by giving an algebraic characterization of compact extensions
between ergodic nilsystems.
By combining a delicate measure-theoretic argument
with properties of rationally equivalent lattices, we then extend the result
to ergodic pro-nilsystems.


\subsection{Compact extensions of ergodic nilsystems}
A nilsystem~$(X,T)$ can be represented as a pair~$(G/\Gamma,\tau)$\footnote{
To be precise, this should be a triple~$(G,\Gamma,\tau)$.
}
with~$G/\Gamma$ a nilmanifold and~$\tau\in G$ an element acting
on~$G/\Gamma$ by left multiplication. Such a representation is not unique.
We think of~$(X,T)$ as an abstract system
up to measurable isomorphism, and of $(G/\Gamma,\tau)$ as a
concrete model; 
for a specific task,
some representations may be better than others.
We will often work with \emph{reduced} representations,
as defined in \cref{sec:prel_nilsystem}.

Let~$X$ be an ergodic nilsystem and~$Y$ a factor.
Abstractly, we identify~$Y$ with the corresponding sub-$\sigma$-algebra.
For our purposes (e.g. compactness of the extension),
this sub-$\sigma$-algebra, and not the specific models of~$X$, $Y$ and factor map~$X\to Y$,
is what matters. On the other hand, having a concrete model allows one to turn
measure-theoretic statements into algebraic statements.
The following result is an algebraic characterization of compact extensions
between ergodic nilsystems. It will be the core ingredient in most results
throughout this section.

\begin{proposition} \label{thm:iso_central}
Let~$(G/\Gamma,\tau)$ be an ergodic nilsystem represented in reduced form,
and let~$L\subseteq  G$ be a closed
subgroup containing~$\Gamma$. Then the following facts are equivalent:
\begin{enumerate} [(i)]
\item The extension~$G/\Gamma\to G/L$ is compact;
\item The identity component~$L^0$ is contained in the center~$Z(G)$.
\end{enumerate}
\end{proposition}

Before proving this result, we present two examples.
The first example illustrates the proof, and shows that it is unavoidable
to consider relatively compact
functions of rank greater than~$1$.

\begin{example} \label{ex:central_higher}
Let~$G=\begin{pmatrix} 1 & \Z & \R \\ 0 & 1 & \R \\ 0 & 0 & 1 \end{pmatrix}$,
$\Gamma=\begin{pmatrix} 1 & \Z & \Z \\ 0 & 1 & \Z \\ 0 & 0 & 1 \end{pmatrix}$,
$L=\begin{pmatrix} 1 & \Z & \R \\ 0 & 1 & \frac{1}{2}\Z \\ 0 & 0 & 1 \end{pmatrix}$,
$\tau=\begin{pmatrix} 1 & 1 & \beta \\ 0 & 1 & \alpha \\ 0 & 0 & 1 \end{pmatrix}$,
with~$1,\alpha,\beta$ rationally independent.
Then~$L^0=Z(G)$, thus by \cref{thm:iso_central} the
extension~$X=G/\Gamma\to Y=G/L$ is compact. It is instructive to prove this fact
directly. Up to isomorphism, we have~$X=\R^2/\Z^2$,  with
transformation~$T(x,y)= (x+\alpha,y+\beta+x) \ \mod 1$, and~$Y=\R/\frac{1}{2}\Z$, with 
the factor~$x\mapsto x+\alpha \ \mod \frac{1}{2}$.
The factor~$Y$ is spanned by functions of the form~$x\mapsto e(2nx)$ with~$n\in \Z$.
The function~$f(x,y):=e(x)$ is rank-1 relatively compact,
that is,~$T^n f \in L^\infty(Y) f$ for all~$n\in \Z$ (in fact~$T^n f \in \C f$).
The function~$f(x,y)=e(y)$ is rank-2 relatively compact,
that is~$T^n f \in L^\infty(Y) f + L^\infty(Y) Tf$ for all~$n\in \Z$. More generally,
the system is spanned by functions of the form~$(x,y)\mapsto e(ax+by)$ with~$a,b\in \Z$,
which are either rank-1 or rank-2 relatively compact
depending on the parity of $b$.
\end{example}

The second example is an instance of an ``abelian'' extension that is not ``central'', in particular not compact.
\begin{example}
Let~$G=\begin{pmatrix}
	1 & \R & \R & \R \\ 0 & 1 & \R & \R \\ 0 & 0 & 1 & \R \\ 0 & 0 & 0 & 1
	\end{pmatrix}$,
$\Gamma=\begin{pmatrix}
	1 & \Z & \Z & \Z \\ 0 & 1 & \Z & \Z \\ 0 & 0 & 1 & \Z \\ 0 & 0 & 0 & 1
	\end{pmatrix}$,
$L=\begin{pmatrix}
	1 & \Z & \R & \R \\ 0 & 1 & \Z & \R \\ 0 & 0 & 1 & \Z \\ 0 & 0 & 0 & 1
	\end{pmatrix}$.
Consider the extension~$G/\Gamma\to G/L$.
Note that~$L^0=[G,G]$ is an abelian group.
Moreover, $L=[G,G]\Gamma$.
Therefore, fibers of~$G/\Gamma\to G/L$ are
isomorphic to the abelian group~$L^0/(L^0\cap \Gamma)\cong \T^3$. However,
as  can be computed explicitly, the action of~$G$ on the fibers is not compact:
it is an affine unipotent transformation with non-trivial unipotent part.
This is due to the fact that~$L^0=[G,G]\not\subseteq  [G,[G,G]]=Z(G)$.
\end{example}

\begin{remark} [Nilcharacters] \label{rem:nilcharacter}
For~$f\in L^2(X)$ and~$\chi$ a character of a rational central subgroup~$Z$,
define~$T:=Z/(Z\cap \Gamma)$. Then~$T$ is a torus.
Define
\[
    f_\chi(x) := \int_T f(tx) \con{\chi(t)} \d \m_T(t).
\]
Then~$f_\chi(zx)=\chi(z)f_\chi(x)$ for all~$z\in Z$ and~$x\in X$, that is,
the function~$f_\chi$ is a \emph{nilcharacter}, with character~$\chi$.
A function~$f$
is $Z$-invariant if and only if~$f_\chi=0$ for all non-trivial characters.

\end{remark}

\begin{proof} [Proof of~\cref{thm:iso_central}]
Throughout, $[a,b]:=aba^{-1}b^{-1}$.
By \cref{sec:prel_nilsystem}, $L^0$ is a normal rational subgroup of~$G$
and~$L/(L^0\Gamma)$ is finite.
Since the intermediate extension~$G/(L^0\Gamma)\to G/L$ is finite-to-one,
the extension~$G/\Gamma\to G/L$ is compact
if and only if the extension~$G/\Gamma \to G/(L^0 \Gamma)$ is compact.
Since~$(L^0\Gamma)^0=L^0$, it is then enough to assume that~$L$ is a connected normal rational subgroup
of~$G$ and prove the equivalence between
\begin{enumerate} [(i')]
\item The extension~$G/\Gamma\to G/(L\Gamma)$ is compact;
\item The subgroup~$L$ is contained in the center~$Z(G)$.
\end{enumerate}

$(i') \Leftarrow (ii')$.
Suppose that~$L\subseteq  Z(G)$.
We will prove that~$\pi \colon  X=G/\Gamma\to Y=G/(L\Gamma)$ is compact.
Since~$L$ is central, $L\cap\Gamma$ is a normal subgroup of~$G$ contained in~$\Gamma$,
thus~$L\cap \Gamma=\{1_G\}$ and~$L$ is a torus.
Let~$F$ be a continuous nilcharacter of~$L$,
with character~$\chi$.
Fix a measurable section~$\sigma\colon G/(L\Gamma) \to G$
(see Section~\ref{prel:sections} for details on measurable sections).
Any~$x\in X$
can be uniquely factorized as~$x=\ell(x) \sigma(y) \Gamma$ with~$y=\pi(x)$ and~$\ell(x)\in L$,
and~$\ell$ is Borel (uniqueness follows from~$L$ central and~$L\cap\Gamma=\{1_G\}$).
Since~$L\subseteq  Z(G)$, we have
\[
	F(\tau^n x) = F(\tau^n \ell(x) \sigma(y)\Gamma) = \chi(\ell(x)) F(\tau^n \sigma(y)\Gamma).
\]
Thus~$T^nF = u\cdot(\phi_n\circ\pi)$, where~$u:=\chi\circ \ell$ has modulus one
and~$\phi_n(y):=F(\tau^n \sigma(y)\Gamma)$ is bounded and Borel.
Moreover, $\tau x = \ell(x)\,\ell(\tau\sigma(y)\Gamma)\,\sigma(\pi(\tau x))\Gamma$,
so that~$u(\tau x)=c(y)\,u(x)$ with~$c(y):=\chi(\ell(\tau\sigma(y)\Gamma))$.
Hence the rank-one module~$L^\infty(Y)u$ is invariant and contains the orbit of~$F$.
Since continuous nilcharacters span a dense subspace of~$L^2(X)$
(see remark above), it follows
from characterization~(iii) in \cref{sec:compact_extensions}
that the extension~$X\to Y$ is compact.

$(i') \Rightarrow (ii')$. Let~$L$ be a connected normal rational
subgroup.
Suppose that the extension~$G/\Gamma \to G/(L\Gamma)$ is compact.
We will prove that~$L\subseteq  Z(G)$. For convenience, we break down the proof in a
series of claims.

\emph{Claim: Without loss of generality~$[G,[L,G]]=1$.}
We obtain this reduction by a ``minimal counterexample''.
Suppose that the extension~$X\to L\backslash X$ is compact but~$[L,G]\neq 1$.
The group~$[G,[L,G]]$ is normal, and it is rational by
the argument in \cite[Chapter 10, Proposition 22]{host2018nilpotent}.
The intermediate extension~$[G,[L,G]]\backslash X \to L\backslash X$ is also compact,
and it can be written as~$X'\to L'\backslash X'$ after replacing~$G,G^0,\Gamma,L,\tau$
by their images~$G',(G^0)',\Gamma',L',\tau'$ in the quotient~$G/[G,[L,G]]$.
Clearly~$[G',[L',G']]=1$. It remains to show that~$X'\to L'\backslash X'$
is a counterexample.
The property~$[L',G']\neq 1$ is preserved:
indeed~$[L',G']=1$ in~$G'$ would imply~$[L,G]\subseteq  [G,[L,G]]$ in~$G$, and
iterating this inclusion along the lower central series of the nilpotent group~$G$ finitely many times
would give~$[L,G]=1$, contrary to our hypothesis.
The representation~$(G'/\Gamma',\tau')$ might be non-reduced:
while the group~$G'$ is certainly generated by~$(G^0)'$ and~$\tau'$,
in general~$\core(\Gamma'):=\bigcap_{g\in G'} g\Gamma' g^{-1}$
might be non-trivial. However, it is enough to replace~$G'$
by~$G'/\core(\Gamma')$, which does not change the systems~$X'$ and~$L'\backslash X'$,
and note that the previous properties
are unaffected. In particular~$[L',G']\neq 1$ survives this quotient
since~$[L',G']$ is connected.

\emph{Claim: Relative eigenfunctions are~$[\tau,L]$-invariant.}
Consider any finite-rank invariant $L^\infty(Y)$-module.
Since~$Y$ is ergodic,
we choose
a relative orthonormal basis~$h=(h_1,\ldots,h_r)$,
see \cref{sec:compact_extensions} in preliminaries.
We have
\begin{equation} \label{eq:unitary}
	h(\tau x) = \Lambda(y) h(x), \qquad h_i(\tau x) = \sum_{j=1}^r \Lambda_{ij}(y) h_j(x)
	\quad \text{for } i=1,\ldots,r,
\end{equation}
where~$y=\pi(x)$
and~$\Lambda(y)=(\Lambda_{i,j}(y))\in \C^{r\times r}$
is a unitary matrix depending measurably on~$y$.
Since~$[\tau,L]\subseteq  Z(G)$, the map~$L\to [\tau,L]$ sending~$g\mapsto [\tau,g]$
is a group homomorphism and~$[\tau,L]$ is a group.
The group~$[\tau,L]$ need not be closed, so we work with the larger group~$Z:=[L,G]$:
by the previous claim it is central, it is connected, it contains~$[\tau,L]$, it is contained in~$L$
(since $L$ is normal), and it is closed and rational by the argument
of~\cite[Chapter 10, Proposition 22]{host2018nilpotent} used in the previous claim.
Let~$T:=Z/(Z\cap\Gamma)$, a torus, and
let~$\chi$ be a character of~$T$ that correlates with the finite-rank module,
that is, for some base element~$h_i$ the nilcharacter~$h_{i,\chi}$
is non-zero.
We have to show that~$\chi$ is trivial on~$[\tau,L]$.
We claim that, since~$Z\subseteq  L$, equation~\eqref{eq:unitary}
holds if we replace each~$h_i$ by~$h_{i,\chi}$.
Indeed, every~$t\in T$ commutes with~$\tau$ (since~$Z$ is central)
and satisfies~$\pi(tx)=\pi(x)$ (since~$Z\subseteq  L$ and~$L$ is normal), thus
\begin{align*}
	h_{i,\chi}(\tau x)
	& = \int_T h_i(\tau t x) \con{\chi(t)} \d \m_T(t)
	= \int_T \sum_{j=1}^r \Lambda_{ij}(\pi(t x)) h_j(t x) \con{\chi(t)} \d \m_T(t) \\
	& = \int_T \sum_{j=1}^r \Lambda_{ij}(y) h_j(t x) \con{\chi(t)} \d \m_T(t)
	= \sum_{j=1}^r \Lambda_{ij}(y) h_{j,\chi}(x).
\end{align*}

To simplify notation, from now on we write~$h_i$
instead of~$h_{i,\chi}$. To summarize our current data, we have
\begin{equation} \label{eq:iso_mid_hyp}
	h(\tau x) = \Lambda(y) h(x), \qquad \exists i, \ h_i\neq 0,
	\qquad \forall i, \forall z\in [\tau,L], \ h_i(zx)=\chi(z)h_i(x),
\end{equation}
and our goal is to prove~$\chi=1$ on~$[\tau,L]$.

Let~$g\in L$ be arbitrary. Consider the function
\[
	H_g(x) := \inp{h(gx)}{h(x)}_{\C^r} = \sum_{i=1}^r h_i(gx) \con{h_i(x)}.
\]
From~\eqref{eq:iso_mid_hyp}, $g\in L$, and~$g\tau=[g,\tau]\tau g$ with~$[g,\tau]\in[\tau,L]$,
we obtain the following eigenvalue identity:
\begin{align*}
	H_g(\tau x) & = \inp{h(g\tau x)}{h(\tau x)}_{\C^r}
		= \chi([g,\tau]) \inp{h(\tau g x)}{h(\tau x)}_{\C^r}
		\\ & = \chi([g,\tau]) \inp{\Lambda(\pi(g x)) h(g x)}{\Lambda(y) h(x)}_{\C^r}
		\\ & = \chi([g,\tau]) \inp{\Lambda(y) h(g x)}{\Lambda(y) h(x)}_{\C^r}
		= \chi([g,\tau]) H_g(x).
\end{align*}
In particular~$\abs{H_g}$ is $\tau$-invariant, hence constant,
in particular~$H_g\in L^\infty(X)$.
If~$H_g\neq 0$, then~$\chi([g,\tau])$ is an eigenvalue of~$\tau$.
We claim that the map~$L\to L^1(X)$ associating~$g\mapsto H_g$ is
continuous. Since
$\norm{H_g-H_{g'}}_{L^1}\leq \sum_{i=1}^r\norm{h_i(g\,\cdot)-h_i(g'\,\cdot)}_{L^2}\norm{h_i}_{L^2}$,
it is enough to fix~$f\in L^2(X)$ and
prove that the map~$G\to L^2(X)$
associating~$g\mapsto (x\mapsto f(gx))$ is continuous.
By density, it is enough to prove this for~$f$ continuous.
Let~$d_G$ be a right-invariant metric on~$G$ inducing its topology
(for instance as in~\cite[Definition 2.2]{green2012quantitative}),
and let~$d_X(g\Gamma,g'\Gamma):=\inf_{\gamma\in\Gamma}d_G(g,g'\gamma)$.
Since~$X$ is compact, then~$f$ is uniformly continuous with respect to~$d_X$.
Therefore, it is enough to prove that~$d_X(g_1 g\Gamma, g_2 g \Gamma)$ is small,
uniformly in~$g\in G$, if~$d_G(g_1,g_2)$ is small.
This follows from
\[
	d_X(g_1 g\Gamma, g_2 g \Gamma)
	= \inf_{\gamma\in \Gamma} d_G(g_1 g, g_2 g \gamma)
	\leq d_G(g_1 g, g_2 g) = d_G(g_1, g_2),
\]
where we used the fact that the metric~$d_G$ is right-invariant.
It follows that~$g\mapsto H_g$ is continuous.

By~\eqref{eq:iso_mid_hyp},
the function~$H_{1_L}(x)=\inp{h(x)}{h(x)}_{\C^r}=\sum \abs{h_i(x)}^2$
is $\tau$-invariant, hence a non-zero constant.
Therefore, for all~$g$ in some connected open neighborhood~$U$ of~$1_L$ in~$L$,
the function~$H_g$ is not the constant zero, and thus~$\chi([g,\tau])$ is an eigenvalue.
That is, the continuous map~$g\mapsto \chi([g,\tau])$
maps the connected set~$U$
into the countable set of eigenvalues, hence it is constant, equal to~$1$, on~$U$.
Since~$g\mapsto \chi([g,\tau])$ is a group homomorphism and~$U$ generates the connected
group~$L$, this implies~$\chi([L,\tau])=1$.
Since~$h_i=\sum_{\chi}h_{i,\chi}$, every~$h_i$ is $[\tau,L]$-invariant.

\emph{Claim: Compact extension implies $[\tau,L]=1$.}
We proved that relative eigenfunctions are~$[\tau,L]$-invariant, and
since the extension is compact, such eigenfunctions span a dense subspace of~$L^2(X)$ with coefficients in~$L^\infty(Y)$.
Since~$[\tau,L]\subseteq  L$, every element~$z\in [\tau,L]$ acts
trivially on~$L^\infty(Y)$.
It follows that every~$z\in[\tau,L]$ acts trivially on~$L^2(X)$,
and thus on~$G/\Gamma$, that is~$z\in\core(\Gamma)$.
Since~$\Gamma$ contains no non-trivial normal subgroups of~$G$,
this implies that~$[\tau,L]=1_G$.

\emph{Claim: The subgroup~$L$ is central.}
By the previous claims, $[L,G]\subseteq  Z(G)$ and~$[\tau,L]=1$.
Hence, for every~$h\in L$ and~$g\in G$, the maps~$G\to Z(G)$, $g'\mapsto[h,g']$,
and~$L\to Z(G)$, $h'\mapsto[h',g]$, are continuous group homomorphisms.

We first show that~$L\cap\Gamma\subseteq  Z(G)$. Let~$\delta\in L\cap \Gamma$.
For every~$\gamma\in\Gamma$ we have~$[\delta,\gamma]\in \Gamma\cap Z(G)$, which is a normal
subgroup of~$G$ contained in~$\Gamma$, hence trivial.
Therefore the map~$g\mapsto[\delta,g]$ is right-$\Gamma$-invariant and induces
a continuous map~$f\colon X\to Z(G)$, $f(g\Gamma):=[\delta,g]$.
Since~$[\delta,\tau]=1$, the map~$f$ is $\tau$-invariant. By ergodicity it is almost
everywhere constant, and since it is continuous and~$\m_X$ has full support,
it is constant, equal to~$f(1_X)=1_G$. Thus~$\delta\in Z(G)$.

We now show that~$L$ commutes with~$G^0$.
Since~$L$ is rational, $L/(L\cap\Gamma)$ is compact, so there is a compact set~$K\subseteq  L$
with~$L=K\cdot(L\cap\Gamma)$.
Let~$U$ be a neighborhood of~$1_G$ containing no non-trivial subgroup of~$G$.
Since~$(h,g)\mapsto[h,g]$ is continuous, $K$ is compact, and~$[K,1_G]=\{1_G\}$,
there is a neighborhood~$V$ of~$1_G$ such that~$[K,V]\subseteq  U$.
Fix~$g\in V$. The homomorphism~$L\to Z(G)$, $h\mapsto[h,g]$, is trivial on~$L\cap\Gamma$
by the previous paragraph, so its image is~$[K,g]$. This is a subgroup of~$G$ contained in~$U$,
hence trivial. Therefore~$L$ commutes with~$V$, and thus with the group generated by~$V\cap G^0$,
which is~$G^0$ since~$G^0$ is connected.

Finally, since~$G$ is generated by~$G^0$ and~$\tau$ (by assumption) and~$[\tau,L]=1$,
we conclude that~$L\subseteq  Z(G)$.
\end{proof}

As a corollary, we obtain  $\nilbor$ returns for  compact extensions
of ergodic nilsystems.

\begin{proposition} \label{prop:key_lem_nilsystem}
Let~$X$ be an ergodic nilsystem and~$Y$ a factor such that the extension~$X\to Y$
is compact. Then there is a countable set of functions~$\CF\subseteq  L^\infty(X)$,
dense in~$L^2(X)$, with the following property: for every~$f\in \CF$ and every~$\ve>0$,
there is a nilmanifold~$H/\Lambda$, a continuous function~$g$ on~$H/\Lambda$,
an element~$a\in H$,
a measurable map~$\alpha \colon Y\to H$, and~$\ve'>0$ such that
\[
	\{n \in \Z\colon\norm{T^n f - f}_{L^2(\mu_y)} \leq \ve\}
	\supseteq
	\{n \in \Z\colon\abs{g(a^n \alpha(y)\Lambda) - g(\alpha(y)\Lambda)} \leq \ve'\}
\]
for almost every~$y\in Y$, i.e., the set on the left contains a measurable varying $\nilbor$ set.
\end{proposition}

\begin{proof}
For exposition purposes, we first present the argument in the case of a quotient
by a connected subgroup.
Let~$X=(G/\Gamma,\tau)$ be in reduced form and~$Y=G/(L\Gamma)$
for some connected rational subgroup~$L\subseteq  G$.
By \cref{thm:iso_central} we have~$L\subseteq  Z(G)$.
As in the proof of \cref{thm:iso_central}, any element~$x\in X$
factorizes as~$x=\ell(x)\sigma(y)\Gamma$, where~$y=\pi(x)$,
for some measurable maps~$\sigma\colon Y\to G$ and~$\ell\colon X\to L$.
Let~$F$ be a nilcharacter, say with (central) character~$\chi$, and~$\ve>0$.
Since~$L\subseteq  Z(G)$, we have
\begin{equation} \label{eq:key_lem_nilsystem_F}
	\abs{F(\tau^n x) - F(x)}
	= \abs{\chi(\ell(x))} \abs{F(\tau^n \sigma(y)\Gamma) - F(\sigma(y)\Gamma)}
	= \abs{F(\tau^n \sigma(y)\Gamma)-F(\sigma(y)\Gamma)}.
\end{equation}
It follows that
\[
	\{n\in \Z\colon\norm{T^n F - F}_{L^2(\mu_y)} \leq \ve \}
	= \{n\in \Z\colon\abs{F(\tau^n \sigma(y)\Gamma)-F(\sigma(y)\Gamma)} \leq \ve \}.
\]
In particular, every nilcharacter~$F$ has the desired $\nilbor$ property.
Since the intersection of
finitely many $\nilbor$ sets is a $\nilbor$ set, any finite linear combination of
nilcharacters~$\sum_{i=1}^m c_i F_i$ has the desired property. 
Choosing the coefficients from a fixed countable dense subfield of $\C$, we obtain a countable set of finite linear combinations of nilcharacters that is dense in $L^2(X)$. This proves the claim in the connected case.

In the general case we have~$X=G/\Gamma$ and~$Y=G/L$ with~$L$
a closed subgroup containing~$\Gamma$ and such that~$L^0\subseteq  Z(G)$.
A small change in the previous argument is necessary to handle the
finite extension~$G/(L^0\Gamma) \to G/L$.
Let~$g_1,\ldots,g_k\in L$
be such that~$L=\bigcup_{i=1}^k g_i L^0 \Gamma$.
An element~$x\in X$ now factorizes as~$x=\ell(x)\sigma(y) g_i \Gamma$,
where~$y=\pi(x)$, $\ell(x)\in L^0$, and the index ~$i=i_x\in \{1,\ldots,k\}$.
Now~\eqref{eq:key_lem_nilsystem_F} takes the form
\[
	\abs{F(\tau^n x) - F(x)}
	= \abs{F(\tau^n \sigma(y) g_i\Gamma)-F(\sigma(y)g_i\Gamma)}
	\leq \sum_{j=1}^k \abs{F(\tau^n \sigma(y) g_j\Gamma)-F(\sigma(y)g_j\Gamma)}.
\]
In other words, to remove the dependence of~$g_i$ on~$x$,
we sum over all finitely many possibilities.
Since the intersection of finitely many $\nilbor$ sets is~$\nilbor$,
the same conclusion holds as before.
\end{proof}

In the remainder of this section, we extend the previous result to the case where $X$ is an inverse limit of ergodic nilsystems.


\subsection{The good intersection property} \label{sec:good}
We take a brief detour to discuss relative independence for $\sigma$-algebras,
and define a particular instance that will be useful throughout the paper.
As usual, we assume all probability spaces standard
and, unless the specific factor map is important, identify factors with sub-$\sigma$-algebras.
The results in \cref{sec:good} will be used not only
for our current objective, \cref{thm:beautiful_main},
but throughout the whole paper.

\begin{definition}
We say that the factors $Y$ and~$Z$  of $X$ have the \emph{good intersection property}
if they are relatively independent
with respect to~$Y\wedge Z$, that is, if any of the following equivalent
properties hold:
\begin{enumerate} [(i)]
\item $\E(fg\mid Y\wedge Z)=\E(f\mid Y\wedge Z)\E(g\mid Y\wedge Z)$
	for all bounded $Y$-measurable functions~$f$ and $Z$-measurable functions~$g$;
\item $\E(\E(f\mid Y)\mid Z)=\E(f\mid Y\wedge Z)$ for all bounded measurable functions~$f$;
\item $\E(f\mid Y)$ is $Z$-measurable for all bounded $Z$-measurable functions~$f$.\end{enumerate}
\end{definition}

Given a system~$(X,\mu,T)$, the invariant factor~$I(X)$ is ``universal'' with respect
to the good intersection property, in the following sense.

\begin{lemma}
Let $(X,\mu,T)$ be a system, $Y$~be an arbitrary factor, and $I(X)$ the invariant factor.
Then~$I(X)$ and~$Y$ have the good intersection property.
\end{lemma}
\begin{proof}
It follows immediately from the Mean Ergodic Theorem: if~$f$ is $Y$-measurable,
then~$\E(f\mid I(X))$ is $Y$-measurable.
\end{proof}

Similarly, $\CZ_s$-factors are universal with respect to the good intersection property.

\begin{lemma}\label{L:ZsGood}
Let $(X,\mu,T)$ be a system, $Y$ be an arbitrary factor,
	and~$\CZ_s(X)$ the usual Host-Kra factor.
	Then~$\CZ_s(X)$ and~$Y$ have the good intersection property.
\end{lemma}

\begin{proof}
	Let $f\in L^\infty(\CZ_s(X))$. In order to show that $\E(f\mid Y)\in L^\infty(\CZ_s(X))$ it suffices to show that $\E(f\mid Y)$ is orthogonal to all $\nnorm{\cdot}_{s+1}$-uniform functions, see
    \cref{sec:prel_HK} for details. So let $h\in L^\infty(\mu)$ with $\nnorm{h}_{s+1}=0$. Since 
	$$
	\langle \E(f\mid Y), h\rangle  = \langle  f, \E(h\mid Y)\rangle,
	$$
	it suffices to show that $\nnorm{\E(h\mid Y)}_{s+1}=0$; indeed, since $f\in L^\infty(\CZ_s(X))$ this would then  imply that $\, f\bot \, \E(h\mid Y)$.
	So it remains to establish the implication 
	$$
	\nnorm{h}_{s+1}=0\implies  \nnorm{\E(h\mid Y)}_{s+1}=0.
	$$
	Let $h_Y:=\E(h\mid Y)$ and 
    $\mathcal{D}_{s+1}(h_Y)\in L^\infty(\CZ_s(X))$ be the order  $s+1$  dual function, see \cref{sec:prel_HK}. We have (\cref{sec:prel_HK})
		$$
	\nnorm{\E(h\mid Y)}^{2^{s+1}}_{s+1}= \langle h_Y,  \mathcal{D}_{s+1}(h_Y)\rangle  =\langle  h, \mathcal{D}_{s+1}(h_Y)\rangle,
	$$
where the last equality holds since $h_Y$ and hence $\mathcal{D}_{s+1}(h_Y)$ are $Y$-measurable.
	Since $\nnorm{h}_{s+1}=0$ we have $h\bot L^\infty(\CZ_s(X))$, hence  the last inner product is $0$, which gives  $\nnorm{\E(h\mid Y)}_{s+1}=0$, as required.
\end{proof}

Next, we show that compact  extensions are well
behaved with respect to the good intersection property.

\begin{lemma}\label{L:compactfactor}
Let $(X,\mu,T)$ be a system and $X\to Y$ be a compact extension.
Let~$Z$ be a factor such that~$Y$ and~$Z$ have the good intersection property.
Then the extension $Z\to Z\wedge Y$ is compact.
\end{lemma}

\begin{proof}	
We will use the following well-known criterion (see preliminaries, \cref{sec:compact_extensions}):
The extension $X\to W$ is compact if and only if 
for every $f\in L^\infty(\mu)$ that is not a.e.~$0$ we have 
\begin{equation}\label{E:positivity}
	\limsup_{N\to \infty}\E_{n\in [N]}\int |\E(T^nf\cdot \bar{f}\mid W)|^2\, d\mu>0.
\end{equation}
 	Let $f\in L^\infty(\mu)$ be $Z$-measurable and not a.e.~$0$. 	
	It suffices to verify \eqref{E:positivity} with  $Z\wedge Y$ in place of $W$ and $\mu_Z$ in place of $\mu$; note though that having $\mu$ or $\mu_{Z}$ makes no difference, as the integrand  is $Z$-measurable. Since the extension $X\to Y$ is compact and $f$ is not a.e. $0$ we have 
	\begin{equation}\label{E:positivity'}
		\limsup_{N\to \infty}\E_{n\in [N]}\int |\E(T^nf\cdot \bar{f}\mid Y)|^2\, d\mu>0.
	\end{equation}
	 Since $f$ is $Z$-measurable and $Z$ is a factor, $T^nf\cdot \bar{f}$ is also $Z$-measurable, and since $Y,Z$ have the good intersection property,   we have 
     $$\E(T^nf\cdot \bar{f}\mid Y)=\E\big(T^nf\cdot \bar{f}|(Z\wedge Y)\big)
     $$
     for every $n\in \N$. 
	Hence,~\eqref{E:positivity'} gives 
	$$
		\limsup_{N\to \infty}\E_{n\in [N]}\int |\E\big(T^nf\cdot \bar{f}\mid(Z\wedge Y)\big)|^2\, d\mu>0,
	 $$	
	 as required.
\end{proof}



	\subsection{Rationally equivalent lattices}
	
	We take a second brief detour to discuss rationally equivalent ``lattices''
	(more precisely, rationally equivalent discrete cocompact subgroups).
	Let~$\Gamma$ and~$\Gamma'$ be discrete cocompact subgroups of a nilpotent
	Lie group~$G$. Following~\cite[Section~1]{leibman2006rational},
	we call~$\Gamma$ and~$\Gamma'$ \emph{rationally equivalent}
	if any of the following equivalent properties hold:
	the group $\Gamma\cap \Gamma'$ has finite index in~$\Gamma$;
	the group $\Gamma\cap \Gamma'$ has finite index in~$\Gamma'$;
	every point in~$\Gamma'$ is $\Gamma$-rational;
	every point in~$\Gamma$ is $\Gamma'$-rational;
	the set of $\Gamma$-rational points and the set of $\Gamma'$-rational points
	coincide.
	
	\begin{lemma} \label{lemma:commensurable}
		Let~$g\in G$ be a $\Gamma$-rational. Then~$\Gamma$ and~$g\Gamma g^{-1}$
		are rationally equivalent.
	\end{lemma}
	
	\begin{proof}
		By~\cite[Corollary 1.15]{leibman2006rational}, the group~$\Gamma':=g\Gamma g^{-1}$
		is $\Gamma$-rational. By~\cite[Lemma 1.11]{leibman2006rational},
		the group~$\Gamma \cap \Gamma'$ is cocompact in~$\Gamma'$.
		Since~$\Gamma'$ is discrete, the set~$\Gamma'/(\Gamma \cap \Gamma')$ is finite.
	\end{proof}
	
	
	\begin{lemma} \label{lem:rational_generated}
		Let~$g_1,\ldots,g_m$ be $\Gamma$-rational.
		Then the group~$\Gamma'$ generated by~$\Gamma$ and $g_1,\ldots,g_m$ is discrete
		and cocompact.
		In particular~$\Gamma$ and~$\Gamma'$ are rationally equivalent.
	\end{lemma}
	
	\begin{proof} 
		By~\cite[Chapter 10, Lemma 24]{host2018nilpotent}, the group~$\Gamma$
		is finitely generated.
		Therefore, the group~$\Gamma'$ is also finitely generated and, in particular,
		it is generated by finitely many $\Gamma$-rational elements.
		By~\cite[Lemma 1.14]{zorin2013ergodic},\footnote{
			This lemma states that, in \textit{any} nilpotent group, if~$H,H'$
			are subgroups with~$H\subseteq  H'$ and~$H'$ is generated by finitely
			many $H$-rational elements, then~$H$ has finite index in~$H'$.
		}
		the group~$\Gamma$ has finite index in~$\Gamma'$. It follows that~$\Gamma'$ is discrete.
		It is clear that~$\Gamma'$ is cocompact since~$\Gamma\subseteq  \Gamma'$.
	\end{proof}

The following lemma presents a useful factorization of rational subgroups.

\begin{lemma} \label{fact:L_factorization}
Let~$G$ be a nilpotent Lie group and~$\Gamma$ a discrete cocompact subgroup.
Let~$L$ be a $\Gamma$-rational subgroup. Then
\[
	L=L^0\cdot \Lambda
\]
where~$\Lambda$ is a discrete group containing~$L\cap \Gamma$ with finite index.
\end{lemma}
\begin{proof}
The quotient~$L/(L^0 \cdot (L\cap \Gamma))$ is both discrete and compact, thus finite.
In particular, there are finitely many elements~$q_1,\ldots,q_k \in L$
such that~$L = \bigsqcup_{i=1}^k q_i \cdot L^0 \cdot (L\cap \Gamma)$.
By \cite[Chapter 10, Lemma 12]{host2018nilpotent}, applied to the group~$L$
and the discrete cocompact subgroup~$L\cap \Gamma$,
there is~$d\geq 1$ such that for every~$q_i$
we have~$(q_i \ell)^d \in L\cap \Gamma$ for some~$\ell \in L^0$.
Up to replacing~$q_i$ by~$q_i \ell$, we may suppose~$q_i^d \in L\cap \Gamma$,
that is~$q_i$ is $(L\cap \Gamma)$-rational.
Then~$\Lambda := \gen{L\cap \Gamma,q_1,\ldots,q_k}$ is generated by~$L\cap \Gamma$
and finitely many $(L\cap \Gamma)$-rational elements,
thus by \cref{lem:rational_generated} it is discrete and contains~$L\cap \Gamma$
with finite index. By construction~$L=L^0 \cdot \Lambda$.
\end{proof}

We will need the following generalization to sequences of rational subgroups.

\begin{lemma} \label{fact:Ln_factorization}
Let~$G$ be a nilpotent Lie group and~$\Gamma$ a discrete cocompact subgroup.
Let~$(L_n)_{n\geq 1}$ be a decreasing sequence of $\Gamma$-rational subgroups.
Then there is a decreasing sequence of discrete
$\Gamma$-rational subgroups~$(\Lambda_n)_{n\geq 1}$
such that
\[
	L_n = L_n^0 \cdot \Lambda_n.
\]
Moreover, for every~$n\geq 1$ the group~$\Lambda_n$ contains~$L_n\cap \Gamma$
with finite index.
Furthermore, there is~$N\geq 1$ and a connected $\Gamma$-rational
subgroup~$Q$ such that~$L_n^0=Q$ for all~$n\geq N$.
\end{lemma}
\begin{proof}
For each~$n\geq 1$, let
$
    L_n:= L_n^0 \cdot \Gamma_n
$
be the factorization given by \cref{fact:L_factorization}. 
The sequence~$n\mapsto L_n^0$ eventually stabilizes. Let~$N$ be a positive
integer and~$Q$ a connected rational subgroup such
that~$L_n^0=Q$ for all~$n\geq N$.
Define
\[
    \Lambda_n :=
    \begin{dcases*}
        \Gamma_N\cap L_n, & if $n>N$ \\
        \Gamma_N, & if $n=N$ \\
        \gen{\Gamma_N,\Gamma_{N-1},\ldots,\Gamma_n}, & if $1\leq n<N$.
    \end{dcases*}
\]
It is clear that the sequence~$n\mapsto \Lambda_n$ is decreasing.

Let us prove that~$\Lambda_n$ is discrete.
For~$n<N$, note that~$\Gamma_m\cap \Gamma_n$ has finite index in~$\Gamma_m$
for~$n\leq m<N$, so~$\gen{\Gamma_N,\Gamma_{N-1},\ldots,\Gamma_n}$
is generated by a rational subgroup of~$\Gamma_n$ and by finitely many
rational elements, thus it is discrete by \cref{lem:rational_generated}.

Let us prove that~$L_n = L_n^0\cdot \Lambda_n$. This is clear for~$n=N$.
Let~$n>N$. Take an element~$\ell \in L_n$. Since~$L_n\subseteq  L_N$,
we have~$\ell=q\gamma$ for some~$q\in Q$ and some~$\gamma\in \Gamma_N$.
Since~$L_n^0=L_N^0=Q$,
we have~$\gamma = q^{-1} \ell \in Q L_n = L_n$, and
thus~$\gamma\in \Gamma_N\cap L_n = \Lambda_n$.
This proves~$L_n\subseteq  L_n^0 \cdot \Lambda_n$.
Conversely, we trivially
have~$L_n^0 \cdot \Lambda_n\subseteq  L_n^0 \cdot \Gamma_n = L_n$.
Now let~$n<N$. On one hand, we have
\[
    L_n^0 \cdot \Lambda_n
    = L_n^0 \cdot \gen{\Gamma_N,\Gamma_{N-1},\ldots,\Gamma_n}
    \subseteq  L_n^0 \cdot \gen{L_N,L_{N-1},\ldots,L_n}
    = L_n.
\]
On the other hand, we have
\[
    L_n = L_n^0 \cdot \Gamma_n \subseteq 
        L_n^0 \cdot \gen{\Gamma_N,\Gamma_{N-1},\ldots,\Gamma_n}
        = L_n^0 \cdot \Lambda_n.
\]

It remains to prove that~$\Lambda_n$ contains~$L_n\cap \Gamma$ with finite index.
Since~$\Gamma_n$ is obtained from \cref{fact:L_factorization},
it contains~$L_n\cap \Gamma$ with finite index. In particular, for~$n<N$,
the group~$\Lambda_n=\gen{\Gamma_N,\Gamma_{N-1},\ldots,\Gamma_n}$,
which contains~$\Gamma_n$ with finite index, also contains~$L_n\cap \Gamma$
with finite index. For~$n\geq N$,
we have~$L_n/(L_n\cap \Gamma)$ compact and
\[
    L_n\cap \Gamma \subseteq  L_n\cap (L_N\cap \Gamma)
        \subseteq  L_n \cap \Gamma_N = \Lambda_n,
\]
thus~$\Lambda_n$ contains~$L_n\cap \Gamma$ with finite index.
\end{proof}

As a corollary, we obtain the following.

\begin{corollary} \label{cor:nilsystems_stabilize}
Let~$X$ be an ergodic nilsystem and~$(X_n)_{n\in \N}$ an ascending
sequence of factors. Then the sequence stabilizes.
\end{corollary}
\begin{proof}
Let~$X=G/\Gamma$. Then~$X_n=G/L_n$ where~$L_n$ is a closed subgroup
with~$\Gamma\subseteq  L_n\subseteq  G$ and~$L_{n+1}\subseteq  L_n$
for all~$n\in \N$. By \cref{fact:Ln_factorization},
we have~$L_n=L_n^0 \Lambda_n$ where~$\Lambda_n$ is a
discrete rational subgroup, containing~$L_n\cap \Gamma$ with finite
index, and~$\Lambda_{n+1}\subseteq  \Lambda_n$ for all~$n\in \N$.
Since by hypothesis~$L_n\cap \Gamma=\Gamma$, we find a decreasing
sequence of positive integers~$[\Lambda_n\colon \Gamma]\in \Z$.
In particular, there is~$N\geq 1$ such that
\[
    [\Lambda_n\colon \Gamma] = [\Lambda_N\colon \Gamma]
\]
for all~$n\geq N$, which implies~$[\Lambda_N\colon \Lambda_n]=1$,
thus~$\Lambda_n=\Lambda_N$ for all~$n\geq N$.
This proves that~$\Lambda_n$ stabilizes. Since~$L_n^0$
stabilizes by \cref{fact:Ln_factorization}, so does~$L_n=L_n^0 \Lambda_n$.
\end{proof}


\subsection{Commuting Subgroups}
The following notion will be used later in connection with the good intersection property.
\begin{definition}
We say that two subgroups~$L$ and~$M$ \emph{commute} if any of the following
equivalent properties hold:
\begin{enumerate} [(i)]
\item $LM = ML$;
\item $LM\subseteq  ML$;
\item the set~$LM$ is a subgroup;
\item $[L,M]\subseteq  LM$.
\end{enumerate}
\end{definition}

\begin{remark}
The condition~$LM=ML$ is strictly weaker than~$\ell M = M\ell$ for~$\ell \in L$.
That is, two subgroups may commute without normalizing each other.
\end{remark}

We will use the following observation about commuting subgroups.

\begin{lemma} \label{lem:generated_commute}
Let~$L$, $M_1$, and~$M_2$ be subgroups. If~$L$ commutes with~$M_1$ and with~$M_2$,
then~$L$ commutes with~$\gen{M_1,M_2}$.
\end{lemma}
\begin{proof}
Define~$M':=\{m\in \gen{M_1,M_2}\colon mL \subseteq  L \gen{M_1,M_2}\}$.
By hypothesis~$M'$ contains~$M_1\cup M_2$. We claim that~$M'$ is closed
under product: for every~$m_1,m_2\in M'$ we have
\[
	m_1 m_2 L \subseteq  m_1 L \gen{M_1,M_2} \subseteq  L \gen{M_1,M_2}.
\]
We conclude that~$M'=\gen{M_1,M_2}$.
\end{proof}


\subsection{Compact extensions in ergodic pro-nilsystems} \label{sec:iso_nil_pro}

By the Ergodic Structure \cref{thm:HK-structure}, every ergodic $\CZ_s$ system
is a pro-nilsystem, that is, it is an inverse limit of ergodic nilsystems.
Our next goal is to extend Proposition~\ref{prop:key_lem_nilsystem}
from ergodic nilsystems to ergodic pro-nilsystems.

Let~$X=\varprojlim_{n} X_n$ be an inverse limit of ergodic nilsystems,
and let~$Y$ be a factor of~$X$.
The factor~$Y$ is also an ergodic pro-nilsystem (\cref{sec:prel_HK}).
Suppose that the extension~$X\to Y$ is compact.
One might hope that
the extensions~$X_n\to X_n\wedge Y$ were compact
(so that Proposition~\ref{prop:key_lem_nilsystem} would apply).
Unfortunately, in general~$X_n\to X_n\wedge Y$ are not compact,
and the obstruction is the finite-to-one nilsystem
extensions discussed in preliminaries, \cref{sec:prel_nilsystem}.
As suggested in Example~\ref{ex:central_higher},
finite-to-one extensions arise from
the presence of elements~$\ell\in L$ with~$\ell^d \in \Gamma$ for some~$d>1$
but~$\ell\notin L^0 \Gamma$. The idea is that the coordinates of such~$\ell$
are rational numbers with non-trivial denominators
in a Malcev basis adapted to~$\Gamma$ and passing through~$L^0$.
Such elements are an obstruction.
However, by a ``killing the denominators'' type of argument,
we are able to slightly enlarge the factors~$X_n$, so that~$X_n$ and~$Y$
satisfy the good intersection property discussed in Section~\ref{sec:good}.
From the good intersection property it follows
that the extension~$X_n \to X_n\wedge Y$ is compact, to which
Proposition~\ref{prop:key_lem_nilsystem} applies.
From this and the fact that the factors~$X_n$ are dense in~$X$,
we will generalize Proposition~\ref{prop:key_lem_nilsystem}
to pro-nilsystems.

\begin{lemma} \label{lem:LM_expected}
Let~$G/\Gamma$ be a nilmanifold. Let~$L,M$ be closed subgroups containing~$\Gamma$.
The following facts are equivalent:
\begin{enumerate} [(i)]
\item The subgroups~$L$ and~$M$ commute;
\item The factors~$G/L$ and~$G/M$ have the good intersection property.
\end{enumerate}
\end{lemma}
\begin{proof}
We recall that,
since~$L$ contains~$\Gamma$, the identity component~$L^0$
is normal in~$G$ by~\cite[Chapter 13, Lemma 2]{host2018nilpotent}.
Let~$f$ be a right-$\Gamma$-invariant bounded measurable function on~$G$. Then
\begin{equation} \label{lem:LM_expected:G/L}
	\E(f\mid G/L)(g) = \int_{\ell \in L/\Gamma} f(g \ell),
\end{equation}
where the integral is with respect to the normalized Haar measure on~$L/\Gamma$.
On the other hand, since~$L^0 \Gamma \subseteq  L^0 (L\cap M) \subseteq  L$
and~$L/(L^0\Gamma)$ is finite, there exist finitely many~$\ell_1,\ldots,\ell_r\in L$
such that
\[
	L = \bigsqcup_{i=1}^r \ell_i L^0 (L\cap M)
\]
where~$\bigsqcup$ denotes a disjoint union. It follows that
\begin{equation} \label{lem:LM_expected:LM}
	LM = \bigsqcup_{i=1}^r \ell_i L^0 M,
\end{equation}
where this last union is disjoint since, for any~$i\neq j$,
whenever~$\ell_j^{-1}\ell_i = \ell^0 m$ for some~$\ell^0\in L^0$
and~$m\in M$, we also have~$m\in L\cap M$.

Let us prove ``$(i)\implies (ii)$''. Let~$f$ be a right-$M$-invariant bounded measurable
function on~$G$. We need to prove that~$\E(f\mid G/L)$ is right-$M$-invariant,
using the fact that~$LM=ML$.
Since~$f$ is right-$(L\cap M)$-invariant,
the average~\eqref{lem:LM_expected:G/L} is the same as
\[
	\E(f\mid G/L)(g) = \int_{\ell \in L/L\cap M} f(g \ell)
		= \frac{1}{r} \sum_{i=1}^r \int_{\ell \in L^0/(L^0\cap M)} f(g \ell_i \ell).
\]
Fix an arbitrary~$m\in M$. From~$LM=ML$, we have~$mLM\subseteq  LM$.
Thus from~\eqref{lem:LM_expected:LM}
and the fact that cosets are either equal or disjoint, it follows that
there is a permutation~$\sigma$ of~$\{1,\ldots,r\}$ (the permutation depends on~$m$)
such that~$m\ell_{\sigma(i)} = \ell_i \ell^0_i m_i$
for some~$\ell^0_i\in L^0$ and~$m_i \in M$.
It follows that 
\begin{align*}
	\E(f\mid G/L)(gm)
	& = \frac{1}{r} \sum_{i=1}^r \int_{\ell \in L^0/(L^0\cap M)} f(g m \ell_{\sigma(i)} \ell) \\
	& = \frac{1}{r} \sum_{i=1}^r
		\int_{\ell \in L^0/(L^0\cap M)}
			f(g \ell_{i} \ell^0_{i} m_{i} \ell) \\
	& = \frac{1}{r} \sum_{i=1}^r
		\int_{\ell \in L^0/(L^0\cap M)}
			f(g \ell_{i} \ell^0_{i} m_{i} \ell m_i^{-1}) \\
	& = \frac{1}{r} \sum_{i=1}^r
		\int_{\ell \in L^0/(L^0\cap M)}
			f(g \ell_{i} \ell)
\end{align*}
where we used the right-$M$-invariance of~$f$ and the fact
that the map~$\ell\mapsto \ell_i^0 m_i \ell m_i^{-1}$ preserves the Haar measure
on~$L^0/(L^0\cap M)$. This proves ``$(i)\implies (ii)$''.

The converse implication ``$(ii)\implies (i)$'' is not hard to prove,
but it is not needed in this paper and is omitted.
\end{proof}

\subsubsection{Killing the denominators}
Let~$L,M$ be closed subgroups of~$G$ containing~$\Gamma$.
The groups~$L$ and~$M$ might not satisfy~$LM=ML$,
and the corresponding factors~$Y=G/L$ and~$Z=G/M$ might not have the good intersection
property.
However, the subgroup~$M^0 \Gamma$,
which has finite index in~$M$,
does satisfy~$L(M^0\Gamma)=(M^0\Gamma)L$, since~$M^0$ is normal and~$\Gamma\subseteq  L$.
Therefore, even if the good intersection property fails for~$Y$ and~$Z$,
it will hold if we replace~$Z=G/M$ by the finite extension~$\tilde Z:=G/(M^0\Gamma)$.
Later we construct a finite extension~$\tilde Z$ that is good
not just for a single factor~$Y$, but simultaneously
for an inverse family of factors~$(Y_n)_{n\geq 1}$.
Unfortunately, the simple choice~$\tilde Z=G/(M^0\Gamma)$ that worked here
for a single factor no longer works, and the argument has to be refined.

\begin{example} \label{ex:denominators}
Let
\[
	G := \begin{pmatrix} 1 & \R & \R \\ 0 & 1 & \R \\ 0 & 0 & 1 \end{pmatrix}, \quad
	M := \begin{pmatrix} 1 & \Z & \Z \\ 0 & 1 & \Z \\ 0 & 0 & 1 \end{pmatrix}.
\]
Let
\[
	\renewcommand \arraystretch{1.2}
	L_n :=
	\begin{pmatrix}
	1 & \frac{2^{n}}{3} \Z & 5 \Z
		\\ 0 & 1 & \frac{15}{2} \Z
		\\ 0 & 0 & 1
	\end{pmatrix}, \quad \text{in particular} \quad
	L_1 :=
	\begin{pmatrix}
	1 & \frac{2}{3} \Z & 5 \Z
		\\ 0 & 1 & \frac{15}{2} \Z
		\\ 0 & 0 & 1
	\end{pmatrix}.
\]
Let~$\tilde M$ be a subgroup of the form
\[
	\tilde M :=
	\begin{pmatrix}
	1 & A \Z & C \Z \\ 0 & 1 & B \Z \\ 0 & 0 & 1
	\end{pmatrix}
\]
for some $A,B,C\in \Z$.
The condition~$\tilde M L_n = L_n \tilde M$ reduces to the following condition on the
``central coordinate'' of the group product:
\[
	A \cdot \frac{15}{2} \Z + C\Z + 5\Z
	= B \cdot \frac{2^n}{3} \Z + C\Z + 5\Z.
\]
The condition is satisfied for all~$n$, for example, by choosing~$(A,B,C)=(2,15,1)$.
\end{example}

Our next goal is to prove the following result.

\begin{proposition} [Killing the Denominators] \label{prop:Znil_Ypro}
Let~$X$ be an ergodic pronilsystem. Let~$Y$ and~$Z$ be factors of~$X$, with~$Z$
a nilfactor. Then there is a nilfactor $\tilde{Z}$ such that $Z\subseteq  \tilde Z \subseteq  X$, $\tilde Z \to Z$ is a finite extension,
and such that~$Y$ and~$\tilde Z$ have the good intersection property.
\end{proposition}

 \cref{prop:Znil_Ypro} follows from the following more detailed  statement.

\begin{proposition} \label{prop:Znil_Ypro_explicit}
Consider an inverse limit of nilmanifolds
\[
	\cdots \to G_3/\Gamma_3
		\to G_2/\Gamma_2 \to G_1/\Gamma_1
\]
with factor maps
induced by surjective continuous group homomorphisms~$\varphi_{n,m}\colon G_n\to G_m$
satisfying~$\varphi_{n,m}(\Gamma_n)\subseteq  \Gamma_m$. Let~$M$ be a closed
subgroup such that~$\Gamma_1 \subseteq  M \subseteq  G_1$. Let~$(L_n)_{n\geq 1}$ be
a sequence of closed subgroups such that~$\Gamma_n\subseteq  L_n \subseteq  G_n$ for all~$n$
and~$L_n \subseteq  \varphi_{n,m}^{-1}(L_m)$ for~$n\geq m$.
Then there exists a $\Gamma_1$-rational subgroup~$\tilde M \subseteq  G_1$ and a positive
integer~$N\geq 1$ satisfying the following properties:
\begin{enumerate} [(i)]
\item $\tilde M \subseteq  M$ with finite index;
\item $\Gamma_n \subseteq  \varphi_{n,1}^{-1}(\tilde M)$ for all~$n\geq N$;
\item $\varphi_{n,1}^{-1}(\tilde M) L_n = L_n \varphi_{n,1}^{-1}(\tilde M)$ for
	all~$n\geq N$.
\end{enumerate}
\end{proposition}

Before proving \cref{prop:Znil_Ypro_explicit}, let us derive \cref{prop:Znil_Ypro}
from it.

\begin{proof} [Proof that \cref{prop:Znil_Ypro_explicit} implies \cref{prop:Znil_Ypro}]
Since every pro-nilsystem is an inverse limit of nilsystem factors,
and since every factor of a pro-nilsystem is a
pro-nilsystem~\cite[Section 3, Chapter 16]{host2018nilpotent},
there is an increasing sequence of nilfactors~$(Y_n)_{n\in \N}$
such that~$Y=\varprojlim_n Y_n$.

Each factor~$Z\vee Y_n$ is an ergodic nilsystem, thus it admits
a reduced representation~$Z\vee Y_n = G_n/\Gamma_n$, and thus up to inductively adjusting
base points
we may assume that the factor maps~$G_n\to G_m$
are given by surjective continuous group homomorphisms~$\varphi_{n,m}\colon G_n\to G_m$
satisfying~$\varphi_{n,m}(\Gamma_n)\subseteq  \Gamma_m$.
The factor~$Z$ can be represented in~$Z\vee Y_1 = G_1/\Gamma_1$ as~$Z=G_1/M$
for some closed subgroup~$M$ with~$\Gamma_1 \subseteq  M \subseteq  G_1$
and in~$Z\vee Y_n = G_n/\Gamma_n$ as~$G_n/\varphi_{n,1}^{-1}(M)$.
Similarly, for every~$n\geq 1$,
the factor~$Y_n$ can be represented in~$Z\vee Y_n$ as~$Y_n=G_n/L_n$,
for some closed subgroup~$L_n$ with~$\Gamma_n \subseteq  L_n \subseteq  G_n$.
Since~$Y_m$ is a factor of $Y_n$ for~$n\geq m$, we have~$L_n\subseteq  \varphi_{n,m}^{-1}(L_m)$
for~$n\geq m$. As each nilfactor is represented as a nilsystem in several ways,
we summarize the current data with a diagram:
for every~$n\geq 1$ we have the following network of systems
\begin{center} \small
\begin{tikzcd} [column sep=11pt, row sep=12pt]
G_n/\Gamma_n \arrow[r] \arrow[dd] \arrow[rrdd]
	& G_n/\varphi_n^{-1}(\Gamma_1) \arrow[rrrr, "\cong"]
	\arrow[ldd,  crossing over] \arrow[rrdd] &                             &                                                           
	&     & G_1/\Gamma_1 \arrow[ldd] \arrow[rdd] &     \\
                                                           &                                                                              &                             &                                                           &     &                                  &     \\
G_n/\varphi_n^{-1}(M) \arrow[rrrr, "\cong", bend right=20] &                                                                              & G_n/L_n \arrow[r] & G_n/\varphi_n^{-1}(L_1)
	\arrow[rrr, "\cong", bend right,  crossing over] & G_1/M &                                  & G_1/L_1
\end{tikzcd}
\end{center}
where $\cong$-labeled arrows are isomorphisms induced by~$\varphi_n:=\varphi_{n,1}$
and unlabeled arrows are canonical projections. The corresponding network of
abstract factors is
\begin{center} \small
\begin{tikzcd} [column sep=36pt, row sep=12pt]
Z\vee Y_n \arrow[r] \arrow[dd] \arrow[rrdd]
	& Z\vee Y_1 \arrow[rrrr, "\cong"]
	\arrow[ldd,  crossing over] \arrow[rrdd,  shorten=5pt,  shift left] &                             &                                                           
	&     & Z\vee Y_1 \arrow[ldd] \arrow[rdd] &     \\
                                                           &                                                                              &                             &                                                           &     &                                  &     \\
Z \arrow[rrrr, "\cong", bend right=25] &                                                                              & Y_n \arrow[r] & Y_1 \arrow[rrr, "\cong", bend right,  crossing over] & Z &                                  & Y_1
\end{tikzcd}.
\end{center}
Note that, while each nilfactor is represented as a nilsystem in several ways,
the ergodic-theoretical properties we are interested in (for example, the good intersection
property) does not depend on the representation.
Let~$N$ and~$\tilde M$ be as given by \cref{prop:Znil_Ypro_explicit}.
Let~$\tilde Z$ be the factor represented by~$G_1/\tilde M$ in~$G_1$,
and by~$G_n/\varphi_{n,1}^{-1}(\tilde M)$ in~$G_n$.
Property~$(i)$ in \cref{prop:Znil_Ypro_explicit} guarantees that~$\tilde Z$
is an extension of~$Z$. Property~$(ii)$ guarantees that~$\tilde Z$ is a factor
of some~$G_n/\Gamma_n$, thus in particular of the ambient system~$X$.
Property~$(iii)$ guarantees that~$\tilde Z:=G_n/\varphi_{n,1}^{-1}(\tilde M)$
and~$Y_n=G_n/L_n$ have the good intersection property for all~$n\geq N$.
We claim that
it follows that~$\tilde Z$ and~$Y$ have the good intersection property.
Let~$f$ be a bounded $\tilde Z$-measurable function.
By the Martingale convergence theorem the function~$\E(f\mid Y_n)$ converges to~$\E(f\mid Y)$.
Since each~$\E(f\mid Y_n)$ is $\tilde Z$-measurable, then~$\E(f\mid Y)$ is $\tilde Z$-measurable. Thus~$\tilde Z$ and~$Y$
have the good intersection property.
\end{proof}

\begin{proof} [Proof of \cref{prop:Znil_Ypro_explicit}]
To shorten notation, we set~$\varphi_n:=\varphi_{n,1}$, $G:=G_1$,
and~$\Gamma:=\Gamma_1$.
Fix~$n\geq 1$.
Since~$\varphi_n(\Gamma_n)\subseteq  \Gamma_1$,
it follows that~$\varphi_n(L_n)$ is a $\Gamma$-rational subgroup. Since~$L_n$ contains~$\Gamma_n$, then~$L_n^0$ is normal in~$G_n$,
and thus~$\varphi_n(L_n^0)$ is normal in~$G$ by surjectivity of~$\varphi_n$.
Moreover, since~$\varphi_n(L_n^0)$ is a connected and clopen subgroup
of~$\varphi_n(L_n)$, we have~$\varphi_n(L_n^0)=\varphi_n(L_n)^0$.\footnote{
This holds in more generality.
Quotient maps between Lie groups are open, and for Lie groups
the connected components are clopen, thus their images are also clopen.}

By applying \cref{fact:Ln_factorization} to the decreasing sequence
of $\Gamma$-rational subgroups~$\varphi_n(L_n)$, we obtain a factorization
\[
	\varphi_n(L_n) = (\varphi_n(L_n))^0 \cdot \Lambda_n
\]
with the following properties: the sequence~$n\mapsto \Lambda_n$ is decreasing,
for every~$n\geq 1$ the group~$\Lambda_n$ is discrete
and contains~$\phi_n(L_n)\cap \Gamma$ with finite index, and there exists
a connected $\Gamma$-rational subgroup~$Q$ such that~$Q=\varphi_n(L_n)^0$
for all~$n\geq N$.
By the previous paragraph, the subgroup~$Q$ is normal in~$G$.

Since~$L_1\cap \Gamma=\Gamma$ has finite index in~$\Lambda_1$, the intersection
\begin{equation}\label{E:tildeGamma}
	\tilde \Gamma := \bigcap_{\lambda\in \Lambda_1} (\lambda \Gamma \lambda^{-1})
\end{equation}
involves only finitely many distinct terms, and thus by~\cref{lemma:commensurable}
the group~$\tilde \Gamma$ is a subgroup of~$\Gamma$ of finite index.
Recall that, since~$M$ contains~$\Gamma$, the identity component~$M^0$
is normal. In particular~$M^0 \cdot \tilde \Gamma$ is a subgroup, and it
has finite index in~$M$.
It follows that
\[
	\tilde M_n:=
	\gen{M^0, \tilde \Gamma, \varphi_n(L_n)\cap M}
\]
is a decreasing sequence of subgroups of~$M$ all containing
the fixed finite index subgroup~$M^0 \cdot \tilde \Gamma$. Thus the
sequence must eventually stabilize: up to increasing the previously defined~$N$,
we may suppose that
\[
	\tilde M :=
	\gen{M^0, \tilde \Gamma, \varphi_N(L_N)\cap M}
\]
satisfies~$\tilde M = \tilde M_n$ for all~$n\geq N$.
Let us verify that~$\tilde M$ satisfies~$(i)$, $(ii)$, and~$(iii)$.

\textit{Proof of $(i)$.}
Recall that~$M^0\cdot \tilde \Gamma \subseteq  \tilde M_n \subseteq  M$
and that~$M^0\cdot \tilde \Gamma$ has finite index in~$M$.
Thus~$\tilde M = \tilde M_N$ is a subgroup of~$M$ of finite index.

\textit{Proof of $(ii)$.}
Fix~$n\geq N$. We
have~$\Gamma_n\subseteq  L_n \subseteq  \varphi_n^{-1}(\varphi_n(L_n))$
and~$\Gamma_n\subseteq  \varphi_n^{-1}(\Gamma) \subseteq  \varphi_n^{-1}(M)$,
from which
\[
    \Gamma_n \subseteq  \varphi_n^{-1}(\varphi_n(L_n)\cap M)
        \subseteq  \varphi_n^{-1}(\varphi_N(L_N)\cap M)
        \subseteq  \varphi_n^{-1}(\tilde M).
\]
This proves~$(ii)$. 

Let us prove~$(iii)$. Let~$n\geq N$. We want to
prove~$\varphi_n^{-1}(\tilde M) L_n = L_n \varphi_n^{-1}(\tilde M)$.
Since~$\varphi_n^{-1}(\varphi_n(L_n)) = L_n \ker(\varphi_n)$
and~$\ker(\varphi_n) \subseteq  \varphi_n^{-1}(\tilde M)$, it is enough to prove
\[
	\tilde M \varphi_n(L_n) = \varphi_n(L_n) \tilde M.
\]
Since~$n\geq N$, we have~$\tilde M = \gen{M^0, \tilde \Gamma, \varphi_n(L_n)\cap M}$.
Since~$M^0$ is normal and~$\varphi_n(L_n)\cap M$ is a subgroup of~$\varphi_n(L_n)$,
they both commute with~$\varphi_n(L_n)$. Thus, by \cref{lem:generated_commute},
it is enough to prove that~$\tilde \Gamma$ and~$\varphi_n(L_n)$ commute.
Since~$\varphi_n(L_n)=Q\cdot \Lambda_n$ and~$Q$ is normal, again
by \cref{lem:generated_commute} it is enough to prove that~$\Lambda_n$ and~$\tilde \Gamma$
commute. By definition, the group~$\tilde \Gamma$ is normalized
by any element of~$\Lambda_1$, thus in particular by any element of~$\Lambda_n$,
concluding the proof.
\end{proof}

\begin{remark}
We use notations as in the proof of \cref{prop:Znil_Ypro_explicit}.
By \cref{lem:generated_commute}, for every~$n$ there exists a unique
largest subgroup~$P_{M}(\varphi_n(L_n))$ of~$M$ that commutes with~$\varphi_n(L_n)$.
We may call this group the ``commutizer'' of~$\varphi_n(L_n)$ in~$M$.
From the proof of \cref{prop:Znil_Ypro_explicit} we see that~$P_{M}(\varphi_n(L_n))$
has uniformly bounded index in~$M$. However,
although~$\varphi_n(L_n)$ is monotone decreasing, we see no reason why
the sequence~$P_{M}(\varphi_n(L_n))$
should stabilize (as opposed to commutizer and normalizer,
monotonicity is not a property of the commutizer).
\end{remark}

As a corollary, we obtain the following.
\begin{proposition} \label{prop:pronil_rectified}
Let~$X$ be an ergodic pronilsystem. Let~$Y$ be a factor.
Then there exists an increasing sequence of nilfactors~$(X_n)_{n\in \N}$ such that:
\begin{enumerate} [(i)]
\item $X=\varprojlim_n X_n$;
\item $X_n$ and~$Y$ have the good intersection property for all~$n\in \N$.
\end{enumerate}
\end{proposition}
\begin{proof}
Let~$(X_n)_{n\in \N}$ be any increasing sequence of nilfactors
with~$X=\varprojlim_n X_n$.
Suppose that~$(ii)$ holds for all~$1\leq n\leq m-1$ for some~$m\geq 1$
(if~$m=1$, this set is empty).
We will show how to modify the tail~$(X_n)_{n\geq m}$ so
that~$(ii)$ holds for all~$1\leq n\leq m$. In countably many steps,
one obtains the desired result.

By \cref{prop:Znil_Ypro}, there is a nilfactor~$\tilde X_m$
with~$X_m\subseteq  \tilde X_m \subseteq  X$
such that~$\tilde X_m$ and~$Y$ have the good intersection property. Define
\[
    X'_n :=
    \begin{dcases*}
        X_n, & if $n\leq m-1$, \\
        \tilde X_m, & if $n=m$, \\
        X_n \vee \tilde X_m, & if $n>m$.
    \end{dcases*}
\]
Then~$(X'_n)_{n\in \N}$ is an increasing sequence of factors converging to~$X$
such that~$(ii)$ holds for all~$1\leq n \leq m$.
\end{proof}

We are now ready to generalize
\cref{prop:key_lem_nilsystem} to pro-nilsystems.

\begin{proposition} \label{prop:key_lem_pro}
Let~$X$ be an ergodic pro-nilsystem and~$Y$
a factor such that the extension~$X\to Y$
is compact. Then the conclusion of \cref{prop:key_lem_nilsystem} holds.
\end{proposition}
\begin{proof}
Let~$(X_n)_{n\in \N}$ be as in \cref{prop:pronil_rectified}.
By \cref{L:compactfactor}, for every~$m\in \N$ the
extension~$X_m\to Y\wedge X_m$ is compact.
Let~$\CF_m \subseteq  L^\infty(X_m)$ be the countable set obtained
by applying \cref{prop:key_lem_nilsystem} to the compact
extension~$X_m\to Y\wedge X_m$. Let~$f\in \CF_m$.
Since $T^n f - f$ is $X_m$-measurable  and the factors $X_m,Y$ have the good intersection property, we have
\[
	\E(\abs{T^n f - f}^2\mid Y\wedge X_m) = \E(\abs{T^n f - f}^2\mid Y).
\]
Therefore, the function~$f$ (seen as an $X_m$-measurable function on~$X$)
has the desired $\nilbor$ property with
respect to the extension~$X\to Y$.
Now let~$\CF = \bigcup_{m} \CF_m$. Since~$\CF_m$ is dense in~$L^2(X_m)$
and~$X=\varprojlim_m X_m$, it follows that~$\CF$ is dense in~$L^2(X)$.
\end{proof}


\section{Compact extensions in non-ergodic $\CZ_s$-systems}
\label{sec:iso_Zs}

In the previous section, we proved \cref{thm:beautiful_main}
under the assumption of ergodicity. In this section, we prove
the theorem in its full generality.

\subsection{The idea of countable control}

As we will see, a compact extension disintegrates,
through the ergodic decomposition, into a family of compact extensions.
In contrast, an extension disintegrating into compact extensions
may not be compact itself. Roughly speaking, we have
\[
\begin{tikzcd}[column sep=tiny, row sep=0pt]
\text{compact extension} & \iff &
{\begin{array}{c}
   \text{compact extension} \\ \text{in each component}
 \end{array}} & + & \text{countable control}.
\end{tikzcd}
\]
We will make this precise in the context of $\CZ_s$-systems.
Before that, we collect some preliminary results and examples.
We begin by proving the disintegration claim.

\begin{lemma} \label{lem:compact_disintegration}
Let~$(X,\CX,T,\mu)$ be a system with ergodic
decomposition $(\mu_{\omega})_\omega$. Let~$\CY\subseteq  \CX$ be a $T$-invariant
sub-sigma-algebra such that the extension
\[
	(X,\CX,T,\mu) \mapsto (X,\CY,T,\mu)
\]
is compact. Then for almost every ergodic component~$\mu_\omega$ the extension
\[
	(X,\CX,T,\mu_\omega) \mapsto (X,\CY,T,\mu_\omega)
\]
is also compact.
\end{lemma}
\begin{proof}
For each ergodic component,
let~$L_\omega \subseteq  L^2(X,\mu_\omega)$ be the closure of the linear span
of all $T$-invariant finite-rank $L^\infty(Y,\mu_\omega)$-modules.
We want to prove that~$L_\omega = L^2(X,\mu_\omega)$ for almost every~$\mu_\omega$.
Since~$X$ is standard Borel, we can fix a countable family of Borel
functions~$\CF$ that is dense in~$L^2(X,\rho)$ for every Borel
probability measure~$\rho$.
By countability,
it is enough to prove that for every~$f\in \CF$
we have~$f \in_{\mu_\omega} L_\omega$ for almost every~$\omega$.
Fix~$f\in \CF$. By compactness of the extension, we can approximate~$f$ by a sequence
of functions~$(f_k)_{k\in \N}$ with~$f_k$ belonging to a $T$-invariant
finite-rank $L^\infty(Y,\mu)$-module
and~$\norm{f-f_k}_{L^2(X,\mu)}\leq 2^{-k}$ for all~$k\in \N$.
By Borel-Cantelli, we have~$\norm{f-f_k}_{L^2(X,\mu_\omega)}\to 0$ for almost every~$\omega$.
It remains to notice that, for almost every~$\omega$, the function~$f_k$
belongs to a $T$-invariant
finite-rank $L^\infty(Y,\mu_\omega)$-module:
let~$f_k =_\mu \sum_{i=1}^r h_i g_i$ with~$h_i\in L^\infty(Y,\mu)$,
then~$f_k =_{\mu_\omega} \sum_{i=1}^r h_i g_i$ for almost every~$\omega$,
and moreover~$h_i \in L^\infty(Y,\mu_\omega)$ for almost every~$\omega$.
\end{proof}

Let~$\pi\colon X\to Y$ be a realization of the factor~$(X,\CY,T,\mu)$,
and let
\[
    \kappa\colon I(X)\to I(Y)
\]
be the induced map between invariant factors.
Throughout this section we use the shorthand
\[
    \Omega:=I(X), \qquad \Lambda:=I(Y).
\]
For almost every~$\omega\in \Omega$, the measure~$\mu_\omega$
is supported on~$X_\omega$, the fiber over~$\omega$ of the map~$X\to \Omega$,
and the factor map~$\pi$ induces
a factor map~$\pi_\omega\colon X_\omega\to Y_{\kappa(\omega)}$
between ergodic systems.
In this language, \cref{lem:compact_disintegration} states that
if the extension~$X\to Y$ is compact, then for almost every~$\omega$
the extension~$X_\omega\to Y_{\kappa(\omega)}$ is compact.
Note that several ergodic components of~$X$ may be mapped to the same
ergodic component of~$Y$; this happens exactly when~$I(Y)\neq I(X)$.
For a factor~$W$ of~$X$ containing~$I(X)$,
we write~$W_\omega$ for the factor of~$X_\omega$
given by the sub-$\sigma$-algebra of~$W$, taken modulo $\mu_\omega$-nullsets.
For this to be meaningful, we require factor maps
to be Borel,
so that the corresponding sub-$\sigma$-algebras are
countably generated;
see \cref{app:measure} for details on this technical point.


The converse of \cref{lem:compact_disintegration} is false.
Simply take~$X=\T^2$
with~$T(x,y)=(x,x+y)$ and let~$Y$ be the trivial factor.
The following example (in which~$Y$ is trivial)
introduces the idea of  countable control and how it implies
our desired $\nilbor$ property.

\begin{example} [Countable Control I] \label{ex:countable_control_I}
Consider the transformation~$T\acts X:=\T^2$ given by~$T(x,y):=(x,y+\sigma(x))$,
for some measurable function~$\sigma\colon\T\to \T$ with irrational values.
Restricted to every component (i.e. fixing~$x$), the transformation is compact.
However, the transformation~$T$ itself is compact if and only if~$\sigma$
has essentially countable range. The ``if'' direction is straightforward.
Let us prove the ``only if'' direction. Assume~$T$ is compact. Then for every non-zero
bounded measurable function~$f$ we have
\[
	0 < \lim_{N\to\infty} \frac{1}{N}\sum_{n=1}^N
			\abs{\int T^n f \cdot \con{f} \d \m_{\T^2}}^2.
\]
Let~$f(x,y):=\1_A(x)e(y)$ with~$A$ a set of positive measure
and~$e(y):=\exp(2\pi i y)$. Substituting, we obtain
\[
	0 < \lim_{N\to\infty} \frac{1}{N}\sum_{n=1}^N
		\abs{\int_A e(\sigma(x))^n \d x}^2
		= \int_{A\times A} \1_{[\sigma(x)=\sigma(x')]} \d x \d x'.
\]
Therefore, every positive measure set~$A$
(which we think of as a set of ergodic components) contains a subset of positive
measure where the restriction of~$\sigma$ is constant. It follows that~$\sigma$
has countable essential range, i.e., there is a countable set~$\{\sigma_m\}_{m\in \N}$
containing~$\sigma(x)$ for almost every~$x$. Let~$E_m := \sigma^{-1}(\sigma_m)$.
The system is spanned by functions~$f(x,y):=\1_{E_m}(x) e(bx+ay)$,
with~$m\in \N$ and~$a,b\in \Z$. For such functions the quantity
\begin{align*}
	\int \abs{T^n f - f}^2 \d \m_{\T^2}
    & = \int \abs{\1_{E_m}(x) e(a n \sigma_m + ay)
        - \1_{E_m}(x) e(ay)}^2 \d \m_{\T^2}(x,y) \\
	& \leq  \abs{e(a\sigma_m)^{n} - 1}^2
\end{align*}
is smaller than a given~$\ve>0$ on a $\nilbor$ set of integers~$n$.
\end{example}

To check countability along fibers when the base~$Y$ is non-trivial,
we will use the following elementary lemma.

\begin{lemma} \label{lem:countable_range}
Let~$\pi\colon (X,\mu)\to (Y,\nu)$ be a factor map between probability spaces
with disintegration~$\mu=\int_Y \mu_y \d\nu(y)$.
Let~$f\colon X\to \C$ be measurable. Suppose
that for every measurable~$B\subseteq  X$ with~$\mu(B)>0$ we have
\[
	\int_Y \int_{B\times B} \1_{\{f(x)=f(x')\}} \d \mu_y(x) \d\mu_y(x') \d\nu(y) > 0.
\]
Then for almost every~$y\in Y$ the restriction of~$f$ to the fiber over~$y$
has essentially countable range (that is, the map~$f$ has $\mu_y$-essentially countable
range).
\end{lemma}
\begin{proof}
For~$x\in X$, consider the level set~$[x]:=\{x'\in X\colon f(x')=f(x)\}$.
Define~$B:=\{x\in X\colon \mu_{\pi(x)}([x])=0\}$. The set~$B$ is measurable:
it is the zero set of the
map~$x\mapsto \int_X \1_{\{(x,x')\colon f(x)=f(x')\}} \d\mu_{\pi(x)}(x')$, which can
be seen to be measurable by approximating~$\1_{\{(x,x')\colon f(x)=f(x')\}}$ by rectangles
and using the monotone class theorem. By definition of~$B$ we have
\[
	\int_Y \int_{B\times B} \1_{\{f(x)=f(x')\}} \d \mu_y(x) \d\mu_y(x') \d\nu(y) = 0.
\]
Therefore, by hypothesis, we must have~$\mu(B)=0$,
in particular~$\mu_y(B)=0$ for almost every~$y$.
In turn~$\mu_y(B)=0$ implies that for $\mu_y$-almost every~$x$ the level set~$[x]$
has positive $\mu_y$-measure. In particular, there are only countably many values.
\end{proof}

\begin{remark}
The converse of \cref{lem:countable_range} seems also true, but we will not need
it.
\end{remark}

We now give an example in which $Y$ is nontrivial: the system
is a ``nilbundle'' (definition given below), that is,
a measurable varying collection
of nilsystems, parametrized by the invariant factor.
In the following example, the system is
a measurable varying collection of affine unipotent transformations
on the torus.

\begin{example} [Countable Control II] \label{ex:countable_control_II}
Let~$T\acts X:=\T^4$ be given by
\[
	(t,u,v,w)\mapsto (t,u+t, v,w+u+\sigma(t,v)),
\]
for some measurable function~$\sigma$.
Note that~$(X,m_{\T^4},T)$ is a non-ergodic $\CZ_2$-system,
since for fixed~$t,v$ the restriction is a $\CZ_2$-system.
Let~$Y$ be the factor generated by $t,u$; the induced map is $(t,u)\mapsto (t,u+t)$.
Suppose that the extension~$X\to Y$ is compact.
We will verify
that the extension~$X\to Y$ has the desired $\nilbor$ property.
To this end, we first prove that for almost every~$t$ the map~$v\mapsto \sigma(t,v)$
has countable essential range.
Let~$A\subseteq \T$ be a subset of positive measure and~$(B_t)_{t\in \T}$
a measurable varying family such that~$B_t\subseteq  \T$
has positive measure if~$t\in A$.
Define~$f(t,u,v,w):=\1_A(t)\1_{B_t}(v)e(w)$.
Since~$X\to Y$ is compact and~$f\neq 0$, we have
\[
	\lim_{N\to\infty} \frac{1}{N}\sum_{n=1}^N
	\int \abs{\int (T^n f\cdot \con{f})(t,u,v,w) \d v \d w}^2 \d t \d u > 0,
\]
which gives~$\lim_{N\to\infty} \frac{1}{N}\sum_{n=1}^N
		\int_A
			\abs{\int_{B_t} e(n\sigma(t,v)) \d v}^2
		\d t >0$,
which in turn gives
\[
		\int_A
			\p{\int_{B_t \times B_t} \1_{[\sigma(t,v)=\sigma(t,v')]}\d v \d v'}
		\d t >0.
\]
By \cref{lem:countable_range},
for almost every~$t$ the map~$v\mapsto \sigma(t,v)$
has countable essential range.
By Lusin-Novikov's Theorem~\cite[Theorem 18.10]{kechris2012classical},
we can track this range measurably,
i.e., there are measurable maps~$(t\mapsto \sigma_m(t))_{m\in \N}$
with the following property:
let~$E_m(t) := \{v\colon \sigma(t,v)= \sigma_m(t)\}$,
then for almost every~$t$ the union~$\bigcup_{m\in \N} E_m(t)$ has full measure.
We now prove the desired $\nilbor$~property.
The system is spanned by countably many functions of the
form
\[
	f(t,u,v,w)= \1_{E_m(t)}(v) e(A t)\, e(B u)e(C v)e(D w)
\]
with 
$m\in \N$ and integers~$A,B,C,D$. Any such
function satisfies
\[
	\abs{T^n f(t,u,v,w) - f(t,u,v,w)}
	\leq \abs{e(Bnt) \, e\Big(D(nu+{n\choose 2}t + n\sigma_m(t))\Big) - 1},
\]
giving the desired $\nilbor$ property.
\end{example}


\subsection{Nilbundles}
We borrow from~\cite{jamneshan2026non} the notion of bundle of nilsystems,
which we call nilbundle for short.

\begin{definition} \label{def:nilbundle}
A \emph{simple nilbundle} is a system of the form
\[
	\Omega\ltimes_\tau G/\Gamma:=(\Omega\times G/\Gamma, \mu\times \m_{G/\Gamma}, T),
	\quad T(\omega,x):=(\omega,\tau(\omega)x),
\]
where~$(\Omega,\mu)$, called \emph{base},
is a standard probability space, $(G/\Gamma,\m_{G/\Gamma})$
is a nilmanifold, and~$\tau\colon \Omega\to G$ is a measurable map
such that the action~$\tau(\omega) \acts G/\Gamma$ is ergodic
for $\mu$-almost every~$\omega\in \Omega$.
A \emph{nilbundle} is a disjoint measurable union of simple nilbundles, that is,
\[
	\bigsqcup_{i\in \N} (\Omega_i \ltimes_{\tau_i} G_i/\Gamma_i)
\]
where~$(\Omega_i)_{i\in \N}$ is a measurable partition of~$\Omega$
and each~$\Omega_i \ltimes_{\tau_i} G_i/\Gamma_i$ is a simple nilbundle.
We call~$\Omega$ the \emph{base} of the nilbundle.
If all~$G_i/\Gamma_i$ are $s$-step, we call the nilbundle \emph{$s$-step}.
\end{definition}

It follows from the definition
that the base of a nilbundle coincides with the invariant factor
of the nilbundle.
An example of  a nilbundle is
the system in \cref{ex:countable_control_II}.

A.~Jamneshan and S.~Machado generalized
the Structure Theorem of B.~Host and B.~Kra as follows.

\begin{theorem} [Non-Ergodic Structure Theorem \cite{jamneshan2026non}] \label{T:NonErgodicStructure}
Every~$\CZ_s$-system is an inverse limit of $s$-step nilbundle factors
with base the invariant factor.
\end{theorem}

In fact, in~\cite{jamneshan2026non} they prove the analogous
result for~$\Z^d$-actions, but we will not need such generality here.
Let us unpack the notion of nilbundle factor.
Let~$\Omega:=I(X)$ be the invariant factor
of a system~$X$. A {\em nilbundle factor with base~$\Omega$}
is given by a factor map
\[
    \xi\colon X\to
        \bigsqcup_{i\in \N} (\Omega_i \ltimes_{\tau_i} G_i/\Gamma_i),
\]
defined as follows. The right hand side is a nilbundle
as in \cref{def:nilbundle}. Let~$\omega(x)\in \Omega$ denote the ergodic
component containing~$x$. The map~$\xi$ is defined by cases, one case for
each partition element~$\Omega_i$:
for~$\omega(x)\in \Omega_i$ we have
\[
    \xi(x) = (\omega(x),g(x)\Gamma_i)
\]
where~$g(x)\in G_i$ is a measurable varying group element, and such that
\[
    \xi(Tx) = (\omega(x), \tau_i(\omega(x)) g(x)\Gamma_i).
\]
All previous identities are understood to hold almost everywhere.
We remark that, in general, a nilbundle factor~$Y$ of~$X$ might have base~$I(Y)$
smaller than~$I(X)$. We will return to this in detail below.

Before discussing countable control and $\nilbor$ properties
for nilbundles, we need to extend to nilbundles some structural results
known to hold for nilsystems. 

\begin{lemma} \label{joint_nilbundle_nilbundle}
If~$X$ a system and~$X_1,\ldots,X_k$ are nilbundle factors,
then~$X_1\vee \cdots \vee X_k$ is a nilbundle factor.
\end{lemma}
\begin{proof}
This is \cite[Lemma 3.5]{jamneshan2026non}.
\end{proof}

\begin{lemma} \label{factor_nilbundle_nilbundle}
Any factor of an $s$-step  nilbundle is an $s$-step nilbundle.
\end{lemma}

\begin{proof}
Let~$(X_\omega)_{\omega\in \Omega}$ and~$(Y_{\lambda})_{\lambda\in \Lambda}$
be ergodic decompositions of~$X$ and~$Y$, where~$\Omega=I(X)$ and~$\Lambda=I(Y)=I(X)\wedge Y$.
The factor map~$X\to Y$ induces factor maps between ergodic components.
Almost every ergodic component~$Y_\lambda$ of~$Y$
is the image of an ergodic component~$X_\omega$ of~$X$.
In turn, almost every ergodic component of~$X$ is isomorphic to an $s$-step nilsystem.
Since factors of $s$-step nilsystems are $s$-step nilsystems,
(see preliminaries, \cref{SS:nilfactors})
it follows that almost every ergodic
component of~$Y$ is isomorphic to an $s$-step nilsystem.
The nontrivial claim is that this isomorphism can be taken measurably across components.

Since~$X$ is a $\CZ_s$-system, so is~$Y$. In particular
\cref{T:NonErgodicStructure} applies: the system~$Y$ is spanned
by an increasing sequence
of nilbundle factors~$\xi_n\colon Y\to Y_n$, with base~$\Lambda$.
For almost every~$\lambda$ we have an increasing sequence~$(Y_n)_\lambda$
of factors of the nilsystem~$Y_\lambda$,
and by \cref{cor:nilsystems_stabilize}
this sequence must stabilize.
In particular, up to removing from~$\Lambda$ a nullset, we have a countable measurable
partition~$\Lambda=\bigsqcup_{i\in \N} \Lambda_i$ such that~$(Y_n)_\lambda=Y_\lambda$
for all~$n\geq i$ if~$\lambda\in \Lambda_i$.
The map~$\xi = \sum_{i\in \N} \1_{\Lambda_i} \xi_i$
is a measurable isomorphism between~$Y$ and a nilbundle.
\end{proof}

For technical reasons, we introduce the notion of reduced nilbundle, analogous
to that of reduced nilsystems (see preliminaries \cref{sec:preliminaries_reduced}).

\begin{definition}
We say that a nilbundle~$X = \bigsqcup_i (\Omega_i \ltimes_{\tau_i} G_i/\Gamma_i)$
is \emph{reduced} if almost every component is a reduced nilsystem,
that is,~$\Gamma_i$ contains no non-trivial normal subgroups of~$G_i$ for all~$i\in \N$
and for almost every~$\omega\in \Omega_i$ the group~$G_i$ is spanned by~$G_i^0$
and~$\tau_i(\omega)$.
\end{definition}

\begin{lemma} \label{lem:reduced_nilbundle}
Every nilbundle admits a reduced representation.
\end{lemma}
\begin{proof}
Since a nilbundle is a disjoint union of simple nilbundles,
it is enough to prove the claim for a simple nilbundle
\[
	X:=\Omega\ltimes_\tau G/\Gamma, \quad \tau\colon \Omega\to G, \quad \Omega:=I(X).
\]
We will show that the argument in \cite[Chapter 11, Section 1.2]{host2018nilpotent},
used to prove that every ergodic nilsystem admits a reduced representation,
extends to nilbundles.

\textit{Generation.}
Fix a generic~$\omega$. Let~$G'(\omega)$ be the group generated by~$G^0$ and~$\tau(\omega)$.
Let us show that~$\omega\mapsto G'(\omega)$ can be trivialized over a countable measurable
partition of~$\Omega$. Recall that~$G/G^0$ is countable and~$\omega\mapsto \tau(\omega)$
is measurable. In particular, there exists a countable measurable partition
trivializing~$\omega\mapsto \tau(\omega) G^0$. Since~$G'(\omega)=\gen{G^0,\tau(\omega)}$
is uniquely determined by~$\tau(\omega) G^0$, such partition trivializes~$G'(\omega)$.
Up to restricting to a single part, we may from now on suppose~$G'(\omega)=G'$
for almost every~$\omega$.
Since~$G^0$ is open, then~$G'$ is open, and thus
the orbit~$G' 1_{G/\Gamma}$ is open. It is also closed, since its complement is also a union of $G'$ orbits, each of which is open.
Since~$\tau(\omega)\acts G/\Gamma$ is ergodic, such orbit is also dense,
that is~$G'\Gamma = G$. In particular, there exists a measurable
isomorphism~$G/\Gamma \cong G'/(G'\cap \Gamma)$ that preserves the action of~$\tau(\omega)$.
Up to replacing~$(G/\Gamma,\tau(\omega))$ by~$(G'/(G'\cap \Gamma),\tau(\omega))$
we obtain~$G=\gen{G^0,\tau(\omega)}$ for almost every~$\omega$.

\textit{Faithfulness.}
Finally, let~$\Gamma'$ be the largest normal subgroup of~$G$ contained in~$\Gamma$.
Up to replacing~$G$, $\Gamma$, $\tau$ by~$G/\Gamma'$, $\Gamma/\Gamma'$,
$\tau\Gamma'$ respectively, we obtain that~$\Gamma$ contains no non-trivial
normal subgroups of~$G$.
\end{proof}

We will also need the following technical lemma about nilmanifold morphisms.

\begin{lemma} \label{hom_nilmanifolds}
Let~$G/\Gamma$ and~$G'/\Gamma'$ be nilmanifolds
with~$\bigcap_{g\in G'} (g^{-1} \Gamma' g) = \{1_{G'}\}$. Then
\[
	\Phi(G/\Gamma,G'/\Gamma'),
\]
the set of all continuous surjective group homomorphisms~$\varphi\colon G\to G'$
with~$\varphi(\Gamma)\subseteq  \Gamma'$, is countable.
Moreover, two such homomorphisms~$\varphi_1,\varphi_2 \in \Phi(G/\Gamma,G'/\Gamma')$
are equal as functions on~$G/\Gamma$ if and only if they are equal as functions on~$G$.
\end{lemma}
\begin{proof}
We start by proving the second claim. Let~$\varphi_1,\varphi_2$ be as in the statement
and suppose that~$\varphi_1(g\Gamma)=\varphi_2(g\Gamma)$ for all~$g\in G$.
Then~$\varphi_1(g)^{-1}\varphi_2(g)\in \Gamma'$ for all~$g\in G$. In particular,
the map~$\delta(g):=\varphi_1(g)^{-1}\varphi_2(g)$ takes values in~$\Gamma'$.
In particular, for any~$g,h\in G$, a short computation shows that
\[
	\varphi_2(g)^{-1} \delta(h) \varphi_2(g)
	= \delta(g)^{-1} \delta(hg) \in \Gamma'.
\]
Thus, by
surjectivity of~$\varphi_2$,
we have~$\delta(h) \in \bigcap_{g\in G'} (g^{-1} \Gamma' g) = \{1_{G'}\}$ for any~$h\in G$,
proving that~$\varphi_1=\varphi_2$ as maps on~$G$.

Let us now prove countability. Let~$\varphi\in \Phi(G/\Gamma,G'/\Gamma')$.
Since~$\varphi(\Gamma)\subseteq  \Gamma'$, rational elements are mapped to rational elements.
In particular, for any rational element~$g\in G$, there
are only countably many possibilities for~$\varphi(g)$. Since~$G$ is generated by~$G^0$
and finitely many rational elements, we can suppose~$G=G^0$, that is~$G$ connected.
Up to moving to
a universal cover (recall that a continuous surjective homomorphism lifts uniquely
to the universal cover), we may also suppose~$G$ simply connected. Then~$\varphi$ is
uniquely determined by its differential~$d\varphi$. In turn~$d\varphi$ is uniquely
determined by its values on~$\log \Gamma$ (which spans the Lie algebra).
Since~$\Gamma$ is finitely generated, there are only countably many such maps.
\end{proof}

A factor map between ergodic systems in reduced form is algebraic:
up to a change of base point, it is induced by a continuous
surjective group homomorphism.
The following result extends this fact to nilbundles.

\begin{proposition} \label{prop:nilbundle_factor_representation}
Let~$X$ be a nilbundle, $\pi\colon X\to Y$ a factor map,
$\kappa\colon I(X)\to I(Y)$ the induced factor map between invariant factors,
and denote~$\Omega:=I(X)$ and~$\Lambda:=I(Y)$. Then there exist nilbundle
representations of~$X$ and~$Y$ such that~$\pi\colon X\to Y$ takes the form
\[
	\pi\colon \bigsqcup_{i\in \N}( \Omega_i\ltimes_{\tau_i} G_i/\Gamma_i)
		\to \bigsqcup_{j\in \N} (\Lambda_j\ltimes_{\tau'_j} G'_j/\Gamma'_j)	
\]
where:
$\Omega=\bigsqcup_{i\in \N} \Omega_i$ and
$\Lambda=\bigsqcup_{j\in \N} \Lambda_j$ are measurable partitions
satisfying~$\Omega_i\subseteq  \kappa^{-1}(\Lambda_{\iota(i)})$
for some map~$\iota\colon \N\to \N$;
for every~$i\in \N$ there is a surjective continuous
group homomorphism~$\varphi_i\colon G_i\to G'_{\iota(i)}$
satisfying~$\varphi_i(\Gamma_i)\subseteq  \Gamma'_{\iota(i)}$,
in particular~$L_i:=\varphi_i^{-1}(\Gamma'_{\iota(i)})$ is a closed subgroup
containing~$\Gamma_i$;
the maps~$\tau_i\colon \Omega_i\to G_i$ and~$\tau'_j\colon \Lambda_j\to G'_j$
are measurable; restricted to any~$\Omega_i$, the map~$\pi$ induces a factor map between
simple nilbundles
\[
	\Omega_i\ltimes_{\tau_i} G_i/\Gamma_i
	\to \Lambda_{\iota(i)} \ltimes_{\tau'_{\iota(i)}} G'_{\iota(i)}/\Gamma'_{\iota(i)},
	\qquad (\omega,g\Gamma_i) \mapsto (\kappa(\omega),\varphi_i(g)\Gamma'_{\iota(i)}),
\]
and for almost every ergodic component~$\omega\in \Omega_i$
it induces a factor map between ergodic nilsystems in reduced form
\[
	(G_i/\Gamma_i,\tau_i(\omega))
		\to (G'_{\iota(i)}/\Gamma'_{\iota(i)},\tau'_{\iota(i)}(\kappa(\omega))), \qquad
	g\Gamma_i \mapsto \varphi_i(g)\Gamma'_{\iota(i)}.
\]
\end{proposition}

\begin{proof} We organize the proof in several steps.

\textit{Partition refinement.}
By \cref{factor_nilbundle_nilbundle}, the factor~$Y$ is a nilbundle.
By \cref{lem:reduced_nilbundle}, there exist
reduced representations~$X=\bigsqcup_{i\in \N} (\Omega_i\ltimes_{\tau_i} G_i/\Gamma_i)$
and~$Y=\bigsqcup_{j\in \N} (\Lambda_j\ltimes_{\tau'_j} G'_j/\Gamma'_j)$.
Since~$\Lambda$ is a factor of~$\Omega$, the partition~$(\kappa^{-1}(\Lambda_j))_{j\in \N}$
is $\Omega$-measurable. Since the joining of two countable partitions
is countable, up to refining~$(\Omega_i)_{i\in \N}$, we obtain the
desired selector~$\iota\colon \N\to \N$
with~$\Omega_i\subseteq  \kappa^{-1}(\Lambda_{\iota(i)})$.

\textit{Change of base points.} We will perform a measurable varying change of base point in
the $X$-representation so that the factor map preserves base points over each ergodic
component. Fix a generic~$\omega\in \Omega$, let~$\lambda:=\kappa(\omega)$
and~$j=\iota(i)$. The map~$\pi$ induces a factor map between ergodic components
\begin{equation} \label{eq:nilbundle_repr_factor}
	(G_i/\Gamma_i,\tau_i(\omega))
		\to (G'_{\iota(i)}/\Gamma'_{\iota(i)},\tau'_{\iota(i)}(\kappa(\omega))), \qquad
	g\Gamma_i \mapsto \pi_\omega(g\Gamma_i).
\end{equation}
Since these are reduced ergodic nilsystems, the map~$\pi_\omega$
is algebraic (see \cref{sec:prel_nilsystem} for details):
it is the composition of a continuous surjective group homomorphism
(which preserves the base point)
with a nilrotation (which may not preserve the base point).
Since~$\pi\vert_{\Omega_i\times G_i/\Gamma_i}$ is measurable,
in particular~$\omega\mapsto \pi_\omega(1_{G_i/\Gamma_i})$, the image of the base
point, is measurable. By fixing a section, we find  a
measurable varying element~$\beta_\omega \in G_i$
such that~$\pi_\omega(\beta_\omega 1_{G_i/\Gamma_i})=1_{G'_j/\Gamma'_j}$.
The transformation
\[
	\Omega_i\ltimes_{\beta^{-1} \tau_i \beta} G_i/\Gamma_i \to
	\Omega_i\ltimes_{\tau_i} G_i/\Gamma_i,
	\quad (\omega,g\Gamma) \mapsto (\omega,\beta_\omega g \Gamma)
\]
is a measurable isomorphism, and it intertwines dynamics since
\[
	\tau_i(\omega)\cdot \beta_\omega g\Gamma
	= \beta_\omega \cdot \beta_\omega^{-1} \tau_i(\omega) \beta_\omega g \Gamma.
\]
Therefore, up to precomposing the factor map~\eqref{eq:nilbundle_repr_factor}
with such an isomorphism, and up to relabeling~$\beta^{-1}_\omega \tau_i(\omega) \beta_\omega$
as~$\tau_i(\omega)$, we can suppose that~$\pi_\omega$ preserves base points.

\textit{Homomorphism selection.}
The map~$\pi_\omega$ is now an algebraic factor map that preserves
base points, thus it is induced by a surjective continuous
group homomorphism~$\varphi_\omega\in \Phi(G_i/\Gamma_i,G'_j/\Gamma'_j)$.
By \cref{hom_nilmanifolds}, the set~$\Phi(G_i/\Gamma_i,G'_j/\Gamma'_j)$
of such homomorphisms is countable. We now show that~$\varphi_\omega$ is
measurable in~$\omega$.
Let~$D\subseteq  G_i$ be a countable set such
that~$\{g1_{G_i/\Gamma_i}\colon g\in D\}$ is dense in~$G_i/\Gamma_i$.
Let~$\varphi\in \Phi(G_i/\Gamma_i,G'_j/\Gamma'_j)$.
By \cref{hom_nilmanifolds}, we have~$\varphi_\omega=\varphi$ (as maps on~$G_i$)
if and only if
\[
	\pi_\omega(g1_{G_i/\Gamma_i})
		= \varphi_\omega(g1_{G_i/\Gamma_i})=\varphi(g1_{G_i/\Gamma_i})
\]
for every~$g\in D$. Since the left hand side is measurable in~$\omega$ and~$D$ is countable,
the set~$\{\omega\colon \varphi_\omega=\varphi\}$ is measurable.
Since~$\Phi(G_i/\Gamma_i,G'_j/\Gamma'_j)$ is countable, we conclude
that, up to refining the $\Omega$-partition, we may suppose~$\varphi_\omega=\varphi_i$
for~$\omega\in\Omega_i$.
\end{proof}

\subsection{Countable control in nilbundles}

Our next goal is to extend to general nilbundles
the countable control of \cref{ex:countable_control_I}
and \cref{ex:countable_control_II}.
Recall that a nilbundle is countable union of simple nilbundles.
We first present countable control in the special case of simple nilbundles. We will use the following lemma. 
\begin{lemma}\label{L:CondExp}
Let $(\Omega, \mu)$ and $(M,m)$ be probability spaces, let $\Lambda$ be a factor of $\Omega$, and let $N$ be a factor of $M$. Then for every $w\in L^\infty(\Omega)$ and $F\in L^\infty(\Lambda\times M)$, we have 
$$
\E_{\mu\otimes m}(w\cdot  F \mid \Lambda\times N)=\E_\mu(w\mid \Lambda)\cdot \E_{\mu\otimes m}(F\mid \Lambda\times N)
$$
\end{lemma}
\begin{proof}
Let $h\in L^\infty(\Lambda \times N)$ and define
$$
\Phi(\omega):=\int_M F(\omega,z)\,  h(\omega, z) \d\m(z). 
$$
Since $F$ and $h$ depend on the $\omega$-coordinate only through the factor $\Lambda$, we have $\Phi\in L^\infty(\Lambda)$.
Moreover we have
\begin{align*}
\int_{\Omega\times M} w(\omega) \, F(\omega,z) \, h(\omega, z) \d(\mu\otimes m)=& \int _\Omega \, w(\omega) \, \Phi(\omega) \d\mu 
= \int_\Omega \E_\mu(w\mid \Lambda) (\omega) \, \Phi(\omega) \d\mu \\
=& \int_{\Omega \times M} \E_\mu(w\mid \Lambda)(\omega) \, F(\omega,z) \, h(\omega,z) \d(\mu\otimes m), 
\end{align*}
where we used Fubini's theorem for the first equality, the $\Lambda$-measurability of $\Phi$ for the second, and Fubini's theorem again  together with the fact that 
$\E(w\mid \Lambda)$ is $\Omega$-measurable for the third.  Since $h\in L^\infty(\Lambda \times N)$ is arbitrary, we deduce that 
$$
\E_{\mu\otimes m}(w\cdot F \mid \Lambda\times N)=\E_{\mu\otimes m}\big( \E_\mu(w\mid \Lambda)\cdot F\mid \Lambda\times N\big).
 $$  
 The asserted identity follows by pulling the $(\Lambda\times N)$-measurable function 
 $\E_\mu(w\mid \Lambda)$ outside the conditional expectation. 
 \end{proof}

\begin{proposition} \label{countable_control_simple_nilbundle}
Consider a simple nilbundle in reduced form
\[
	X:=\Omega\ltimes_\tau G/\Gamma, \quad \tau\colon \Omega\to G, \quad \Omega:=I(X),
\]
and a factor~$Y$ of the form
\[
	Y := \Lambda \ltimes_{\tau'} G/L, \quad \tau'\colon \Lambda\to G,
	\quad \Lambda:=I(Y)=I(X)\wedge Y
\]
where~$L$ is a closed subgroup of~$G$ containing~$\Gamma$.
Suppose that the extension~$\pi\colon X\to Y$ is compact and the factor map has the form $\pi(\omega,g\Gamma)=(\kappa(\omega),gL)$, where $\kappa \colon \Omega\to \Lambda$ is the induced map between the invariant factors. 
Let~$(\mu_\lambda)_{\lambda\in \Lambda}$ denote a disintegration
of~$(\Omega,\mu)\to (\Lambda,\nu)$.
Then for almost every~$\lambda \in \Lambda$
the restriction of~$\tau$ to the fiber~$\kappa^{-1}(\{\lambda\})$
has $\mu_\lambda$-essentially countable range.
\end{proposition}
\begin{proof}
Since~$X\to Y$ is compact, by \cref{lem:compact_disintegration},
for almost every component~$\omega$ the
extension~$(G/\Gamma,\tau(\omega))\mapsto (G/L,\tau'(\kappa(\omega)))$
is compact. By \cref{thm:iso_central} it follows that~$L^0\subseteq  Z(G)$.
By \cref{fact:L_factorization}, we have~$L=L^0 K$ for some discrete countable
subgroup~$K$ containing~$\Gamma$.

Fix a generic~$\omega$ and let~$\lambda=\kappa(\omega)$.
Since $\pi(\omega,g\Gamma)=(\kappa(\omega),gL)$,
we have~$\tau(\omega) = \tau'(\lambda) \rho(\omega) \gamma(\omega)$
for some~$\rho(\omega)\in L^0 \subseteq  Z(G)$ and~$\gamma(\omega) \in K$
depending measurably on~$\omega$ (see \cref{prel:sections} on Borel sections
in nilmanifolds).
Since~$K$ is countable,
up to dividing~$\Omega$
into countably many measurable parts and restricting to one of them
(note that doing this does not affect the conclusion of the statement),
we may suppose~$\gamma(\omega)=\gamma$ for all~$\omega\in \Omega$ for
some fixed~$\gamma\in K$.
Let~$\tilde\tau(\lambda):=\tau'(\lambda) \gamma$. We just obtained a measurable factorization
\[
	\tau(\omega) = \rho(\omega) \tilde\tau(\lambda)
\]
with~$\rho(\omega)\in L^0 \subseteq  Z(G)$
and~$\tilde\tau(\lambda)\in G$ depending only on~$\lambda=\kappa(\omega)$.

Next, we need a factorization of the elements in the nilmanifold. It is convenient
to introduce some notation.
For $x\in X$, we write 
$$
x=(\omega(x),g(x)\Gamma), \qquad \lambda(x):=\kappa(\omega(x)),
$$
where~$g(x)\in G$ and $\omega(x)\in \Omega$, and hence $\lambda(x)\in \Lambda$,
depend measurably on $x$. We let 
$$
y(x):=(\lambda(x),g(x)L),
$$
so that $y(x)=\pi(x)$ for almost every $x\in X$.

We are now ready for the main argument.
Note that $L^0/(L^0\cap \Gamma)$ is a finite dimensional torus.
Moreover, since~$G/\Gamma$ is reduced and~$L^0\cap \Gamma$ is a subgroup of~$\Gamma$
normal in~$G$, we have~$L^0\cap \Gamma=\{1_G\}$, thus in fact~$L^0$ is a finite dimensional
torus.
In order to show that $\tau$ has essentially countable range
over a generic $\lambda$-fiber, we claim that it is enough to show that  for every set $B\subseteq  \Omega$
of positive measure and character $\chi$ of $L^0$ we have 
\begin{equation} \label{eq:BB_simple}
	\int_{\Lambda} \int_{B\times B} \1_{\{\chi(\rho(\omega))=\chi(\rho(\omega'))\}} \d \mu_\lambda(\omega) \d\mu_\lambda(\omega') \d\nu(\lambda) > 0.
\end{equation}
Indeed, using \cref{lem:countable_range}, this implies
that for almost every $\lambda \in \Lambda$ the function $\omega\mapsto \chi(\rho(\omega))$ has 
$\mu_\lambda$-essentially countable range. Since finitely many characters~$\chi$ separate
points in the finite dimensional torus $L^0$, this implies that for almost every $\lambda\in \Lambda$
the map $\omega\mapsto \rho(\omega)$
itself has $\mu_\lambda$-essentially countable range.
As a consequence, the same holds for $\tau(\omega)=\rho(\omega)\, \tilde{\tau}(\lambda)$. 

Thus fix arbitrarily a positive-measure set~$B$ and a character~$\chi$ of~$L^0$.
Let $F\colon G/\Gamma\to \C$ be
a corresponding nilcharacter, that is, a measurable map with modulus $1$  that satisfies 
$$
F(ux)=\chi(u) F(x) \qquad \text{for all } u\in L^0.
$$
The existence of such a function can be shown by making use of the  Borel sections mentioned above. Let
\[
	f(x) := \1_{B}(\omega(x)) F(g(x)\Gamma)
\]
Fix an arbitrary~$n\in \N$. Since $L^0 \subseteq Z(G)$ and $\rho(\omega(x))\in L^0$,
we have
\begin{align}
	T^n f(x) & = \1_{B}(\omega(x)) F(\tau(\omega(x))^n g(x)\Gamma) \nonumber \\
		& = \1_{B}(\omega(x)) \chi(\rho(\omega(x)))^n 
			F(\tilde\tau(\lambda(x))^n g(x) \Gamma). \label{prop:simple_nilbundle_countable_range_Tn}
\end{align}
It follows that 
$$
	 f(T^nx) \cdot \con{f}(x)= w_n(\omega(x)) \cdot F_n(\lambda(x),g(x)\Gamma) 
$$
where 
    $$
    w_n(\omega):= \1_{B}(\omega) \cdot \chi(\rho(\omega)^n), \qquad 
		F_n(\lambda, g\Gamma) :=F(\tilde{\tau}(\lambda)^n g \Gamma)
			\cdot \overline{F(g \Gamma)}.
$$
We will apply \cref{L:CondExp}  for $\Omega:=I(X)$, $\Lambda:=I(Y)$, $(M,m)=(G/\Gamma,m_{G/\Gamma})$, $N:=G/L$, $w:=w_n$, $F:=F_n$, respectively.
Since $w_n\in L^\infty(\Omega)$ and $F_n\in L^\infty(I(Y)\times (G/\Gamma))$, we have 
$$
	\E(T^n f \cdot \con{f} \mid Y)(x)
	 = \E(w_n\mid I(Y)) \cdot \E(F_n\mid Y)
 $$ 
  Since $|F_n|=1$ and $|\E(F_n\mid Y)|\leq 1 $, we deduce that 
  \begin{align}\label{E:idf}
  \abs{\E(T^n f \cdot \con{f} \mid Y)(x)}^2
\leq & |\E(w_n\mid I(Y))|^2= \Big|\int_B \notag \chi(\rho(\omega))^n\d\mu_\lambda(\omega)\Big|^2 \\
=&\int_{B\times B}
		\chi(\rho(\omega) \, \rho(\omega')^{-1})^n
		\d\mu_\lambda(\omega)\d\mu_\lambda(\omega').
\end{align}
 
On the other hand, since the extension $X\mapsto Y$ is compact and $f\neq 0$, we have 
$$ 
	 \limsup_{N\to\infty} \frac{1}{N}\sum_{n=1}^N
		\E\big(\abs{\E(T^n f \cdot \con{f} \mid Y)}^2\big) >0.
  $$
     Using  this positivity, equation \eqref{E:idf}, Fubini's theorem,  the bounded convergence theorem, and the geometric sum identity, we get that the preceding positive limsup is bounded above by 
     $$
     \int_{\Lambda}
		\int_{B\times B} \1_{\{\chi(\rho(\omega))=\chi(\rho(\omega'))\}}
			\d\mu_\lambda(\omega) \d\mu_\lambda(\omega')
		\d \nu(\lambda).  
 $$
 This proves the needed positivity and  completes the proof. 
\end{proof}

By \cref{countable_control_simple_nilbundle}, over a generic~$\lambda\in \Lambda$
the range of~$\tau\vert_{\kappa^{-1}(\lambda)}$ is countable.
As a consequence of Lusin-Novikov's Theorem, we can track this countable range
measurably as~$\lambda$ varies, as follows.

\begin{proposition} \label{cor:countable_range_selection}
Under the assumptions of \cref{countable_control_simple_nilbundle},
there is a countable collection~$\CB$ of measurable subsets
of~$\Omega$ such that~$\{\1_B\}_{B\in \CB}$ has dense span in~$L^2(\Omega,\mu)$,
and such that for every~$B \in \CB$ and for almost every~$\lambda$
the restriction of~$\omega \mapsto \tau(\omega)$ to~$B$
is $\mu_\lambda$-almost everywhere constant, with constant value~$\tau_{B}(\lambda)$
depending measurably on~$\lambda$.
\end{proposition}
\begin{proof}
Let 
$$
E:=\{
(\lambda,u)\in\Lambda\times G\colon
\tau_*\mu_\lambda(\{u\})>0\}.
$$
The set $E$ is Borel, since the map $(\lambda,u)\mapsto \int_\Omega \1_{\tau(\omega)=u}\d\mu_\lambda(\omega)$ is Borel, 
and every section
$E_\lambda:=\{u\in G\colon (\lambda,u)\in E\}$
is countable.
By the Lusin-Novikov  theorem~\cite[Theorem 18.10]{kechris2012classical},
there are Borel sets
$D_m\subseteq \Lambda$
and measurable maps
$$
\tau_m\colon D_m\to G,
\qquad m\in\N,
$$
whose graphs are pairwise disjoint
and whose union is $E$.
Let
$$
B_m:= \{
\omega\in\Omega\colon \kappa(\omega)\in D_m
\text{ and }
\tau(\omega)=\tau_m(\kappa(\omega))
\}.
$$
The sets $B_m$ are measurable
and pairwise disjoint.
For almost every $\lambda\in\Lambda$,
the measure $\tau_*\mu_\lambda$
is purely atomic, hence 
$
\mu_\lambda
\big(
\bigcup_{m\in\N}B_m
\big)
=1$, and it follows  from the disintegration formula that $\mu \big(\bigcup_{m\in\N}B_m\big)=1.$

Let~$\CA$ be a countable Boolean algebra generating the $\sigma$-algebra of~$\Omega$.
Consider the countable collection~$\CB := \{B_m\cap A\colon A\in \CA, m\in \N\}$.
Since~$\CA$ is a generating algebra and~$(B_m)_{m\in\N}$ is a measurable partition of $\Omega$, modulo nullsets, for every $A\in \CA$ we have 
$
\1_A=\sum_{m\in \N}\1_{A\cap B_m}$, with convergence in $L^2(\Omega,\mu)$. It follows that 
the span of~$\{\1_B\}_{B\in \CB}$ is dense
in~$L^2(\Omega,\mu)$. 
Lastly, for $B=B_m\cap A$, define
$\tau_B(\lambda):=\tau_m(\lambda)$
for $\lambda\in D_m$, and define $\tau_B$
to be an arbitrary constant outside $D_m$. Then $\tau_B$ is measurable, and for
almost every $\lambda\in\Lambda$,
the restriction of $\tau$ to $B$
is $\mu_\lambda$-almost everywhere constant,
with value $\tau_B(\lambda)$.
\end{proof}

\cref{cor:countable_range_selection} proved countable control for simple nilbundles
and an algebraic factor map. \cref{prop:nilbundle_factor_representation} shows
that any factor map~$\pi\colon X\to Y$ with~$X$ a nilbundle can be broken down
into countably many pieces where \cref{cor:countable_range_selection} applies,
obtaining the general case.

\begin{proposition} \label{prop:countable_range_selection_general}
Let~$X$ be a nilbundle and~$Y$ a factor with~$\pi\colon X\to Y$ compact.
With notation as in~\cref{prop:nilbundle_factor_representation}, define~$\tau$
as~$\tau\vert_{\Omega_i}=\tau_i$, taking values in the disjoint union $\bigsqcup_{i\in \N}G_i$. Then
there is a countable collection~$\CB$ of measurable subsets
of~$\Omega$ such that~$\{\1_B\}_{B\in \CB}$ has dense span in~$L^2(\Omega,\mu)$,
and such that for every~$B \in \CB$, and almost
every~$\lambda \in \Lambda$,
the restriction~$\tau\vert_B$
is $\mu_\lambda$-almost everywhere constant, with constant
value~$\tau_{B}(\lambda)$ depending measurably on~$\lambda$. 
\end{proposition}
\begin{proof}
We use notations as in~\cref{prop:nilbundle_factor_representation}.
It is enough to prove the following claim: for almost every~$\lambda\in \Lambda$, the restriction
of~$\tau$ to the fiber~$\kappa^{-1}(\{\lambda\})$ has essentially countable
range. In fact, this was the analogous claim in \cref{countable_control_simple_nilbundle}, in the case
of simple nilbundles. The construction of~$\CB$ from this claim
is the same as that carried out in the proof
of \cref{cor:countable_range_selection}. Let us prove the claim.
The key observation is that, since~$\Omega=\bigsqcup_{i\in \N}\Omega_i$ is a countable partition and~$\tau\vert_{\Omega_i}=\tau_i$,
the following facts are equivalent:
\begin{enumerate} [(i)]
\item for almost every~$\lambda\in \Lambda$ the restriction
of~$\tau$ to the fiber~$\kappa^{-1}(\{\lambda\})$
has $\mu_\lambda$-essentially countable range;
\item for every~$i$ and almost every~$\lambda\in \Lambda_{\iota(i)}$ the restriction
of~$\tau_i$ to the fiber~$\Omega_i\cap \kappa^{-1}(\{\lambda\})$
has $\mu_\lambda$-essentially countable range.
\end{enumerate}
In particular, we can fix an arbitrary~$i\in \N$ and consider the extension
\[
	\Omega_i\ltimes_{\tau_i} G_i/\Gamma_i
	\to \Lambda_{\iota(i)} \ltimes_{\tau'_{\iota(i)}} G'_{\iota(i)}/\Gamma'_{\iota(i)}.
\]
We have thus reduced to the case of a simple nilbundle.
In order to apply \cref{countable_control_simple_nilbundle}, it remains to
show that~$G'_{\iota(i)}/\Gamma'_{\iota(i)}$ can be replaced by its isomorphic
copy~$G_i/L_i$. This requires lifting~$\tau'_{\iota(i)}\in G'_{\iota(i)}$
to an element of~$G_i$.
Let~$\omega\in \Omega_i$, $\lambda=\kappa(\omega)$, and~$j=\iota(i)$.
From
$
	\varphi_i(\tau_i(\omega)) \varphi_i(g) \Gamma'_j
	= \tau'_j(\lambda) \varphi_i(g) \Gamma_j'
$
for all~$g\in G_i$, surjectivity of~$\varphi_i$,
and the fact that~$\Gamma_j'$ contains no non-trivial normal subgroups of~$G'_j$, it
follows that~$\varphi_i(\tau_i(\omega))=\tau'_j(\lambda)$.
By fixing a section, we obtain a $\Lambda_j$-measurable element $\tilde \tau_i(\lambda)\in G_i$
such that~$\tau_i(\omega) \ker\varphi_i = \tilde \tau_i(\lambda) \ker\varphi_i$.
Since~$\ker\varphi_i$ is normal and~$\ker\varphi_i\subseteq  \varphi_i^{-1}(\Gamma'_j)=L_i$,
we see that~$\tau_i(\omega)$ and~$\tilde \tau_i(\lambda)$ have the same action
on~$G_i/L_i$. In particular, we can replace the previous extension by
\[
	\Omega_i\ltimes_{\tau_i} G_i/\Gamma_i
	\to \Lambda_{\iota(i)} \ltimes_{\tilde \tau_i} G_i/L_i,
\]
for which \cref{countable_control_simple_nilbundle} now applies.
\end{proof}

We are now ready to prove the $\nilbor$ property
for compact extensions of nilbundles.

\begin{proposition} \label{beautiful_nilbundles}
Let~$X$ be a nilbundle and~$Y$ a factor with~$X\to Y$ compact. Then there
is a countable set~$\mathcal F\subseteq  L^\infty(X)$, dense in~$L^2(X)$,
with the $\nilbor$ property.
\end{proposition}

\begin{proof}
The argument is a generalization of \cref{ex:countable_control_I},
with characters replaced by nilcharacters,
and using the countable control guaranteed by \cref{prop:countable_range_selection_general}.

We use the notation introduced in \cref{prop:nilbundle_factor_representation}
and \cref{prop:countable_range_selection_general}.
By \cref{lem:compact_disintegration} and \cref{thm:iso_central} we have $L^0_i\subseteq Z(G_i)$. The system is spanned by functions of the form~$f=\1_A F$,
where~$A=\Omega_i\cap B$ for some~$i\in \N$ and~$B\in \CB$,
and~$F$ is a nilcharacter of~$G_i/\Gamma_i$, say with character~$\chi$, that is,
a continuous map~$G_i/\Gamma_i\to \C$ satisfying
\[
	F(\ell g\Gamma_i) = \chi(\ell) F(g\Gamma_i)
\]
for every~$\ell\in Z(G_i)$. Since the $\nilbor$ property is preserved under finite linear combinations, it suffices to establish it for each such function $f$.

Let $\tau_B\colon \Lambda\to \bigsqcup_{i\in \N} G_i$ be as in \cref{prop:countable_range_selection_general}. 
Replacing $\tau_B(\lambda)$ by $1_{G_i}$ whenever its value lies outside $G_i$, we obtain a measurable map $\tau_B\colon \Lambda \to G_i$ that still satisfies $\tau_B(\kappa(\omega))=\tau_i(\omega)$ for almost every $\omega\in A$. For~$(\omega,g\Gamma_i)\in \Omega_i\times G_i/\Gamma_i$ we have
\[
	f(\omega,g\Gamma_i) = \1_A(\omega) F(g\Gamma_i),
\]
while~$f=0$ outside~$\Omega_i \times G_i/\Gamma_i$.
For simplicity, we assume~$L_i=L_i^0 \Gamma_i$. The general case can be obtained
exactly as done in the proof of \cref{prop:key_lem_nilsystem}.
As in the proof of \cref{countable_control_simple_nilbundle}, we write
$$
x=(\omega(x),g(x)\Gamma_i), \qquad \lambda(x):=\kappa(\omega(x)),
\qquad y(x):=(\lambda(x),g(x)L_i),
$$
and we have a factorization
\[
	g(x)\Gamma_i = \ell(x)\,  \sigma(y(x)) \Gamma_i,
\]
for some measurable maps $\ell\colon \Omega_i\times G_i/\Gamma_i \to L^0_i$
and~$\sigma\colon \Lambda_{\iota(i)} \times G'_{\iota(i)}/\Gamma'_{\iota(i)} \to G_i$.
For convenience, we extend~$\ell$ and~$\sigma$ to measurable maps defined on
the whole~$X$ and~$Y$
respectively, in an arbitrary way: due to the presence of~$\1_A=\1_{B\cap \Omega_i}$ in
the following argument, the specific extension is irrelevant.

By definition, for~$x=(\omega(x),g(x)\Gamma_i)\in \Omega_i\times G_i/\Gamma_i$ we have
\[
	f(x)=\1_A(\omega(x)) \, \chi(\ell(x)) \, F(\sigma(y(x)) \Gamma_i).
\]
By definition of~$B$ in \cref{prop:countable_range_selection_general},
for almost every~$\omega\in A$ we have~$\tau_i(\omega)=\tau_B(\kappa(\omega))$,
thus
\[
	T^n f(x) = \1_A(\omega(x)) \, \chi(\ell(x))\,  F(\tau_B(\lambda(x))^n \, \sigma(y(x)) \Gamma_i).
\]
By combining these identities, we obtain
\[
	\abs{T^n f(x)-f(x)} \leq \abs{F(\tau_B(\lambda(x))^n \sigma(y(x)) \Gamma_i)
		- F(\sigma(y(x)) \Gamma_i)}.
\]
The right hand side is $Y$-measurable
(it depends only on~$y$ and~$\lambda=\lambda(y)$), thus
\begin{equation}\label{E:esty}
	\E(\abs{T^n f-f}^2\mid Y)(y)
		\leq \abs{F(\tau_B(\lambda(y))^n \sigma(y) \Gamma_i)
			- F(\sigma(y) \Gamma_i)}^2.
\end{equation}
Fix an arbitrary~$\ve>0$, let~$\delta:=\sqrt{\ve}$, and define
\[
    \CM(y,n)
    :=
    \1_{[0,\delta)}\!\left(
    \abs{
    F(\tau_B(\lambda(y))^n\sigma(y)\Gamma_i)
    -F(\sigma(y)\Gamma_i)}
    \right).
\]
By \cref{lem:F_nilbohr},
$\CM$ is a $Y$-measurable varying $\nilbor$ sequence.
If $\CM(y,n)=1$, the estimate \eqref{E:esty} gives
$$
    \E(\abs{T^n f-f}^2\mid Y)(y)\leq \ve.
$$
This completes the proof.
\end{proof}


\subsection{Killing denominators in nilbundles}

For convenience, we state the following elementary observation.

\begin{lemma} \label{measurable_stabilization}
Let~$(\Omega,\mu)$ be a measure space and~$\CM$ a countable set
endowed with the discrete topology. Let~$(f_n \colon \Omega\to \CM)_{n\in \N}$
be a sequence of measurable maps such that
for $\mu$-almost every~$\omega$ the sequence~$(f_n(\omega))_{n\in \N}$
eventually stabilizes. Then~$f:=\lim_{n\to\infty} f_n$ exists and it is measurable.
\end{lemma}
\begin{proof}
For every~$M\in \CM$ the set
$
	\{\omega\colon f(\omega)=M\}
							=\bigcup_{N\geq 1}\bigcap_{n\geq N} \{\omega\colon f_n(\omega)=M\}
$
is measurable. Since~$\CM$ is countable discrete, this makes~$f$ measurable.
\end{proof}

\begin{remark} \label{rem:nilbundle_alternative_notation}
A nilbundle has the form~$\bigsqcup_{i\in \N} \Omega_i \ltimes_{\tau_i} G_i/\Gamma_i$.
Alternatively, it can be written as~$\bigsqcup_{\omega\in \Omega}
		(G_{\omega}/\Gamma_{\omega},\tau_{\omega})$,
provided that the map~$\omega\mapsto (G_\omega,\Gamma_\omega)$ is trivialized
by a countable measurable partition of~$\Omega$,
and~$\omega\mapsto \tau_\omega$ is measurable with values~$\tau_\omega \in G_\omega$.
In the proof of the next result we will use this latter notation.
\end{remark}

The next result is a non-ergodic generalization of \cref{prop:Znil_Ypro}.
\begin{proposition} [Killing the Denominators in the Non-Ergodic Case] \label{killing_nonergodic}
Let~$(X,\mu,T)$ be a $\CZ_s$-system. Let~$Y$ and~$Z$ be factors of~$X$, with~$Z$ a nilbundle-factor
of~$X$. Then there is a nilbundle-factor~$\tilde Z$ such that~$Z\subseteq  \tilde Z\subseteq  X$,
and such that~$Y\vee I(X)$ and~$\tilde Z$ have the good intersection property. 

Furthermore, letting $\Omega:=I(X), \Lambda:=I(Y)$, $(\mu_\omega)_{\omega\in \Omega}$ be the ergodic decomposition of $\mu$,  and $\kappa\colon \Omega\to \Lambda$ be the induced map, 
given an increasing sequence of nilbundle factors $(Y_n)_{n\in\N}$ with  $Y=\bigvee_{n\in \N} Y_n$ and $I(Y_n)=I(Y)$, 
the factor $\tilde{Z}$ can be chosen so that for almost every $\omega\in \Omega$ the factors $\tilde{Z}_\omega$ and $(Y_n)_{\kappa(\omega)}$ have the good intersection property with respect to $\mu_\omega$ for all sufficiently large $n$ depending on $\omega$. 
\end{proposition}
\begin{remark}
A sequence $(Y_n)$ as in the statement exists by \cref{T:NonErgodicStructure},
since $Y$, being a factor of the $\CZ_s$-system $X$,
is itself a $\CZ_s$-system.
\end{remark}
\begin{proof}
Our plan is to use \cref{prop:Znil_Ypro_explicit} on each ergodic component and then combine the resulting factors measurably to construct $\tilde{Z}$. Throughout, a countably valued family of groups, subgroups, homomorphisms,
or nilmanifolds indexed by $\Omega$ is called measurable
if each level set is measurable.
Equivalently, the family is {\em trivialized by a countable measurable partition
of $\Omega$,} meaning that it is constant on each cell of the partition.

{\em  Preparation.}
We carry out some preparatory work to apply \cref{prop:Znil_Ypro_explicit}
on almost every ergodic component.
To work over a common base $\Omega$,
we replace $Z$ with $Z\vee I(X)$,
which is a nilbundle by \cref{joint_nilbundle_nilbundle}.
We may therefore assume that $I(X)\subseteq  Z$.
We use the sequence $(Y_n)$ from the statement,
whose common nilbundle base is $I(Y)$.
By \cref{joint_nilbundle_nilbundle},
both $Z\vee Y_n$ and $Y_n\vee I(X)$
are nilbundles with base $\Omega$.

By \cref{lem:reduced_nilbundle} and 
applying inductively the change of base points and homomorphism selection in the proof of
\cref{prop:nilbundle_factor_representation}, we can assume that the factor maps $Z\vee Y_{n+1}\to Z\vee Y_n$ are induced on almost every ergodic component by  surjective continuous group homomorphisms
$$
\psi_{n,\omega}\colon G_{n+1,\omega} \to G_{n,\omega}.
$$
With notation
as in the preceding remark,
we have representations  
\begin{align*}
	& Z\vee Y_n
	= \bigsqcup_{\omega\in \Omega}
		(G_{n,\omega}/\Gamma_{n,\omega},\tau_{n,\omega}), \\
	& Y_n \vee I(X) = \bigsqcup_{\omega \in \Omega}
		(G_{n,\omega}/L_{n,\omega},\tau_{n,\omega}), \\
	& Z = \bigsqcup_{\omega\in \Omega}
		(G_{1,\omega}/M_{\omega},\tau_{1,\omega}).
\end{align*}
We let $\varphi_{1,\omega}=\text{id}_{G_{1,\omega}}$ and for $n\geq 2$ we let 
$$
\varphi_{n,\omega}:=\psi_{1,\omega}\circ \cdots \circ \psi_{n-1,\omega}\colon G_{n,\omega}\to G_{1,\omega}.
$$
By \cref{hom_nilmanifolds} and the selection argument in the proof of \cref{prop:nilbundle_factor_representation}, for every $n\in \N$
the data
$$
	G_{n,\omega},\Gamma_{n,\omega},L_{n,\omega},M_{\omega},
		\varphi_{n,\omega} 
$$
are  trivialized by a countable measurable partition of~$\Omega$
(note that~$n$ here denotes
the index of the nilbundle in the inverse limit, not the elements
of a countable partition of~$\Omega$, which are instead left unlabeled to not
clutter notation). 
On each cell the rotation $\omega\mapsto \tau_{n,\omega}$ is measurable and the countability of the groups and maps will make subsequent subgroup choices measurable. Moreover since $
\psi_{n,\omega}(L_{n+1,\omega})
\subseteq  L_{n,\omega}
$
we get that $$
\varphi_{n+1,\omega}(L_{n+1,\omega})
=
\varphi_{n,\omega}\!\left(\psi_{n,\omega}(L_{n+1,\omega})\right)
\subseteq 
\varphi_{n,\omega}(L_{n,\omega}),
$$ so the sequence 
$\varphi_{n,\omega}(L_{n,\omega})$ is decreasing in $n$. Also, as in the first paragraph of the proof of \cref{prop:Znil_Ypro_explicit}, these groups are $\Gamma_{1,\omega}$-rational with normal identity components.  

\medskip

{\em Use of \cref{prop:Znil_Ypro_explicit}.}
The previous discussion guarantees that 
for almost every $\omega\in \Omega$ we can apply  \cref{prop:Znil_Ypro_explicit} to the sequence
$(G_{n,\omega}/\Gamma_{n,\omega})$,
with successive maps induced by $\psi_{n,\omega}$,
and with $M:=M_\omega$ and $L_n:=L_{n,\omega}$.
Its compatibility hypotheses   follow from
$\psi_{n,\omega}(\Gamma_{n+1,\omega})\subseteq \Gamma_{n,\omega}$
and $\psi_{n,\omega}(L_{n+1,\omega})\subseteq  L_{n,\omega}$.
We get that there exist closed subgroups $\tilde{M}_\omega\subseteq M_\omega$
of finite index and an integer $N_\omega$ such that for every $n\geq N_\omega$ we have 
\begin{equation}\label{E:commute}
\Gamma_{n,\omega}\subseteq \varphi^{-1}_{n,\omega}(\tilde{M}_\omega), \qquad 
\varphi^{-1}_{n,\omega}(\tilde{M}_\omega) L_{n,\omega}=L_{n,\omega} \varphi^{-1}_{n,\omega}(\tilde{M}_\omega). 
\end{equation}
The inclusion $\tilde{M}_\omega\subseteq M_\omega$ will give an extension of $Z_\omega$, the lattice inclusion makes it a factor of $(Z\vee Y_n)_\omega$, and the subgroup equality gives good intersection 
with $(Y_n)_{\kappa(\omega)}$ by \cref{lem:LM_expected}. To combine these factors into a nilbundle we must make some choices measurably; this is the goal of the next step. 

\medskip

{\em Measurable selection.}
The collection~$\{\varphi_{n,\omega}(L_{n,\omega})^0\}_{\omega\in \Omega,n\in \N}$
is countable, for every~$n$ the map~$\omega\mapsto \varphi_{n,\omega}(L_{n,\omega})^0$
is measurable
(since it is trivialized by a countable measurable partition),
and for almost every~$\omega$ the sequence~$n\mapsto \varphi_{n,\omega}(L_{n,\omega})^0$
eventually stabilizes by \cref{fact:Ln_factorization}. Thus, by \cref{measurable_stabilization}, 
the limit $Q_\omega$ is measurable. Choose $N_\omega$ to be the least integer such that $\varphi_{n,\omega}(L_{n,\omega})^0=Q_\omega$ for all $n\geq N_\omega$. Then $\omega \mapsto N_\omega$ is measurable. 

Let $\Lambda_{1,\omega}$ be the first subgroup constructed in the proof
of \cref{fact:Ln_factorization}, applied to
$(\varphi_{n,\omega}(L_{n,\omega}))_{n\geq 1}$
using the stabilization index $N_\omega$.
 By its construction in the proof of \cref{fact:Ln_factorization}, the map $\omega\mapsto\Lambda_{1,\omega}$ is countably valued and measurable.
We define $\tilde\Gamma_\omega$ by \eqref{E:tildeGamma}
in the proof of \cref{prop:Znil_Ypro_explicit}, with $\Gamma:=\Gamma_{1,\omega}$ and $\Lambda_1:=\Lambda_{1,\omega}$.
Since $\Lambda_{1,\omega}$ and $\Gamma_{1,\omega}$ are constant on each part
of a common countable measurable partition of $\Omega$,
the intersection in \eqref{E:tildeGamma} is also constant on each part, hence it is trivialized by this partition.
Similarly, the group 
$$
	\tilde{M}_{n,\omega}:=\gen{M_\omega^0,\tilde \Gamma_\omega,
		\varphi_{n,\omega}(L_{n,\omega})\cap M_\omega}
$$
varies in a countable set, it varies measurably for each~$n$, and it stabilizes
for almost every~$\omega$. Indeed, it decreases with $n$ and contains the fixed finite-index subgroup $M^0_\omega\tilde{\Gamma}_\omega$ of $M_\omega$. Thus again by \cref{measurable_stabilization}, and the same argument for the least stabilization index we previously used,  we can 
find measurable maps~$\tilde M_\omega$ and~$\tilde N_\omega$
such that~$\tilde M_{n,\omega} = \tilde M_\omega$ for all~$n\geq \tilde N_\omega$.
Up to replacing~$N_\omega$ and~$\tilde N_\omega$ by their maximum (which is measurable)
we may assume that $N_\omega=\tilde N_\omega$.

\medskip 
{\em Construction of the nilbundle factor and good intersection property.} 
 Note that, by the previous construction,
both~$\tilde M_\omega$ and~$N_\omega$
are trivialized by a countable measurable partition, thus  
\[
	\tilde Z := \bigsqcup_{\omega\in \Omega}
		(G_{N_\omega,\omega}/\varphi^{-1}_{N_\omega,\omega}(\tilde{M}_\omega),\tau_{N_\omega,\omega})
\]
is a nilbundle. For each $r$, the families $G_{r,\omega}$, $\Gamma_{r,\omega}$,
and $\varphi_{r,\omega}^{-1}(\tilde{M}_\omega)$
are trivialized by a common countable measurable partition of $\{N_\omega=r\}$. The natural quotient maps therefore define, over $\{N_\omega=r\}$,
a measurable factor map $Z\vee Y_r\to\tilde{Z}$.
Composing with $X\to Z\vee Y_r$ and combining over $r$,
we obtain a global measurable factor map $X\to\tilde{Z}$.
Since $\tilde{M}_\omega\subseteq  M_\omega$, this factor contains $Z$.

 It remains to verify the good intersection properties. We start with the componentwise property. 
The proof of \cref{prop:Znil_Ypro_explicit}
shows that for almost every ergodic component~$\omega$
the system~$\tilde Z_\omega$ is a factor of~$X_\omega$ and  contains~$Z_\omega$. 
We  verify that $\tilde Z_\omega$ and $(Y_n)_{\kappa(\omega)}$
have the good intersection property for every $n\geq N_\omega$; 
the corresponding property with $Y_{\kappa(\omega)}$
will then follow by martingale convergence.
Fix $\omega$ outside the exceptional nullset and $n\geq N_\omega$.
By construction, we have the representations
$$
\tilde Z_\omega
=
G_{n,\omega}/\varphi_{n,\omega}^{-1}(\tilde M_\omega),
\qquad
(Y_n)_{\kappa(\omega)}
=
G_{n,\omega}/L_{n,\omega}.
$$
Combining  \eqref{E:commute} and \cref{lem:LM_expected} we get that $\tilde Z_\omega$
and $(Y_n)_{\kappa(\omega)}$ have the good intersection property
with respect to $\mu_\omega$.
 Since $(Y_n)_{\kappa(\omega)}$ increases to $Y_{\kappa(\omega)}$,
martingale convergence implies that $\tilde Z_\omega$
and $Y_{\kappa(\omega)}$ also have the good intersection property.
Since both $\tilde{Z}$ and $Y\vee I(X)$ contain $I(X)$,
\cref{lem:good_transfer}
transfers the componentwise good intersection property
to the global factors $\tilde Z$ and $Y\vee I(X)$.
This proves both asserted good intersection properties. 
\end{proof}


\subsection{Concluding the proof}
We are finally ready to prove \cref{thm:beautiful_main}, which we restate again for reader's convenience. 
\begin{restatedbeautiful}
Let~$X$ be a $\CZ_s$-system and~$X\to Y$ a compact extension.
Then there is a countable subset~$\CF\subseteq  L^\infty(X)$, dense in~$L^2(X)$,
with the following property:
for every~$f\in \CF$ and~$\ve>0$
there is a measurable varying $\nilbor$ sequence~$\CM\colon Y\times \Z\to \{0,1\}$ satisfying
\[
	\CM(y,n)=1 \qquad\implies\qquad \E(\abs{T^n f - f}^2\mid Y)(y) \leq\ve
\]
for almost every $y\in Y$ and every $n\in \N$.
\end{restatedbeautiful}

\begin{proof}
Let~$X$ be a $\CZ_s$-system and~$X\to Y$ a compact extension.
Our first goal is to obtain nilbundle approximants of $X$ that have good intersection with $Y$; this will then enable us to use \cref{beautiful_nilbundles}.

To ease notation, let $I:=I(X)$, $\Omega:=I(X)$, $\Lambda:=I(Y)$, and let $\kappa\colon \Omega\to \Lambda$ be the induced map.
Since $Y$ is also a $\CZ_s$-system, by \cref{T:NonErgodicStructure} it is an inverse limit of nilbundle-factors~$Y_m$ of $Y$, hence, $Y=\bigvee_{m\in \N} Y_m$ and $I(Y_m)=\Lambda$. We let 
$$
Y':=Y\vee I, \qquad Y'_m:=Y_m\vee I.
$$

By \cref{T:NonErgodicStructure},
the system~$X$ is an inverse limit of nilbundle-factors~$X_k$ with base $\Omega$.
Fix $k\in \N$. We  apply \cref{killing_nonergodic} to $X_k$ and $Y$, with the prescribed sequence $(Y_m)$. We get that there exist nilbundle factors $X'_k$ such that $X_k\subseteq X'_k\subseteq X$, the factors $X'_k$ and $Y'$ have good intersection, and for almost every $\omega$ the factors $(X'_k)_\omega$ and $(Y_m)_{\kappa(\omega)}$ have good intersection for all sufficiently large $m$ depending on $\omega$. Let 
$$
W_k:=X'_k\wedge Y', \qquad W_{k,m}:=X'_k\wedge Y'_m.
$$
 Then we have \begin{equation}\label{E:Wkm} (W_k)_\omega=\bigvee_{m\in\N}(W_{k,m})_\omega
 \end{equation}
 for almost every $\omega$, by \cref{lem:good_wedge_limit} (applied for $\CI:=I(X)$, $\CA:=X'_k$, $\CB:=Y'$, $\CB_n:=Y'_n$).
 
A crucial observation is that due to the nilmanifold structure of the ergodic components of nilbundles, for almost every $\omega$ the increasing sequence of factors $(W_{k,m})_\omega$ stabilizes  (as a function of $m$).
Indeed, for fixed~$k$ and generic~$\omega$,
the sequence~$m\mapsto (W_{k,m})_\omega$ is an increasing sequence
of factors of the ergodic nilsystem~$(X'_k)_\omega$,
and thus by \cref{cor:nilsystems_stabilize}
it must stabilize.

Combining the previous observation with  \eqref{E:Wkm} gives that  $(W_{k,m})_\omega$  stabilizes to $(W_k)_\omega$ as $m$ increases. 
Consider the sets 
\begin{equation}\label{E:Akm}
A_{k,m}:=\{x\in X\colon  (W_{k,m})_{\omega(x)}=(W_k)_{\omega(x)}\},
\end{equation}
where $x\mapsto \omega(x)$ is the invariant-factor map. These sets are easily shown to be measurable, since equality on components can be tested by equality of the conditional expectation operator on a fixed countable family of bounded functions. The previous stabilization property gives that 
$$
X=\bigcup_{m\in\N}A_{k,m},
$$
modulo nullsets.
Let 
\begin{equation}\label{E:Zkm}
Z_{k,m}:=Y'_m \vee \sigma( L^\infty(X'_k)\cap L^\infty(A_{k,m}))
\end{equation}
where $L^\infty(A)$ denotes the bounded measurable functions on $X$ that vanish almost everywhere outside $A$ and $\sigma(\CF)$ denotes the factor generated by the functions in $\CF$. 

\smallskip 

{\em Claim 1: $Z_{k,m}$ is a nilbundle-factor of $X$.} The restriction of $Z_{k,m}$ to $A_{k,m}$  coincides with that of  $X'_k\vee Y_m$, while its restriction to $A^c_{k,m}$ coincides with that of $Y'_m$.
Since $X'_k\vee Y_m=X'_k\vee Y'_m$ and  both $X'_k, Y'_m$ are nilbundle-factors of $X$, closure under finite joins (\cref{joint_nilbundle_nilbundle})  
gives that $X'_k\vee Y_m$ is also a nilbundle-factor of $X$. 
 Finally, gluing these nilbundles along the measurable partition of
the invariant base induced by $A_{k,m}\in I(X)$ shows that
$Z_{k,m}$ is also a nilbundle-factor of $X$. 

\smallskip 

{\em Claim 2: $Z_{k,m}\wedge Y=Y_m$ and $Z_{k,m}$ and $Y$ have good intersection.} It suffices to show that 
\begin{equation}\label{E:condexp}
\E(h\mid Y)=\E(h\mid Y_m) \quad \text{for every } h\in L^\infty(Z_{k,m}).  
\end{equation}
Indeed, 
$Y_m\subseteq (Z_{k,m}\wedge Y)$ by construction. Conversely, for $h\in L^\infty(Z_{k,m}\wedge Y)$, the identity gives $h=\E(h\mid Y)=\E(h\mid Y_m),$ so $h$ is $Y_m$-measurable. Therefore, $Z_{k,m}\wedge Y=Y_m$. Substituting this equality into \eqref{E:condexp} yields the asserted good intersection property. So it remains to prove \eqref{E:condexp}. (The key point is that on $A_{k,m}$ conditioning on $W_k$ gives the same result as conditioning  on $W_{k,m}$.)  We first handle functions supported on $A_{k,m}$. For $f\in L^\infty(X'_k)$ we have 
\begin{equation}\label{E:condexpAkm}
\E(\1_{A_{k,m}} f\mid Y')=\1_{A_{k,m}} \E(f\mid W_k)= \1_{A_{k,m}} \E(f\mid W_{k,m}),
\end{equation}
where the first equality holds since $A_{k,m}\in I\subseteq Y'$ and the good intersection property of $X'_k$ and $Y'$, while the second equality follows from \cref{lem:relative_conditioning} applied to $W_k$ and $W_{k,m}$; both factors contain $I$, so their conditional expectations can be computed componentwise and for $x\in A_{k,m}$ we get by \eqref{E:Akm} that $(W_{k,m})_{\omega(x)}=(W_k)_{\omega(x)}$. Next, we first prove that 
$$
\E(h|Y')=\E(h|Y'_m), \qquad \text{for  every } h\in L^\infty(Z_{k,m}). 
$$
 By the definition of $Z_{k,m}$ in \eqref{E:Zkm} it suffices to verify this  when $h=g$ and when $h=g\cdot \1_{A_{k,m}} f$, where $g\in L^\infty(Y'_m)$ and $f\in L^\infty(X'_k)$ are  arbitrary. For $h=g$ the identity is immediate. For $h=g\cdot \1_{A_{k,m}} f$, multiplying 
\eqref{E:condexpAkm} by $g$ and using that $g$ is $Y'$-measurable (since $Y'_m$ is factor of $Y'$) we get 
\begin{align}\label{E:Y'Y'm}
\notag \E(h\mid Y')&=g \E(\1_{A_{k,m}} f\mid Y') =g \, \1_{A_{k,m}} \E(f\mid W_{k,m})\\ 
&=\E(g\, \1_{A_{k,m}}\E(f\mid W_{k,m})\mid Y'_m)=\E(\E(h\mid Y')\mid Y'_m)=\E(h\mid Y'_m),   
\end{align}
where the first equality follows since $g$ is $Y'_m$ measurable and $Y'_m$ is a factor of $Y'$, the second follows from \eqref{E:condexpAkm}, the third because all functions are $Y'_m$-measurable, the fourth by using the first two identities, and the last identity since $Y'_m$ is a factor of $Y'$.  

To remove the primes in the previous identity, we first note that since $I(Y)$ is a factor of $Y_m$ and $Y'_m=Y_m\vee I(X)$,  we have 
\begin{equation}\label{E:YY'm}
\E(g\mid Y)=\E(g\mid Y_m), \quad \text{ for every } g\in L^\infty(Y'_m).
\end{equation}
Indeed, it suffices to verify the identity for $g=g_1 g_2$, where $g_1\in L^\infty(I(X))$ and $g_2\in L^\infty(Y_m)$. 
To this end, note that 
$$
\E(g_1\mid Y)= \E(\E(g_1\mid Y)\mid Y_m)=\E(g_1\mid Y_m),
$$
where the first equality follows since the $T$-invariance of $g_1$ implies that  $\E(g_1\mid Y)$ is $I(Y)$-measurable, hence $Y_m$-measurable, and the second   since $Y_m$ is a factor of $Y$.
Using this and that $g_2$ is $Y_m$-measurable hence $Y$-measurable, we deduce that 
$$
\E(g_1 g_2\mid Y)=g_2\E(g_1\mid Y)=g_2\E(g_1\mid Y_m)=\E(g_1g_2\mid Y_m), 
$$ 
establishing  \eqref{E:YY'm} for $g=g_1g_2$. 
Finally, note that for  $h\in L^\infty(Z_{k,m})$ we have 
\begin{align*}
\E(h\mid Y)&=\E(\E(h\mid Y')\mid Y)=
\E(\E(h\mid Y'_m)\mid Y)\\ &=\E(\E(h\mid Y'_m)\mid Y_m)
=\E(h\mid Y_m),
\end{align*}
where the first equality uses that $Y$ is a factor of $Y'$,  the second uses \eqref{E:Y'Y'm}, the third uses \eqref{E:YY'm} for $g:=\E(h\mid Y'_m)$, and the last uses that $Y_m$ is a factor of $Y'_m$.
This completes the proof of \eqref{E:condexp}.

\smallskip

 {\em Claim 3: There exists a countable  subset of $L^\infty(Z_{k,m})$, dense in $L^2(Z_{k,m})$, that has the $\nilbor$ property.}

By Claim 1, the factor~$Z_{k,m}$ is a nilbundle.
By Claim 2, the factors~$Z_{k,m}$ and~$Y$ have the good intersection
property. In particular,
for any $Z_{k,m}$-measurable function~$f$ and every~$n\in \N$ we have
\begin{equation} \label{eq:good_nilbohr_final}
    \E(\abs{T^n f - f}^2\mid Y) = \E(\abs{T^n f - f}^2\mid Z_{k,m}\wedge Y).
\end{equation}
This thus reduces the problem to the extension
of nilbundles~$Z_{k,m}\to Z_{k,m}\wedge Y$.
From the properties of good intersections, \cref{L:compactfactor},
the extension~$Z_{k,m}\to Z_{k,m}\wedge Y$ is compact.
In particular \cref{beautiful_nilbundles}
applies to this extension, giving a dense set of functions in~$L^\infty(Z_{k,m})$
with the $\nilbor$ property with respect to~$Z_{k,m}\wedge Y$, and
thus with respect to~$Y$ by~\eqref{eq:good_nilbohr_final}.

For $f\in L^\infty(X_k)$ the functions $\1_{A_{k,m}} f$ belong to $L^\infty(Z_{k,m})$ and converge to $f$ in $L^2(\mu)$. Since $\bigvee_{k\in \N}X_k=X$, the union of the countable dense families obtained in Claim 3 is dense in $L^2(X)$. 
\end{proof}

From the good intersection property,
we can replace ``$\CZ_s$'' by ``$\CZ_s$-measurable'', in the following sense.

\begin{corollary} \label{cor:beautiful_main}
Let~$X$ be a system (not necessarily ergodic) and~$X\to Y$ a compact extension.
Let~$\CZ_s(X)$ be the usual Host-Kra factor.
Then there is a countable subset~$\CF\subseteq  L^\infty(\CZ_s(X))$, dense in~$L^2(\CZ_s(X))$,
such that for every~$f\in \CF$ and~$\ve>0$ there is a measurable varying
$\nilbor$ sequence~$\CM:Y\times \Z\to \{0,1\}$ satisfying
\[
	\CM(y,n)=1 \qquad\implies\qquad \E(\abs{T^n f - f}^2\mid Y)(y) \leq\ve.
\]
\end{corollary}
\begin{proof}
If~$f$ is $\CZ_s(X)$-measurable, so is~$\abs{T^n f - f}$, and by \cref{L:ZsGood} 
we have
\[
	\E(\abs{T^n f - f}^2\mid Y) = \E(\abs{T^n f - f}^2\mid \CZ_s(X) \wedge Y).
\]
Moreover, by \cref{L:compactfactor},
the extension~$\CZ_s(X)\to \CZ_s(X) \wedge Y$ is compact. 
\cref{thm:beautiful_main} applies.
From the proof of \cref{thm:beautiful_main}, we immediately see that
the resulting $\nilbor$-bundle is $s$-step.
\end{proof}



	\section{Partial characteristic factors}

In this section we prove \cref{cor:partial_factor_use}.
We will use dual functions and dual sequences, defined
in \cref{sec:prel_HK}, and nilsequences,
defined in \cref{sec:prel_nil_definitions}.
We start with some simple preliminary results.

\begin{lemma}\label{L:1}
	Let $(\CN(n))_{n\in\N}$ be an $r$-step nilsequence. Then for every $\varepsilon>0$ there exists a sequence $(\CD(n))_{n\in\N}$, which is a finite linear combination of dual sequences of level $(r+1)$, such that 
	$$
	\norm{\CN-\CD}_\infty\leq \varepsilon
	$$
\end{lemma}
\begin{proof}
	We have that  $\CN(n)=F(b^n\Gamma)$ for some $r$-step nilmanifold $X=G/\Gamma$, $b \in G$ ergodic, and $F\in C(X)$. It is known (\cref{sec:prel_HK}) that there exists $\tilde{F}\in C(X)$, a linear combination of dual functions of level $(r+1)$,  such that  $\norm{F-\tilde{F}}_{C(X)}\leq \varepsilon$
	Then 
	$
	\CD(n):=\tilde{F}(b^n\Gamma)$, $n\in\N$, is a linear combination of dual sequences of level $(r+1)$ and 
    satisfies
	$
	\norm{\CN-\CD}_\infty\leq \varepsilon
	$
\end{proof}

We immediately deduce the following.
\begin{corollary}\label{C:1}
	Let $N_k\to \infty $ and $(A(n))$ be a bounded sequence such that 
	$$
	\lim_{k\to \infty} \E_{n\in [N_k]}\, \CD(n)\cdot A(n)=0
	$$
	for every dual sequence $(\CD(n))_{n\in\N}$ of level $r$. Then 
	$$
	\lim_{k\to \infty} \E_{n\in [N_k]}\, \CN(n)\cdot A(n)=0
	$$
	for every $(r-1)$-step nilsequence $(\CN(n))_{n\in\N}$.
\end{corollary}
We use the previous corollary to prove the following intermediate result. 
\begin{lemma}\label{L:2}
	Let  $d,\ell, r\in \N$ and  $p_1,\ldots, p_\ell\in \Z[t]$ have degree at most $d$. Suppose that either $\ell=1$  and $p_1$ is non-constant, or $\ell\geq 2$ and the  polynomials are   pairwise independent  and $\deg(p_\ell)\geq 2$.  Then there exists $s\in \N$, dependent only on  $d,\ell,r$,  such that the following holds: If $(X,\mu,T_1,\ldots, T_\ell)$ is a system  and $f_1,\ldots,f_\ell \in L^\infty(\mu)$ are $1$-bounded and   $f_\ell \bot \CZ_s(T_\ell)$, then there exists a sequence $N_k\to\infty$ (depending on the previous data) such that for $\mu$-a.e. $x\in X$
	we have 
	$$
	\lim_{k\to\infty} \E_{n\in [N_k]} \, \CN_x(n) \cdot \prod_{i=1}^\ell f_i(T_i^{p_i(n)}x) =0
	$$ 
	for every $r$-step nilsequence $(\CN_x(n))_{n\in \N}$.  Furthermore, $(N_k)_{k\in\N}$ can be taken to be a subsequence of any other given sequence $N_k'\to\infty$. 
\end{lemma}
\begin{proof}
By \cref{C:1}, it is enough to prove the claim with an arbitrary
collection of level $r+1$ dual sequences $(\CD_x(n))$ in place of the $r$-step
nilsequences~$(\CN_x(n))_{n\in\N}$.
	
Let $N_k\to \infty$ be a sequence of integers, to be chosen shortly.

Using the argument in the proof of
\cite[Proposition 6.1]{frantzikinakis2015}, we get that for every
$x\in X$ there exist a constant $C_x=C_x(d,\ell,r)$ and $w_0=w_0(d,\ell,r)$ such
that for some $w\in [w_0]$ the expression
$$
\limsup_{k\to\infty}\Big|\E_{n\in [N_k]} \, \CD_x(n) \cdot \prod_{i=1}^\ell f_i(T_i^{p_i(n)}x)\Big|^{2^w}  
$$
is bounded by 
$$
C_x\cdot \limsup_{H_w\to \infty} \E_{h_w\in [H_w]}\cdots \limsup_{H_1\to \infty} \E_{h_1\in [H_1]} \limsup_{k\to\infty} \Big| \E_{n\in [N_k]} \prod_{\epsilon \in \{0,1\}^w} \CC^{|\epsilon|}A_x(n+\epsilon \cdot h)\Big| 
$$
where
$$
A_x(n):=\prod_{i=1}^\ell f_i(T_i^{p_i(n)}x).  
$$
	

{\bf Case 1.} Suppose that $\ell=1$ and $p_1(n)=an+b$ for some $a,b\in \Z$ with $a\neq 0$.
We claim that the iterated limit vanishes for $\mu$-a.e. $x\in X$
if~$f_1\bot \CZ_{w_0}(T_1)$.
Since we want to show convergence to~$0$ for $\mu$-a.e. $x\in X$, it is enough
to do so for each $T_1^a$-ergodic component,
and since Host-Kra factors commute with ergodic decompositions
(\cite[Chapter 8, Proposition 18]{host2018nilpotent}), without loss of generality
we may assume~$T_1^a$ ergodic and replace~$f_1$ by~$T^b_1 f_1$.
By the ergodic theorem for $T_1^a$, for $\mu$-a.e.\ $x\in X$ and every $h\in\N^w$
the innermost limit
equals~$\abs{\int\prod_{\epsilon\in\{0,1\}^w}
	\CC^{\abs{\epsilon}}T_1^{a\epsilon\cdot h}f_1\,d\mu}$.
By the inductive definition of seminorms and Cauchy-Schwarz, the iterated limit
is bounded by~$\snorm{f_1}^{2^{w}}_{w+1,T_1^a}$.
By \cite[Theorem A.2]{bergelson2006combinatorial}, which
proves~$\CZ_s(T^a)=\CZ_s(T)$ for~$s\geq 1$ and~$a\neq 0$,
it follows that we can take~$s=w_0$.
	
	{\bf Case 2.} 
Suppose that either $\ell=1$ and $\deg(p_1)\geq 2$, or $\ell\geq 2$,
$p_1,\ldots, p_\ell$ are pairwise independent, and $\deg(p_\ell)\geq 2$.
In the expression $\CC^{|\epsilon|}A_x(n+\epsilon \cdot h)$, we can incorporate
all linear iterates into a single term. After doing this, it is easy to verify
that outside a negligible set of $h\in \N^w$ for our averaging, the polynomial
iterates that appear in the expression $\CC^{|\epsilon|}A_x(n+\epsilon \cdot h)$
are pairwise independent. Indeed, one only has to note that if the polynomials
$p(n)$, $q(n)$ are independent and $\deg(p)\geq 2$, then so are the polynomials
$p(n+h)$, $q(n)$ for all but at most one value of $h\in \N$. It follows from
\cite[Theorem~2.8]{frantzikinakis2025joint} that if $s$ is large enough, depending only on
$d,\ell,r$, and $f_\ell \bot \CZ_s(T_\ell)$, then
	\begin{equation}\label{E:fgh}
	\E_{n\in [N]} \prod_{\epsilon \in \{0,1\}^w} \CC^{|\epsilon|}A_x(n+\epsilon \cdot h)
	\end{equation}
converges to $0$ in $L^2(\mu)$. Hence, using a diagonal argument, we get that there exists a subsequence $N_k\to\infty$
along which the convergence holds $\mu$-a.e. for every $w\in [w_0]$  all those good values of
$h\in \N^w$. This is the subsequence we seek.

The same argument shows that 	$(N_k)_{k\in \N}$ can be taken to be a subsequence of any other given sequence $N_k'\to\infty$. 
\end{proof}

\cref{cor:partial_factor_use}
follows by taking~$\CN_x(n)=\CN(x,q_1(n),\ldots,q_m(n))$
in the following result.


\begin{proposition}
		Let  $d,\ell, r\in \N$ and  $p_1,\ldots, p_\ell$ 
		 be pairwise independent polynomials of degree at most $d$.
Then there exists $s\in \N$, dependent only on $d, \ell,r$, such that the following holds: If $(X,\mu,T_1,\ldots, T_\ell)$ is a system,  and $f_0,f_1,\ldots, f_\ell \in L^\infty(\mu)$ are $1$-bounded and   $f_i \bot \CZ_s(T_i)$ for some $i\in [\ell]$, then 
	\begin{equation}\label{E:int1}
		\lim_{N\to\infty} \E_{n\in [N]} \, \int \CN_x(n) \cdot f_0(x) \cdot \prod_{i=1}^\ell f_i(T_i^{p_i(n)}x)  \d\mu =0
	\end{equation} 
	for every measurable varying collection of  $1$-bounded $r$-step nilsequences $(\CN_x(n))_{n\in \N}$. 
\end{proposition}

\begin{proof}
If $\ell=1$ or $\ell\geq 2$ and all the polynomials have degree at least $2$,
then the conclusion follows easily from \cref{L:2}. Indeed, if the limit in
\eqref{E:int1} is non-zero, then it is bounded away from zero on a subsequence
$N_k\to\infty$, and we get a contradiction from \cref{L:2} and the bounded
convergence theorem.

We now deal with the complementary case. Since $p_1,\ldots, p_\ell$ are pairwise
independent, it remains to treat the case where $\deg(p_1)=1$ and
$\deg(p_i)\geq 2$ for $i=2,\ldots, \ell$.
	Using \cref{L:2} again, we deduce that in order to establish \eqref{E:int1},
	we may assume that $f_i$ is $\CZ_s(T_i)$-measurable for $i=2,\ldots, \ell$,
	for some sufficiently large $s\in \N$, depending only on $d,\ell,r$. Using this and \cite[Proposition 3.1]{chu2011ergodic},
	we get that, for $i=2,\ldots, \ell$ and every $\varepsilon>0$, we have a decomposition
	 $$
	 f_i=f_{i,1}+f_{i,2}
	 $$
	where $\norm{f_{i,1}}_{L^1(\mu)}\leq \varepsilon$ and
	$(f_{i,2}(T_i^n x))_{n\in\N}$ is an $s$-step nilsequence for $\mu$-a.e. $x\in X$.
	Note also that $(f_{i,2}(T_i^{p_i(n)}x))_{n\in\N}$ is then a
	$(ds)$-step nilsequence (\cref{sec:prel_nil_definitions}).
It follows that, in proving \eqref{E:int1}, we can assume that, for $\mu$-a.e. $x\in X$
and for $i=2,\ldots, \ell$, the sequences $(f_i(T_i^{p_i(n)}x))_{n\in\N}$ are
$(ds)$-step nilsequences, varying measurably with respect to $x\in X$. Hence,
since the product of a $d_1$-step nilsequence and a $d_2$-step nilsequence is a
$(d_1+d_2)$-step nilsequence, we will be done if we show that there exists
$s'\in \N$, depending only on $d,s,\ell,r$, such that if $f_1\bot \CZ_{s'}(T_1)$ and $a,b\in \Z$ with $a\neq 0$, we have
	$$
	\lim_{N\to\infty} \E_{n\in [N]} \, \int \CN'_x(n) \cdot f_0(x) \cdot f_1(T_1^{an+b}x) \d\mu=0
	$$ 
for every measurable varying collection of $1$-bounded $(d s(\ell-1)+r)$-step nilsequences
$(\CN'_x(n))_{n\in \N}$. This follows again from the $\ell=1$ case of \cref{L:2} and the bounded convergence theorem.
\end{proof}

	


\section{Proof of lifting property} \label{sec:proof_lifting}
In this section we combine the tools developed in the previous sections in order to prove a key lifting result. This result allows us to pass multiple recurrence properties from a factor to an extension, under suitable compactness assumptions. We first define the recurrence property that will be lifted.

\begin{definition}[PSZ Property]\label{D:PSZ}
Let  $d,\ell,r\in \N$. We say that a commuting system
	$(X,\mu,T_1,\ldots,T_\ell)$ has the \emph{$PSZ_r$ property}
	if the following property holds: For all
	$m\in \N$, all jointly intersective families of integral polynomials 
	$ p_1,\ldots,p_\ell,q_1,\ldots,q_m$ with $p_1,\ldots, p_\ell$ pairwise independent and of degree at most $d$, all 
	measurable varying positive nilsequences
	$\CN\colon X\times\Z^m\to[0,1]$
    such that~$\CN(x,q_1(n),\ldots,q_m(n))$ is $r$-step,
    and all sets~$A\subseteq  X$ of positive
	measure, one has
	$$
	\liminf_{N\to\infty} \frac1N\sum_{n=1}^N
	\int_X \CN(x,q_1(n),\ldots,q_m(n))\cdot \1_A(x) 	\cdot 	\prod_{i=1}^{\ell}
	T_i^{p_i(n)}\1_A(x)	\,d\mu(x)>0.
	$$
	We say that the system has the \emph{$PSZ$ property} if it has the
	$PSZ_r$ property for all~$r\in \N$.
\end{definition}

\begin{remark}
	Throughout this section the parameters $d$ and $\ell$ are fixed. Thus the 
 properties $PSZ_r$ and~$PSZ$ are understood with respect to these parameters.
\end{remark}

\begin{remark}
We recall that if~$\CN\colon X\times\Z^m\to[0,1]$ is an $s$-step
measurable varying nilsequence (\cref{def:nil_meas_step}),
that is~$(\CN(x,n_1,\ldots,n_m))_{n\in \N^m}$ is $s$-step
for almost every~$x\in X$,
then~$(\CN(x,q_1(n),\ldots,q_m(n)))_{n\in \N}$
is $r$-step for almost every~$x\in X$,
where~$r:= 2s\cdot \max_{1\leq i\leq m} \deg q_i$,
see remark after \cref{def:nil_meas_step}.
Thus,
instead of bounding~$r$ in \cref{D:PSZ},
one can alternatively bound the step~$s$ of~$\CN$
together with the maximum degree among the~$q_i$.
\end{remark}

\begin{proposition}\label{P:compactlift}
	Let $d, \ell,r \in \N$ and let  $s=s(d,\ell,r)$ be as in \cref{cor:partial_factor_use}.
	Let $(X,\mu,T_1,\ldots, T_\ell)$ be a system 	and~$Y,W$ be  factors with $Y\subseteq W$. Suppose that the factor $Y$  has the $PSZ$ property
	and that the extension 
	$$
	Y\vee (\CZ_s(X,T_i)\wedge W)\to Y
	$$
	is $T_i$-compact  for $i=1,\ldots, \ell$. Then the  factor  $W$ has the $PSZ_r$ property.
\end{proposition} 

\begin{proof}
	Take any~$m\in \N$,
	any jointly intersective  polynomial family 
	$$
	p_1,\ldots,p_\ell,q_1,\ldots,q_m
	$$
such that the polynomials $p_1,\ldots, p_\ell$ are pairwise independent and of degree at most $d$,
any measurable varying positive nilsequence~$\CN\colon W\times \Z^m \to [0,1]$ such that
\[
    \CN(x,q_1(n),\ldots,q_m(n))
\]
is at most $r$-step,
and any $W$-measurable set~$A$ of positive measure.
	We want to prove that
	\[
	\liminf_{N\to\infty} \E_{n\in[N]}
	\int_X \CN(x,q_1(n),\ldots,q_m(n))
	\cdot \1_{A}(x) \cdot \prod_{i=1}^\ell T_i^{p_i(n)} \1_{A}(x)
	\d \mu(x) > 0.
	\]
	For convenience, we divide the argument into several steps.
	
	\vspace{.5em}
	\textit{1. Partial characteristic factors.} 
	By \cref{cor:partial_factor_use} applied to the system~$X$ and our choice of~$s=s(d,\ell,r)$,
	the limiting value will not change if we  replace the integral in the previous
	average by
	\[
	I_n:=\int_X \CN(x,q_1(n),\ldots,q_m(n)) \cdot \1_{A}(x)
	\cdot \prod_{i=1}^\ell T_i^{p_i(n)} f_i(x)
	\d \mu(x),
	\]
	where~$f_i:=\E(\1_A\mid \CZ_s(X,T_i))$ for~$i=1,\ldots,\ell$.
	For notational convenience, we define~$f_0:=\1_A$.
	Let $\mu=\int \mu_y\d\nu(y)$ be the disintegration of the measure $\mu$ with respect to the factor $Y$ and for $y\in Y$ and $n \in \N$ let 
	 $$
	I_n(y):=\int_X \CN(x,q_1(n),\ldots,q_m(n)) \cdot \1_{A}(x)
	\cdot \prod_{i=1}^\ell T_i^{p_i(n)} f_i(x)
	\d \mu_y(x).
	$$
	Then 
\begin{equation}\label{E:In}
	I_n=\int I_n(y)\d\nu(y).
	\end{equation}
	Our plan is  to obtain lower bounds for the expressions 
	$I_n(y)$ using \cref{P:nil+nul=useful}, and \cref{cor:beautiful_main}  
	and then insert these bounds into \eqref{E:In}.
	This will reduce matters to proving   a recurrence property for the factor $Y$, which
	will follow from our assumption that the factor $Y$ has the $PSZ$ property.

	\vspace{.5em}
	\textit{2. Lower bound for~$f_i$.} 
	By Chu's Inequality~\cite[Lemma 1.6]{chu2011multiple},
	applied to~$f_0:=\1_A$ and~$f_i:=\E(\1_A\mid \CZ_s(X,T_i))$,
	we have~$\int_X \prod_{i=0}^\ell f_i \d\mu \geq \mu(A)^{\ell+1}$.
	Let~$c:=\mu(A)^{\ell+1}$.
	Let~$B\subseteq  Y$ be the set of~$y\in Y$ such that
	\[
	\int_X
	\prod_{0\leq i\leq \ell}  f_i
	\d\mu_y \geq \frac{1}{2}c.
	\]
	By Markov's Inequality we have~$\nu(B)\geq \frac{1}{2}c$.
	For every~$y\in B$, let~$B_y$ be the set of~$x\in X$
	such that~$\prod_{0\leq i\leq \ell} f_i(x) \geq \frac{1}{4}c$.
	Note that~$B_y$ depends measurably on~$y \in B$.
	By Markov's Inequality, we have~$\mu_y(B_y)\geq \frac{1}{4}c$
	for every~$y\in B$.
	Now, up to replacing~$c$ by~$\frac{1}{4}c$ and weakening the previous inequalities,
	we have the following: there is a measurable set~$B\subseteq  Y$
	and measurable varying sets~$B_y\subseteq  X$ such that
	\begin{align}
		& \nu(B)\geq c, \nonumber \\
		& \mu_y(B_y)\geq c, && \text{for all } y\in B, \nonumber \\
		& \prod_{0\leq i\leq \ell} f_i(x) \geq c,
		&& \text{ for all } y\in B \text{ and  all } x\in B_y. \label{eq:Bc}
	\end{align}
	Note that these estimates  remain valid if we further shrink~$B$ and~$c$
	accordingly (which we will do).
	
	\vspace{.5em}
	\textit{3. Lower bound for~$\CN$.} 
Let $G/\Gamma$ be the nilmanifold  and $f\colon X\times G/\Gamma\to [0,1]$ be the Borel function associated to $\CN$ as in  \cref{def:mv_nilsequence}. We apply  \cref{P:nil+nul=useful} (see also the second remark after its statement) to this nilmanifold, this function, and the factor $Y$. We deduce that  
up to shrinking~$B$ to a subset of
	$\nu$-positive measure contained in a single piece of the countable
	partition in \cref{P:nil+nul=useful}, we get that there exist integers $\alpha\in \N$ and $\beta \geq 0$ and integral polynomials $r_1,\ldots, r_k$ such that if we set
	$$
	\tilde{p}_i(n):=p_i(\alpha n+\beta), \qquad \tilde{q}_i(n):=q_i(\alpha n +\beta), 
	$$
	 we have 
	\[
	\int_{B_y} \CN(x, \tilde{q}_1(n),\ldots,\tilde{q}_m(n)) \d \mu_y(x)
	\geq \CN'(y,r_1(n),\ldots,r_k(n))+e(y,n)
	\]
	for every~$y\in B$, where~$\CN'$ is a $Y$-measurable varying positive nilsequence, $e$ is a measurable varying nullsequence, and
	the polynomials 
	$$\tilde{p}_1,\ldots,\tilde{p}_\ell, \tilde{q}_1,\ldots,\tilde{q}_m,r_1,\ldots, r_k
	$$ are jointly intersective.
	(To be precise, the sets~$B_y$, and thus~$\CN'(y,\cdot)$ and $e(y,\cdot)$, are only defined for~$y\in B$.
	For~$y\in Y\setminus B$, we simply let~$\CN'(y,\cdot) := 1$ and $e(y,\cdot)=0$.)
	
	Since the function~$y\mapsto \CN'(y,0,\ldots,0)$ is measurable and positive,
	up to further shrinking~$B$ (and~$c$ accordingly) we can suppose~$\CN'(y,0,\ldots,0)\geq c$
	for all~$y\in B$ (and thus for all~$y\in Y$).
	By \cref{lem:nilsequence=nilbohr_measurable}, 
	there is a measurable varying $\nilbor$ sequence~$\CN''$
	such that~$\CN''(y,n_1,\ldots,n_k)=1$ implies~$\CN'(y,n_1,\ldots,n_k)\geq c/2$.
	Up to replacing~$c$ by~$c/2$, we have
	\begin{equation} \label{eq:N''}
		\CN''(y,r_1(n),\ldots,r_k(n))=1
		\quad \implies \quad \int_{B_y} \CN(x,\tilde{q}_1(n),\ldots,\tilde{q}_m(n)) \d \mu_y \geq c+e(y,n).
	\end{equation}
	
	\textit{4. Approximating~$f_i$ by $\nilbor$ return sets.}
	We now choose~$\ve_1$ and~$\ve_2$ as functions of the constant~$c$
	(which, therefore, can no longer be modified). For reasons that will be clear
	later, we choose~$\ve_1,\ve_2>0$ such that
	\begin{equation} \label{eq:ve12}
		c^2	- 2\ell \ve_2- \ell \ve_1>0, \quad c-\ell \ve_1^2\ve_2^{-2}>0.
	\end{equation}
	To see that this system of inequalities has positive solutions, first
    choose~$\ve_2$ such that~$2\ell \ve_2 < c^2/2$,
    then choose~$\ve_1\in (0,1)$
	small enough so that both inequalities hold.  
	Fix~$i\in \{1,\ldots,\ell\}$. By definition~$f_i=\E(\1_A\mid \CZ_s(X,T_i))$.
	Since~$A$ is $W$-measurable, by the good
	intersection property (\cref{L:ZsGood}) $f_i$ is $(\CZ_s(X,T_i)\wedge W)$-measurable.
	It follows that
		\begin{align*}
		\CZ_s(Y\vee (\CZ_s(X,T_i)\wedge W),T_i)
		& \supseteq \CZ_s(\CZ_s(X,T_i)\wedge W,T_i) = \CZ_s(\CZ_s(W,T_i),T_i) \\
		& = \CZ_s(W,T_i) = \CZ_s(X,T_i)\wedge W  \ni f_i.
	\end{align*}
	In particular \cref{cor:beautiful_main} applies to the compact extension
	$
	Y\vee (\CZ_s(X,T_i)\wedge W) \to Y
	$
	and the function~$f_i$, giving a function
	$
	g_i\in L^\infty(\CZ_s(Y\vee (\CZ_s(X,T_i)\wedge W),T_i)),
	$
	taking values in  $[0,1]$
	with~$\norm{f_i-g_i}_{L^2(\mu)}\leq \ve_1$, and  such that the indicator of the set 
	\[
	\{(y,n)\in Y\times \Z \colon \norm{T_i^{n} g_i - g_i}_{L^2(\mu_y)} \leq \ve_1\}
	\]
	is bounded from below by a measurable varying $\nilbor$
	sequence~$\widetilde \CN_i''' \colon Y\times\Z\to \{0,1\}$.\footnote{This fact is crucial, and it is what allows us to bypass the need to prove  a still unavailable coloring theorem, as in the Furstenberg-Katznelson~\cite{furstenberg1978ergodic} and Bergelson-Leibman~\cite{bergelson1996polynomial} approach. It enables us, in the next step of the argument, to replace the iterates $T_i^{p_i(n)}g_i$ with $g_i$ when $n$ is restricted to a $\nilbor$ subset of the integers. To handle the additional difficulty caused by the introduction of these $\nilbor$ subsets, we use the nil-plus-null decomposition result of \cref{P:nil+nul=useful}.}
	In particular
	\[
	\CN'''(y,n_1,\ldots,n_l) := \prod_{i=1}^\ell \widetilde{\CN}'''_i(y,n_i)
	\]
	is a measurable varying $\nilbor$ sequence and for all $y\in Y$ we have 
	\begin{equation} \label{eq:N'''}
		\CN'''(y,\tilde{p}_1(n),\ldots,\tilde{p}_\ell(n))=1 \implies
		\norm{T_i^{\tilde{p}_i(n)} g_i - g_i}_{L^2(\mu_y)} \leq \ve_1, \qquad \text{for } i=1,\ldots,\ell.
	\end{equation}
	Moreover,
	let~$Y_i\subseteq  Y$ be the set of~$y$ where~$\norm{f_i-g_i}_{L^2(\mu_y)} \leq \ve_2$.
	Define~$Y':=\bigcap_{i=1}^\ell Y_i$. Then
	\[
	\nu(Y') \geq 1 - \ell \ve_1^2\ve_2^{-2}
	\]
	and
	\begin{equation} \label{eq:Yi}
		\1_{Y'}(y)=1 \quad \implies \quad
		\norm{f_i-g_i}_{L^2(\mu_y)}\leq \ve_2, \qquad \text{for } i=1,\ldots,\ell.
	\end{equation}
	For notational convenience, we define~$g_0:=f_0=\1_A$.
	
	\vspace{.5em}
	\textit{5. Lower bounds for $I_n(y)$.}  
	For all pairs~$(y,n)\in Y\times \Z$
	satisfying
	\begin{multline} \label{eq:indicator}
	R_n(y):=	\CN''(y,r_1(n),\ldots,r_k(n)) \cdot \CN'''(y,\tilde{p}_1(n),\ldots,\tilde{p}_\ell(n))\\
		\cdot \1_B(y)
		\cdot \1_{Y'}(y) \cdot \prod_{1\leq i \leq \ell} \1_{Y'}(T_i^{\tilde{p}_i(n)}y) > 0
	\end{multline}
	we have, by the previous estimates and Cauchy-Schwarz (CS),
	the following chain of inequalities:
	{\allowdisplaybreaks
		\begin{align*}
			&
			I_{\alpha n+\beta}(y)=\int_X
			\CN(\cdot,\tilde{q}_1(n),\ldots,\tilde{q}_m(n)) \cdot \prod_{0\leq i\leq \ell} T_i^{\tilde{p}_i(n)} f_i
			\d\mu_y
			\\
			& \geq
			\int_X
			\CN(\cdot,\tilde{q}_1(n),\ldots,\tilde{q}_m(n)) \cdot \prod_{0\leq i\leq \ell} T_i^{\tilde{p}_i(n)} g_i
			\d\mu_y
			-\ell \ve_2
			&& \text{by~\eqref{eq:Yi} + CS}\\
			& \geq
			\int_X
			\CN(\cdot,\tilde{q}_1(n),\ldots,\tilde{q}_m(n)) \cdot \prod_{0\leq i\leq \ell} g_i
			\d\mu_y
			-\ell\ve_2 - \ell\ve_1
			&& \text{by~\eqref{eq:N'''} + CS}\\
			& \geq
			\int_X
			\CN(\cdot,\tilde{q}_1(n),\ldots,\tilde{q}_m(n)) \cdot \prod_{0\leq i\leq \ell} f_i
			\d\mu_y
			- 2\ell\ve_2 - \ell\ve_1
			&& \text{by~\eqref{eq:Yi} + CS}\\
			& \geq
			c \cdot
			\int_{B_y}
			\CN(\cdot,\tilde{q}_1(n),\ldots,\tilde{q}_m(n))
			\d\mu_y
			- 2\ell\ve_2 - \ell\ve_1
			&& \text{by~\eqref{eq:Bc}}\\
			& \geq
			c^2	- 2\ell\ve_2 - \ell\ve_1+c\,  e(y,n)
			&& \text{by~\eqref{eq:N''}}.
		\end{align*}
	}
	By the choice of~$\ve_1$ and~$\ve_2$ (equation~\eqref{eq:ve12}),
		$a:=c^2	- 2\ell\ve_2 - \ell\ve_1>0$  is strictly positive. Hence, 
	\begin{equation}\label{E:aIower}
	I_{\alpha n+\beta}(y)\geq a R_n(y)+c\,  e(y,n) R_n(y), \qquad \text{ for all } y\in Y, 
\end{equation}
 where $R_n(y)$ is as in \eqref{eq:indicator}. 
	
	\vspace{.5em}
	\textit{6. Concluding the proof.}  	
	Recall our goal is to show that 
	$$
	\liminf_{N\to\infty}\E_{n\in [N]}I_n>0, 
	$$
	and since $I_n\geq 0$, it suffices to show that 
	$$
	\liminf_{N\to\infty}\E_{n\in [N]}I_{\alpha n+\beta}>0. 
	$$
	Using \eqref{E:In} it suffices to show that 
	$$
	\liminf_{N\to\infty}\E_{n\in [N]}\int I_{\alpha n+\beta}(y)\d\nu(y)>0. 
	$$
	Since \eqref{E:aIower}  holds  and  $0\leq R_n\leq 1$, 	 we get that 
	 $$
	 	\int I_{\alpha n+\beta}(y) \d \nu(y) \geq a \int_Y R_n(y)\d\nu(y)-c \int_Y |e(y,n)|\d\nu(y)
	 	\qquad \text{ for all } n\in \N.
	 $$
Moreover, since  $e(y,n)$ is  a measurable varying nullsequence with values bounded by $1$, the bounded convergence theorem gives 
	$$
	\lim_{N\to\infty} \E_{n\in[N]} \int_Y|e(y,n)|\d\nu(y)=0.
	$$
Combining the above, we get that it is  enough to show  that 
$$
\liminf_{N\to\infty} \E_{n\in [N]} \int_Y R_n \d\nu>0.
$$
	Up to replacing both~$B$ and~$Y'$
	by their intersection~$B\cap Y'$ (which has positive measure
	since~$c-\ell \ve_1^2\ve_2^{-2}>0$),
	and defining~$\CM :=\CN''\cdot \CN'''$, we are left with establishing 
that 
	\begin{multline}\label{E:PSZY}
		\liminf_{N\to\infty} \E_{n\in[N]} \int_Y
\CM(y,r_1(n),\ldots,r_k(n),\tilde{p}_1(n),\ldots,\tilde{p}_\ell(n))
	\cdot
	\1_{B}(y)
	\cdot\\ 
	\prod_{1\leq i \leq \ell}
	\1_{B}(T_i^{\tilde{p}_i(n)}y)
	\d\nu(y)>0.
	\end{multline}
	By \cref{lem:nilsequence=nilbohr_measurable}, we have 
	$\CM(y,n_1,\ldots,n_{k+\ell}) \geq \widetilde{\CN}(y,n_1,\ldots,n_{k+\ell})$
	for $\nu$-almost every $y\in Y$
	for some measurable varying positive nilsequence~$\widetilde{\CN}$. Note that the  polynomials
    $\tilde{p}_1,\ldots,\tilde{p}_\ell$ are pairwise independent with   degree at most $d$,  and the family 
	$$
r_1,\ldots,r_k,\tilde{p}_1,\ldots,\tilde{p}_\ell,\tilde{p}_1,\ldots,\tilde{p}_\ell
	$$
is jointly intersective. Using this and that 
	 the factor $(Y,\nu)$ has the  PSZ property, we deduce that \eqref{E:PSZY} holds,  concluding the proof.
\end{proof}
 
   We will also use  the next corollary.
  \begin{corollary}\label{C:compactlift}
 	Let  $(X,\mu,T_1,\ldots, T_\ell)$ be a system 	and $Y,W$ be  factors with $Y\subseteq W$. Suppose that $Y$  has the $PSZ$ property 
 	and that the extension 
 	$
 	 W\to Y
 	$
 	is $T_i$-compact  for all $i\in [\ell]$. Then the  factor  $W$ has the $PSZ$ property. 
 \end{corollary} 
 \begin{proof}
  It suffices to show that $W$ has the $PSZ_r$ property for every $r\in \N$. So fix $r\in \N$ and  let $s\in \N$ be arbitrary. 
Since  	 $W\to Y$ is $T_i$-compact
and $Y\vee (\CZ_s(X,T_i)\wedge W)$ is an intermediate extension,  
using   property~\eqref{i}  of  \cref{L:CompactProperties}  we get 
 	$Y\vee (\CZ_s(X,T_i)\wedge W)\to Y$ is $T_i$-compact for all $i\in [\ell]$. 
Since $Y$ has the $PSZ$ property, for  $s=s(d,\ell,r)\in \N$ we get that  \cref{P:compactlift} applies and gives that $W$ has the $PSZ_r$ property.  	 
\end{proof}

\section{Proof of the PSZ property} \label{sec:final_proof}

In this section, we describe
the induction scheme that,
via the lifting property of the previous section,
completes the proof of the main theorem.
We will use the simple observation
that if~$T_1,\ldots,T_\ell$ are commuting
transformations, then for every~$1\leq i\leq \ell$
and~$s\geq 1$ the Host-Kra factor~$\CZ_s(T_i,X)$ is a
common factor of all transformations.
More generally, if~$W$ is a common factor,
then~$\CZ_s(T_i,W)=\CZ_s(T_i,X)\wedge W$ is a common factor.
To shorten notation,
we will denote~$Z_s^i(W):=\CZ_s(T_i,W)$.

\subsection{A motivating example}
Our initial idea was to combine the lifting property with the following
finite-tower construction for $\CZ_s$-systems.
\begin{proposition} \label{P:MotivationTower}
Let $s\geq 0$ be an integer and let $(X,\mu,T_1,\ldots,T_\ell)$
be a commuting system such that $X=\CZ_s(T_i)$ for every
$i\in[\ell]$. Then there exists   $N\leq (s+1)^\ell$ and a tower of common factors
$$
\{\text{pt}\}=Y_0\subseteq  Y_1\subseteq \cdots\subseteq  Y_N=X,
$$
such that $Y_j\to Y_{j-1}$ is $T_i$-compact for every $j\in[N]$
and every $i\in[\ell]$.
\end{proposition}
\begin{remark}
For $\ell=s=2$, the construction above gives the following explicit tower of $T_i$-compact extensions for $i=1,2$.
Write $A_r:=\CZ_r(T_1)$ and $B_r:=\CZ_r(T_2)$ for $0\leq r\leq 2$,
then  $A_2=B_2=X$. 
\begin{align*}
\{\text{pt}\}
&\subseteq  A_0\wedge B_0
 \subseteq  A_0\wedge B_1
 \subseteq  A_0\\
&\subseteq  A_0\vee\big(A_1\wedge B_0\big)
 \subseteq  A_0\vee\big(A_1\wedge B_1\big)
 \subseteq  A_1\\
&\subseteq  A_1\vee B_0
 \subseteq  A_1\vee B_1
 \subseteq  X.
\end{align*}
\end{remark}
\begin{proof}
We induct on $\ell$, using only joins and intersections of Host-Kra
factors $\CZ_s(T_i)$ and the point factor; these are common factors by commutativity. We use the convention $\CZ_{-1}(T_i)=\{\text{pt}\}$.

For $\ell=1$, take
$$
\{\text{pt}\}\subseteq \CZ_0(T_1)
\subseteq \cdots\subseteq \CZ_s(T_1)=X.
$$
Its $s+1$ steps are compact by part~\eqref{iv} of \cref{L:CompactProperties}.
Suppose $\ell\geq 2$ and suppose we have  a tower
$Y_0\subseteq \cdots\subseteq  Y_M=X$ with the required compactness property for $T_1,\ldots,T_{\ell-1}$,
with $M\leq(s+1)^{\ell-1}$.
For each $j\in[M]$, replace $Y_{j-1}\subseteq  Y_j$ by the $s+1$
step tower 
$$
V_{j,r}:=Y_{j-1}\vee\big(Y_j\wedge\CZ_r(T_\ell)\big),
\qquad -1 \leq r\leq s.
$$
We then have $V_{j,-1}=Y_{j-1}$ and $V_{j,s}=Y_j$.

For $i<\ell$, each extension $V_{j,r}\to V_{j,r-1}$ is
$T_i$-compact by parts~\eqref{i} and~\eqref{ii} of
Lemma~\ref{L:CompactProperties}, since it is an intermediate
step of $Y_j\to Y_{j-1}$.
For $i=\ell$,  part~\eqref{iv} of \cref{L:CompactProperties}
gives compactness of
$ Y_j\wedge\CZ_r(T_\ell) \to  Y_j\wedge\CZ_{r-1}(T_\ell)$.
Joining both factors with $Y_{j-1}$ and applying
part~\eqref{iii} of \cref{L:CompactProperties} shows that
$V_{j,r}\to V_{j,r-1}$ is $T_\ell$-compact.
Thus the refined tower has  the required properties and at most $(s+1)M\leq(s+1)^\ell$ steps,
completing the induction.
\end{proof}
Since the point system has the $PSZ$ property by
\cite[Proposition~2.4]{bergelson2008intersective},
\cref{P:MotivationTower} and \cref{C:compactlift} establish this property
in the structured setting after at most $(s+1)^\ell$ lifting steps.
Unfortunately, the characteristic factor reduction in \cref{cor:partial_factor_use} does not place us in this setting;  it projects the individual
functions onto their respective Host-Kra factors $\CZ_s(T_i)$, but does not replace
the ambient system by one satisfying $X=\CZ_s(T_i)$ for every
$i\in[\ell]$. To overcome this we plan to retain the finite construction and  adapt the induction statement, 
using the finite
tower supplied by \cref{L:OneCompactLadder} below to  induct on the number of directions in which compactness
is not assumed, retaining the full $PSZ$ property as the conclusion. The details are given in the next two subsections. 

\subsection{Auxiliary  lemmas} The following lemma collects properties used below.
Some of these facts are standard and we recall their proofs for convenience.
 \begin{lemma}\label{L:CompactProperties}
	Let $(X,\mu,T)$ be a system and $Y,Z,W$ be factors with $Z\subseteq  Y$. 
	\begin{enumerate}[(i)]
		\item \label{i}  If $Y\to Z$ is $T$-compact and $Z\subseteq  Y'\subseteq   Y$, then  $Y'\to Z$ is $T$-compact.
		\item  \label{ii} If $Y\to Z$ is $T$-compact and $Z\subseteq   Z'\subseteq  Y$, then  $Y\to Z'$ is $T$-compact.
		
		\item \label{iii} If $Y\to Z$ is $T$-compact, then  $Y\vee W \to Z\vee W$ is $T$-compact.
		
		\item  \label{iv}  
		 $\CZ_s(T)\wedge Y\to \CZ_{s-1}(T)\wedge Y$ is $T$-compact
		for every  $s\geq 0$,  where $\CZ_{-1}(T):=\{pt\}$.
		
	\end{enumerate}
    
\end{lemma}
\begin{proof}
Parts~$(i)$ and~$(ii)$ are \cite[Lemma 9.12]{glasner2003ergodic}.

Let us prove~$(iii)$.
Let~$f\in L^\infty(Y)$ and~$g\in L^\infty(W)$ with~$\norm{g}_\infty\leq 1$.
Since the extension~$Y\to Z$ is compact, for any~$\ve>0$ there is a
finite-rank $T$-invariant $L^\infty(Z)$-module,
say~$M:=\sum_{i=1}^r L^\infty(Z)f_i$, approximating~$f$:
there are~$h_1,\ldots,h_r\in L^\infty(Z)$
such that~$\norm{f-\sum_{i=1}^r h_i f_i}_2 \leq \ve$.
Since~$\norm{g}_\infty\leq 1$,
we have~$\norm{gf-\sum_{i=1}^r g h_i f_i}_2 \leq \ve$.
Note that~$\sum_{i=1}^r g h_i f_i$ belongs
to~$\sum_{i=1}^r L^\infty(Z\vee W)f_i =: \tilde M$.
Clearly~$\tilde M$ is a finite-rank $L^\infty(Z\vee W)$-module.
From~$T f_i \in M\subseteq  \tilde M$
and~$T L^\infty(Z\vee W) \subseteq  L^\infty(Z\vee W)$
it follows that~$\tilde M$ is $T$-invariant.
Since functions of the form~$gf$ above have dense linear span in~$L^2(Y\vee W)$,
we conclude that finite-rank $T$-invariant $L^\infty(Z\vee W)$-modules
have dense span in~$L^2(Y\vee W)$, proving the claim.

Let us prove~$(iv)$. Since~$\CZ_k(X,T)\wedge Y = \CZ_k(Y,T)$,
it is enough to prove that the extension~$\CZ_s(T)\to \CZ_{s-1}(T)$
is $T$-compact. The case~$s=0$ is trivial, as~$\CZ_0(T)$ is the invariant
factor, so let~$s\geq 1$.
Thus, let~$(X,\mu,T)$
be a $\CZ_s$-system and~$f\in L^\infty(X)$ with~$f\neq 0$.
Since~$f=\E(f\mid \CZ_s)\neq 0$, in particular
\[
    0 < \snorm{f}_{s+1}^{2^{s+1}} = \limsup_{N\to\infty} \frac{1}{N} \sum_{n=1}^N
                            \snorm{T^n \con f \cdot f}_{s}^{2^s},
\]
see \cref{sec:prel_HK}. Let~$g:=T^n \con f \cdot f$
and~$h:=\E(g\mid \CZ_{s-1})$.
By the triangle inequality
we have~$\snorm{g}_s \leq \snorm{g-h}_{s} + \snorm{h}_s$.
Since~$\E(g-h\mid \CZ_{s-1})=0$, we
have~$\snorm{g-h}_{s}=0$.
On the other hand, by \cite[Chapter 8, Proposition 15]{host2018nilpotent}
we have~$\snorm{h}_s\leq \norm{h}_{2^s}$,
thus
\[
    \snorm{h}_s^{2^s} \leq \int \abs{h}^{2^s} \d\mu
        \leq \norm{h}_\infty^{2^s-2} \norm{h}_2^2.
\]
Note that~$\norm{h}_\infty^{2^s-2}
    =\norm{T^n \con f \cdot f}_\infty^{2^s-2}
        \leq \norm{f}^{2(2^s-2)}_\infty =: C_f$
is a positive constant independent of~$n$.
We conclude that
\[
    0 < C_f \cdot \p{\limsup_{N\to\infty} \frac{1}{N} \sum_{n=1}^N
        \int \abs{\E(T^n \con f \cdot f \mid \CZ_{s-1})}^2 \d\mu}.
\]
Since this holds for any bounded measurable function,
the extension is compact.
\end{proof}

\begin{lemma}\label{L:OneCompactLadder}
	For $s,\ell \in \N$, let  $(X,\mu,T_1,\ldots, T_\ell)$ be a system,  let $I\subseteq [\ell]$ be  non-empty,  let $Y\subseteq W$ be factors,   and 
	$$
	Z_I:= \bigvee_{j\in I} Z_s^j(W).
	$$
	Then  there exists a tower of factors
	$$
	Y_0:=Y\subseteq Y_1\subseteq \cdots \subseteq Y_N:=Y\vee Z_I
	$$
	such that for every  $n\in [N]$ the extension $Y_n\to Y_{n-1}$ is $T_i$-compact for 
	at least one  $i\in I$.
\end{lemma}
\begin{proof}
	After relabeling, we can assume that $I=[l]$ for some $l\in [\ell]$. 
	For $k\in [l]$, let
	$$
	Z_{[k]}:= Z_s^1(W)\vee \cdots \vee Z_s^k(W).
	$$
	We use finite induction on $k\in [l]$ to prove that the conclusion holds with $Z_{[k]}$ in place of $Z_I$. Since $Z_I=Z_{[l]}$ this implies the asserted statement. 
	
	For $k=1$, we  augment $Y$ by the Host-Kra tower 
	$$
	Y\subseteq Y\vee Z_0^1(W)\subseteq Y \vee  Z_1^1(W)\subset\cdots \subseteq Y\vee  Z_s^1(W).
	$$
	Each extension in this tower is $T_1$-compact  by parts~\eqref{iii} and \eqref{iv} of \cref{L:CompactProperties}. 
	
	Suppose  that $l\geq 2$, $k\in [l-1]$,  and we have proved the statement for $	Z_{[k]}$. We shall show that  it holds  for $Z_{[k+1]}$.
	To do this we expand the tower that reaches $Y_k':=Y\vee Z_{[k]}$ as follows
	$$
	Y_k'\subseteq Y_k'\vee  Z_0^{k+1}(W)\subseteq Y_k' \vee Z_1^{k+1}(W)\subseteq \cdots \subseteq Y_k'\vee Z_s^{k+1}(W)=Y\vee 	Z_{[k+1]}.
	$$
Arguing as in the $k=1$ case we get that 	each extension in this tower is $T_{k+1}$-compact.
\end{proof}

\subsection{Proof of the $PSZ$ property} We combine the auxiliary lemmas with the lifting result to establish the $PSZ$ property for every system.
\begin{theorem}
	For every  $d,\ell,r\in \N$,  every  system $(X,\mu,T_1,\ldots, T_\ell)$ has the $PSZ$ property in the sense of \cref{D:PSZ}.
\end{theorem}
\begin{proof}
Let  $d,\ell \in \N$ be fixed. 
	We will prove the following statement by finite induction on $k$: 
	\smallskip 
	
	\emph{ 	$PSZ(k)$: Let $(X,\mu,T_1,\ldots, T_\ell)$  be a system, let   $Y\subseteq W$ be factors, and   let $k\in \{0,\ldots, \ell\}$. Suppose that $Y$ 
		has the $PSZ$ property and  the extension $W\to Y$ is $T_i$-compact for at least $\ell-k$ values of $i\in [\ell]$. Then $W$ also has the $PSZ$ property.}
	
	\smallskip 

Note that for $k:=\ell$ and $Y:=\{pt\}$, $W:=X$, we get the $PSZ$ property for	every system, since the point system has the $PSZ$ property. To justify the last claim
	we need to prove that 
	$$
	\lim_{N\to \infty} \E_{n\in [N]} f(a_1^{q_1(n)} \cdots a_m^{q_m(n)} 1_{G/\Gamma}) > 0
	$$
	where $G/\Gamma$ is a nilmanifold, $f$ is a non-negative continuous function positive at $1_{G/\Gamma}$,
	and $q_1,\ldots,q_m$ is a jointly intersective family of polynomials.
	This claim follows immediately from the recurrence property of jointly
	intersective sequences \cite[Proposition 2.4]{bergelson2008intersective}.
	
\smallskip 
	{\bf Basis of the induction.} For $k=0$   we have that $W\to Y$ is $T_i$-compact for all $i\in [\ell]$,
so  the statement $PSZ(0)$  follows  from \cref{C:compactlift}.

    \smallskip 
    
{\bf Step of the induction.}	Suppose that $k\in [\ell]$ and  the statement $PSZ(k-1)$ holds. We shall show that  the statement $PSZ(k)$ holds. Let $h:=\ell-k$ and note that $h\in \{0,\ldots, \ell-1\}$.  After relabeling, we can assume that  $W\to Y$ is $T_i$-compact  for all $i\in [h]$,  and recall that  the factor $Y$ has the PSZ property. Our goal is to show that $W$ has the $PSZ$ property.
	It suffices to show that  $W$ has the $PSZ_r$ property for every $r\in \N$. So let $r\in \N$ be fixed
	and $s:=s(d,\ell,r)$ be as in \cref{cor:partial_factor_use}. 
	
	We apply \cref{L:OneCompactLadder} for the extension $W\to Y$ and the non-empty set
	$$
	I:=\{h+1,\ldots, \ell\}.
	$$
	We get that there exists a finite tower of factors
	$$
	Y_0:=Y\subseteq Y_1\subseteq \cdots \subseteq Y_N:=Y\vee Z_I, \quad 
	Z_I= \bigvee_{i\in I } Z_s^i(W)
	$$
	such that for every  $n\in [N]$ the extension $Y_n\to Y_{n-1}$ is $T_i$-compact for 
	at least one  $i\in I$. Note also  that   for all  $i\in [h]$ 
	the extension $Y_n\to Y_{n-1}$ is $T_i$-compact; indeed, 
	we know that 
 $W\to Y$ is $T_i$-compact  
  so the asserted compactness for the intermediate extensions $Y_n\to Y_{n-1}$  follows from 
 parts~\eqref{i} and \eqref{ii} of  \cref{L:CompactProperties}. We deduce that 
the extensions  $Y_n\to Y_{n-1}$  are $T_i$-compact for at least $h+1$ values of $i\in [\ell]$. Hence, the induction hypothesis applies to all these extensions and since the factor $Y$ has the $PSZ$ property, we deduce that  the factor $Y\vee Z_I$ has the $PSZ$ property. 
	
	 Finally, we obtain the $PSZ_r$ property for~$W$ by lifting
    the $PSZ$ property of~$Y\vee Z_I$. 
    We verify that the hypotheses of \cref{P:compactlift} are satisfied for the extension $W\to Y\vee Z_I$, i.e. we show that 
	\begin{equation}\label{E:YZJ}
	 Y\vee Z_I \vee   Z^i_s(W)\to   Y\vee Z_I
	\end{equation}
	is $T_i$-compact for all $i\in [\ell]$.
		For $i\in[h]$, this follows because the extension in
		\eqref{E:YZJ} is an intermediate extension of the
		$T_i$-compact extension $W\to Y$; hence compactness follows from
		parts~\eqref{i} and~\eqref{ii} of \cref{L:CompactProperties}.
		  For $i\in I=[\ell]\setminus [h]$, since $Z^i_s(W)\subseteq  Z_I$ the left and right hand sides of 
		  \eqref{E:YZJ} coincide, so we trivially get $T_i$-compactness.  
 Since $Y\vee Z_I$ has the $PSZ$ property, we deduce that $W$ has the $PSZ_r$ property, as required.  Since $r\in \N$ was arbitrary, this completes the proof. 
\end{proof}

\appendix

\section{Measure-theoretic lemmas} \label{app:measure}

Throughout this section~$(X,\CX,\mu)$ is a standard probability space, as in \cref{sec:prel_systems}.
For a probability measure~$\lambda$ on~$(X,\CX)$, a measurable function~$f$, and
sub-$\sigma$-algebras~$\CA$ and~$\CB$,
we write~$f\in_\lambda\CA$ if~$f$ coincides $\lambda$-almost everywhere
with an $\CA$-measurable function, and we use~$\subseteq _\lambda$ and~$=_\lambda$
accordingly. We denote by~$\CA\wedge_\lambda\CB$ the $\sigma$-algebra of sets~$E\in\CX$
with~$\1_E\in_\lambda\CA$ and~$\1_E\in_\lambda\CB$, so that
$L^2(\CA\wedge_\lambda\CB,\lambda)=L^2(\CA,\lambda)\cap L^2(\CB,\lambda)$.
The good intersection property of \cref{sec:good} with respect to~$\lambda$ is denoted
by~$\CA\perp^\lambda_{\CA\wedge_\lambda\CB}\CB$, meaning~$\CA$ and~$\CB$ are conditionally independent
given their meet.

Fix a Borel map~$\pi\colon X\to\Omega$ to a standard Borel space, let~$P:=\pi_*\mu$,
let~$\CI:=\pi^{-1}(\CB(\Omega))$,
and let~$(\mu_\omega)_{\omega\in\Omega}$ be the disintegration of~$\mu$
over~$\pi$. Since~$\mu_\omega(\pi^{-1}(\{\omega\}))=1$ for almost every~$\omega\in\Omega$,
we have
\begin{equation} \label{eq:trivial_components}
	\mu_\omega(\pi^{-1}(A))=\1_A(\omega)
	\text{ for every Borel } A\subseteq \Omega,
	\text{ for all } \omega \text{ outside a single nullset.}
\end{equation}
In \cref{sec:iso_Zs}, $\CI$ is (a version of) the invariant $\sigma$-algebra and~$(\mu_\omega)$ is the ergodic decomposition.

\noindent\emph{Countable generation.}
The results below concern the measures~$\mu_\omega$, which are sensitive to the choice of a
sub-$\sigma$-algebra within its class modulo~$\mu$. For instance, let~$X=[0,1]^2$ with Lebesgue measure,
$\pi(x,y)=x$, so that~$\mu_\omega=\delta_\omega\otimes\m_{[0,1]}$,
and let~$\CA$ be generated by~$\CI$ and the $\mu$-null Borel sets.
Then~$\CA=_\mu\CI$, but~$\CA=_{\mu_\omega}\CX$ for every~$\omega$, since every Borel set~$E$ is $\mu_\omega$-equal
to the nullset~$E\cap(\{\omega\}\times[0,1])$.
For this reason, the sub-$\sigma$-algebras to which we apply the results below are countably generated.
This costs nothing:
the preimage of the Borel $\sigma$-algebra under a factor map~$X\to Y$ is countably generated;
every sub-$\sigma$-algebra is $\mu$-equal to a countably generated one, since~$L^2(\mu)$ is separable;
and, replacing a countable family of generators~$(A_k)_{k\in \N}$
by~$(T^nA_k)_{n\in\Z,k\in\N}$,
the representatives can be taken $T$-invariant, and nested when the factors are nested.
In particular, in writing~$W_\omega$ for the $\omega$-component of a factor~$W$
in \cref{sec:iso_Zs}, we identify~$W$ with a countably generated invariant
sub-$\sigma$-algebra and~$W_\omega$ with said $\sigma$-algebra taken modulo~$\mu_\omega$.
By \cref{lem:disintegrating_inclusion}~$(ii)$ below, for~$\omega$ outside a nullset
this does not depend on the chosen representative.

\begin{lemma} \label{lem:disintegrating_inclusion}
Let~$\CA,\CB$ be sub-$\sigma$-algebras.
\begin{enumerate} [(i)]
\item If~$f\in_\mu\CA$, then~$f\in_{\mu_\omega}\CA$ for almost every~$\omega$.
\item If~$\CA$ is countably generated and~$\CA\subseteq _\mu\CB\vee\CI$, then~$\CA\subseteq _{\mu_\omega}\CB$
	for all~$\omega$ outside a single nullset.
\end{enumerate}
\end{lemma}
\begin{proof}
$(i)$ If~$f=a$ $\mu$-almost everywhere with~$a\in\CA$, then~$\int_\Omega\mu_\omega(\{f\neq a\})\d P=\mu(\{f\neq a\})=0$.
$(ii)$ Apply~$(i)$ to the indicators of a countable algebra generating~$\CA$.
For~$\omega$ outside a nullset, they all belong to~$\CB\vee\CI$ modulo~$\mu_\omega$,
and every set of~$\CI$ is $\mu_\omega$-trivial by~\eqref{eq:trivial_components}.
Since the sets~$E$ with~$\1_E\in_{\mu_\omega}\CB$ form a $\sigma$-algebra, it contains~$\CA$.
\end{proof}

\begin{lemma} \label{lem:relative_conditioning}
Let~$\CA$ be countably generated with~$\CI\subseteq _\mu\CA$, and let~$f\in L^2(\mu)$.
Then~$\E_\mu(f\mid\CA)=_{\mu_\omega}\E_{\mu_\omega}(f\mid\CA)$ for almost every~$\omega$.
\end{lemma}
\begin{proof}
Let~$g$ be an $\CA$-measurable version of~$\E_\mu(f\mid\CA)$ and~$\CG$ a countable algebra generating~$\CA$.
Then~$g=_{\mu_\omega}\E_{\mu_\omega}(f\mid\CA)$ holds
if and only if~$\varphi_A(\omega):=\int_A(f-g)\d\mu_\omega$
vanishes for every~$A\in\CG$, so it suffices to show~$\varphi_A=0$ almost everywhere for a fixed~$A \in \CG$.
For a Borel set~$G\subseteq \Omega$, let~$I:=\pi^{-1}(G)$. Then
\[
	\int_G\varphi_A\d P=\int_\Omega\int\1_{I\cap A}\cdot (f-g)\d\mu_\omega\d P
		=\int\1_{I\cap A}\cdot (f-g)\d\mu=0,
\]
where the last equality holds because~$I\cap A\in_\mu\CA$, as~$\CI\subseteq _\mu\CA$.
\end{proof}

\begin{corollary} \label{lem:good_inclusion}
Let~$\CA$ be countably generated with~$\CI\subseteq _\mu\CA$, and let~$f\in L^2(\mu)$.
Then~$f\in_\mu\CA$ if and only if~$f\in_{\mu_\omega}\CA$ for almost every~$\omega$.
\end{corollary}
\begin{proof}
One direction is \cref{lem:disintegrating_inclusion}~$(i)$. Conversely, by \cref{lem:relative_conditioning} we
have the identity~$f=_{\mu_\omega}\E_{\mu_\omega}(f\mid\CA)
	=_{\mu_\omega}\E_\mu(f\mid\CA)$ for almost every~$\omega$,
so~$\norm{f-\E_\mu(f\mid\CA)}^2_{L^2(\mu)}=\int_\Omega\norm{f-\E_\mu(f\mid\CA)}^2_{L^2(\mu_\omega)}\d P=0$.
\end{proof}

\begin{corollary} \label{lem:good_transfer}
Let~$\CA,\CB$ be countably generated with~$\CI\subseteq _\mu\CA$ and~$\CI\subseteq _\mu\CB$.
If~$\CA$ and~$\CB$ have the good intersection property with respect to~$\mu_\omega$
for almost every~$\omega$, then they have it with respect to~$\mu$.
\end{corollary}
\begin{proof}
Let~$a$ be bounded and $\CA$-measurable. By \cref{lem:relative_conditioning} and the hypothesis,
$\E_\mu(a\mid\CB)=_{\mu_\omega}\E_{\mu_\omega}(a\mid\CB)\in_{\mu_\omega}\CA$ for almost every~$\omega$.
Hence~$\E_\mu(a\mid\CB)\in_\mu\CA$ by \cref{lem:good_inclusion}.
\end{proof}

The next result compares the intersection of factors, restricted to~$\mu_\omega$,
with the intersection computed with respect to~$\mu_\omega$.

\begin{corollary} \label{lem:good_wedge_limit}
Let~$\CA,\CB$ and~$\CB_1\subseteq \CB_2\subseteq \cdots\subseteq \CB$ be countably generated,
all containing~$\CI$ modulo~$\mu$, with~$\CB=_\mu \bigvee_{n\in \N} \CB_n$.
Suppose that~$\CA$ and~$\CB$ have the good intersection property with respect to~$\mu$,
and that for almost every~$\omega$ there is~$n_0(\omega)$ such that
$\CA \perp^{\mu_\omega}_{\CA\wedge_{\mu_\omega}\CB_n} \CB_n$ for every~$n\geq n_0(\omega)$.
Then
\[
	\CA\wedge_\mu \CB =_{\mu_\omega} \bigvee_{n\in \N} (\CA\wedge_{\mu_\omega} \CB_n)
\]
for almost every~$\omega$, where~$\CA\wedge_\mu\CB$ denotes a countably generated representative.
\end{corollary}
\begin{proof}
All statements below involve countably many functions and sub-$\sigma$-algebras,
so we may fix~$\omega$ outside a single nullset.
By \cref{lem:disintegrating_inclusion}$(ii)$, we have~$\CB=_{\mu_\omega}\bigvee_{n\in\N}\CB_n$
and~$\CA\wedge_\mu\CB\subseteq _{\mu_\omega}\CA\wedge_{\mu_\omega}\CB$ for almost every~$\omega$.
Conversely, let~$A\in\CA$. By the good intersection property and \cref{lem:relative_conditioning},
$\E_{\mu_\omega}(\1_A\mid\CB)=_{\mu_\omega}\E_\mu(\1_A\mid\CB)\in_\mu\CA\wedge_\mu\CB$.
Letting~$A$ range over a countable algebra generating~$\CA$, these functions span~$L^2(\CA\wedge_{\mu_\omega}\CB,\mu_\omega)$.
Hence~$\CA\wedge_\mu\CB=_{\mu_\omega}\CA\wedge_{\mu_\omega}\CB$ for almost every~$\omega$.
We are left to prove the identity
\[
	\CA\wedge_{\lambda} \CB =_{\lambda} \bigvee_{n\in \N} (\CA\wedge_{\lambda} \CB_n)
\]
for any probability measure~$\lambda$ provided that~$\CB=_{\lambda} \bigvee_{n\in \N} \CB_n$
and~$\CA \perp^{\lambda}_{\CA\wedge_{\lambda}\CB_n} \CB_n$ for all large~$n$.
The inclusion~$\supseteq_\lambda$ is trivial.
Conversely, let~$f\in L^2(\CA\wedge_{\lambda} \CB,\lambda)$.
By the martingale convergence theorem, $f_n:=\E_\lambda(f\mid \CB_n)$ converges to~$f$.
On the other hand, as~$f\in_\lambda \CA$, the good intersection property gives~$f_n \in_\lambda \CA \wedge_{\lambda}\CB_n$
for all large~$n$, concluding the proof.
\end{proof}


	
	\appendix


	\bibliographystyle{abbrv}
	\bibliography{refs}

\end{document}